\documentclass[11pt,twoside]{article}

\usepackage{amssymb,amsmath,amsfonts}
\usepackage{mathptmx}
\usepackage{euler,palatino}
\usepackage{latexsym, amssymb, amsmath, theorem, subfigure}

\theoremstyle{plain} \theorembodyfont{\itshape}
\newtheorem{theorem}{Theorem}[section]
\newtheorem{proposition}[theorem]{Proposition}
\newtheorem{lemma}[theorem]{Lemma}
\newtheorem{corollary}[theorem]{Corollary}

\theorembodyfont{\rmfamily}
\newtheorem{definition}[theorem]{Definition}

\newtheorem{remark}[theorem]{Remark}
\newtheorem{conjecture}[theorem]{Conjecture}

\numberwithin{equation}{section}

\usepackage{amscd,amsfonts}
\usepackage{pst-all}
\usepackage{amsfonts}
\usepackage[dvips]{graphicx}
\usepackage{epic}
\usepackage{verbatim}
\usepackage{amssymb, upref}
\usepackage[dvips]{graphicx}
\usepackage{epsf}

\def\a{{\alpha}}
\def\b{{\beta}}
\def\c{{\gamma}}

\def\lam{{\lambda}}
 \def\eps{{\epsilon}}
 \def\G{{\Gamma}}
\def\C{{\mathbb C}}

\def\N{{\mathbb N}}
\def\S{{\mathbb S}}
\def\CP{{\mathbb C\mathbb P}}
\def\del{{\delta}}
\def\ov{\overline}

\title{\textbf{The moduli space of  germs of generic families of analytic diffeomorphisms unfolding  a
parabolic fixed point}\thanks{This work is supported by NSERC in
Canada.}}
\author{
C. Christopher, School of Mathematics and Statistics, University of Plymouth\\
C. Rousseau, DMS and CRM, Universit\'e de Montr\'eal}
\date{July 2008}

\begin{document}
\maketitle

\begin{abstract}
In this paper we describe the moduli space of germs of generic families of
analytic diffeomorphisms which unfold a parabolic fixed point of
codimension 1.

In \cite{MRR} (and also \cite{R}), it was shown that the
Ecalle-Voronin modulus can be unfolded to give a complete modulus
for such germs. The modulus is defined on a ramified sector in the
canonical perturbation parameter $\eps$. As in the case of the
Ecalle-Voronin modulus, the modulus is defined up to a linear
scaling depending only on $\eps$.

Here, we characterize the moduli space for such unfoldings by finding
the compatibility conditions on the modulus which are necessary and
sufficient for realization as the modulus of an unfolding.

The compatibility condition is obtained by considering the region of
sectorial overlap in $\eps$-space.  This lies in the Glutsyuk sector
where the two fixed points are hyperbolic and connected by the orbits of the
diffeomorphism.  In this region we have two
representatives of the modulus which describe the same dynamics. We
identify the necessary compatibility condition between these two
representatives by comparing them both with their common Glutsyuk
modulus.

The compatibility condition implies the existence of a linear
scaling for which the modulus is $1/2$-summable in $\eps$, whose
direction of non-summability coincides with the direction of real
multipliers at the fixed points.  Conversely, we show that the
compatibility condition (which implies the summability property) is
sufficient to realize the modulus as coming from an analytic
unfolding, thus giving a complete description of the space of
moduli. The terminology \lq\lq space" of moduli is justified by the
fact that the moduli depend analytically on extra parameters.

\end{abstract}


\section{Introduction}

The analytic classification of  germs of analytic diffeomorphisms
with a parabolic fixed point of codimension 1 was given by Ecalle
\cite{E} and Voronin \cite{V}. A complete modulus for a germ of
diffeomorphism $f:(\C,0)\rightarrow (\C,0)$ of the form $$f(z)= z+
z^2+o(z^2)$$ is given by a formal invariant $a\in\C$ and an
equivalence class of a pair of germs $(\psi^0,\psi^\infty)$ where
$\psi^0:(\C,0)\rightarrow (\C,0)$,
$\psi^\infty:(\C,\infty)\rightarrow (\C,\infty)$, and where the
equivalence relation is defined as follows:
$$(\psi^0,\psi^\infty)\sim (\breve{\psi}^0,\breve{\psi}^\infty)
\Longleftrightarrow \exists C,C'\in \C^*
\begin{cases} \breve{\psi}^0= L_C\circ \psi^0\circ L_{C'}\\
\breve{\psi}^\infty= L_C\circ \psi^\infty\circ L_{C'}\end{cases}$$
where $L_C$ (resp. $L_{C'}$) is the linear map $w\mapsto Cw$ (resp.
$w\mapsto C'w$). Moreover all tuples
$(a,[\psi^0,\psi^\infty])$ are realizable, where $[\psi^0,\psi^\infty]$
represents the equivalence class of $(\psi^0,\psi^\infty)$.

The paper \cite{MRR} addresses the similar question for the analytic
classification of generic 1-parameter families of analytic
diffeomorphisms unfolding a parabolic fixed point. It was
shown that it is possible to prepare the family so that the
parameter becomes an analytic invariant. Then a conjugacy between
two germs of prepared families must preserve the canonical
parameter. The main result of \cite{MRR} is that the unfolding of
$(a,[\psi^0,\psi^\infty])$ is a complete modulus of analytic
classification for a prepared germ $f_\eps:(\C,0)\rightarrow (\C,0)$
of the form
$$f_\eps(z)=z+(z^2-\eps)(1+b(\eps)+c(\eps)z+O(z^2-\eps)),$$
such that $\frac{\partial f_\eps}{\partial \eps}\neq0$ and $f_0$ has
formal invariant $a$. The paper \cite{MRR} also allows an
explanation of the meaning of the coefficients which form the
Ecalle-Voronin modulus. Indeed the formal invariant $a$ indicates a
shift between the multipliers of the two fixed points in the limit
$\eps=0$. To interpret the coefficients of $\psi^{0,\infty}$ it is
better to split the parameter space in two regions: in the Glutsyuk
region where the two fixed points are hyperbolic and there is an
orbit connecting them, then the coefficients of the unfolded
$\psi_\eps^{0,\infty}$ measure the non compatibility of the two
\lq\lq models" at the fixed points. In the Lavaurs region, they
control the complicated dynamics of the fixed points. In particular
the ``parametric resurgence" phenomenon allows one to predict from
the coefficients of $\psi^{0}$ (resp. $\psi^\infty$) some discrete
sequences $\{\eps_n\}$ converging to the origin  for which the fixed
point $-\sqrt{\eps_n}$ (resp. $\sqrt{\eps_n}$) of $f_{\eps_n}$ is
resonant and nonlinearizable. Moreover it was shown in \cite{MRR}
that it is possible to take a representative of the equivalence
class $[\psi_{\hat{\eps}}^0,\psi_{\hat{\eps}}^\infty]$ depending
analytically on $\hat{\eps}$, for $\hat{\eps}$ in a sector $V$ of
opening less than $4\pi$ of the universal covering of $\eps$ space
punctured at $0$.

While it is easily shown that a function $a(\eps)$ is realizable as
the formal modulus of the family if and only if it is analytic, the
other part of the necessary and sufficient conditions for
realizability of a modulus and the determination of the moduli
space was completely open. The difficulty comes from the fact that
the construction leading to the modulus
$[\psi_{\hat{\eps}}^0,\psi_{\hat{\eps}}^\infty]$ of a family cannot
be extended to make a full turn in $\sqrt{\eps}$. This is because
the unfolded $\psi_{\hat{\eps}}^0$ (resp.
$\psi_{\hat{\eps}}^\infty$) is attached to $-\sqrt{\hat{\eps}}$
(resp. $\sqrt{\hat{\eps}}$), which gives two completely different
descriptions of the same dynamics of $f_\eps$  when $\hat{\eps}$
makes a full turn.  This fact is precisely what we need to exploit to
identify the sufficient condition for realizability. Indeed, in the
Glutsyuk region, i.e. the region where the fixed points are
hyperbolic,  the renormalized return map near $-\sqrt{\hat{\eps}}$
(resp. $\sqrt{\hat{\eps}}$) is the composition of
$\psi_{\hat{\eps}}^0$ (resp. $\psi_{\hat{\eps}}^\infty$) with a
linear map. Since the fixed points are hyperbolic, these
renormalized return maps are linearizable. The comparison of the
linearizing maps is an analytic invariant, thus allowing one to
derive a compatibility condition between
$(\psi_{\hat{\eps}}^0,\psi_{\hat{\eps}}^\infty)$ and
$(\psi_{\hat{\eps}e^{2 \pi i}}^0,\psi_{\hat{\eps}e^{2\pi
i}}^\infty)$ so that they describe the same dynamics.

One important consequence of the compatibility condition is that it
is possible to choose a representative of the equivalence class
$[\psi_{\hat{\eps}}^0,\psi_{\hat{\eps}}^\infty]$ such that
$\psi_{\hat{\eps}}^0$ and $\psi_{\hat{\eps}}^\infty$ are both
$1/2$-summa\-ble in $\eps$, with $\mathbb R^+$ as direction of
non-summability. This property, together with the compatibility
condition, is sufficient for a germ of family
$(a(\eps),[\psi_{\hat{\eps}}^0,\psi_{\hat{\eps}}^\infty])$ to be
realizable.

The realization is done in two steps. We first realize locally by a
family $f_{\hat{\eps}}$ ramified in $\hat{\eps}$. We do this by
first giving the realization for a fixed $\hat{\eps}$: we construct
the realization on an abstract manifold and use the Ahlfors-Bers
theorem to show that this manifold is indeed an open set of $\C$. We
then show that the construction can be performed so as to depend
analytically on $\hat{\eps}$.  We call this part the local
realization. The second step is to correct the ramification. Indeed
using the local realization and the compatibility condition, this
allows us to construct a realization on an abstract 2-dimensional
manifold. The Newlander-Nirenberg theorem can be applied to show
that this manifold is indeed an open set of $\C^2$ containing a
product of a neighborhood of the origin in $\eps$-space with an open
set of $\C$.

The compatibility condition puts very strong constraints on the
families $(a(\eps),[\psi_{\hat{\eps}}^0,\psi_{\hat{\eps}}^\infty])$
that can be realized. Indeed, we have already mentioned that this
forces the existence of a  representative
$\psi_{\hat{\eps}}^{0,\infty}$ which is $1/2$-summable in $\eps$.
But this is far from being sufficient. For instance, we analyze in
detail the case of the Riccati equation and prove that the
compatibility condition
  implies in that case that there exists representatives of the modulus
  $\psi_{\hat{\eps}}^{0,\infty}$ which are
analytic in $\eps$. This allows us to completely characterize the
modulus space in this special case. We also exhibit an example of
family $(a(\eps),[\psi_{\eps}^0,\psi_{\eps}^\infty])$ depending
analytically on $\eps$ which cannot be realized as a modulus.

The identification of the moduli space opens great possibilities.
Indeed, while the knowledge of the Ecalle-Voronin modulus of $f_0$
allows one to deduce the nonlinearizability of the fixed points of
$f_\eps$ when special kinds of resonance occurs (the ``parametric
resurgence" phenomenon mentioned earlier), the dependence in $\eps$
is crucial to be able to draw similar conclusions in the case of
fixed points whose multipliers are irrational rotations, or, in the
case of resonance, when we consider the more complex question of the
convergence of the change of coordinate to normal form. For
instance, it is known that the quadratic map $f_\eps(z)=z(1+\eps)
+z^2$ is never linearizable when $1+\eps=e^{2\pi i \alpha}$ with
$\alpha\notin \mathbb Q$ not a Brjuno number. The system is also
never normalizable when $1+\eps$ is a root of unity. But what can be
said of a map $g_\eps(z)= f_\eps(z)+h_\eps(z)$ with
$h_\eps(z)=o(z^2)$? We hope that our results will give tools to
answer such questions.

Another potential application is in the spirit of Hilbert's 16th
problem. This problem deals with the maximum number $H(n)$ and
relative positions of limit cycles of polynomial vector fields of
degree $\leq n$. The finiteness subproblem deals with the existence
of a uniform upper bound for the number of limit cycles of
polynomial vector fields of degree at most $n$ for each integer $n$:
$H(n)<\infty$. In the paper \cite{DRR} it is shown how the
finiteness part for $n=2$ can be reduced to 121 local problems,
namely showing that 121 graphics have finite cyclicity: let us call
this the DRR-program. A significant step in the DRR-program was
performed in the paper \cite{DIR} where it is shown how the use of
the Martinet-Ramis invariant of a saddle-node allows to prove the
finite cyclicity of several generic graphics of this program. The
most difficult graphics of the DRR-program are graphics surrounding
centers. An efficient method to prove their cyclicity is to divide
the displacement map in the Bautin ideal. This method requires a
deep understanding and a fine control of the dependence on the
parameters.  The compatibility condition is a natural candidate for
obtaining further results in this direction.

The paper is organized as follows. In
Section~\ref{sect:preliminaries} we recall the definition of the
Ecalle-Voronin modulus, the preparation of the family  and  the
results of \cite{MRR}. In Section~\ref{sect:local} we prove the
local realization theorem. In Section~\ref{sect:compatibility} we
derive the compatibility condition and we prove the
$1/2$-summability of $\psi_{\hat{\eps}}^{0,\infty}$  in $\eps$. In
Section~\ref{sec:global} we prove the global realization theorem.
Finally in Section~\ref{sect:Riccati} we study  examples including
 the unfolding of a Riccati equation with a saddle-node, and give
a complete analytic classification of its local unfoldings.

\section{Preliminaries}\label{sect:preliminaries}

\subsection{Notations }\label{notations}

The notations collected here are often referred to in the paper.
\begin{itemize}
\item $L_C$:  the linear map
\begin{equation}L_C(w)=Cw;\label{def_L}\end{equation}
\item $m_A$: the M\"obius transformation \begin{equation}m_A(w)=
\frac{w}{1+Aw};\label{def_Mobius}\end{equation}
\item $T_B$: the translation
\begin{equation}T_B(W)=W+B;\label{def_T}\end{equation}
\item $E$: the map \begin{equation}E(W)=\exp(-2\pi i W)\label{def_E}\end{equation} with inverse
$E^{-1}(w)=-\frac1{2\pi i}\ln (w)$;
\item $R^0$ and $R^\infty$ are the domains of $\mathbb C$ defined
by:
\begin{equation}
\begin{cases} R^0=\{W|Im W< -Y_0\},\\
R^\infty=\{W|Im W>Y_0\},
\end{cases}\label{domain_R}\end{equation}
where $Y_0$ is some sufficiently large constant.
\item We will be dealing with fixed points $\pm\sqrt{\eps}$ of a
diffeomorphism $f_\eps$. In order to make this well-defined, we work
on the universal covering of $\eps$-space punctured at $0$
parameterized by $\hat{\eps}$. The function $\sqrt{\hat{\eps}}$ is
defined by $\arg\sqrt{\hat{\eps}}=\frac{\arg\hat{\eps}}{2}$. In
particular $\sqrt{\hat{\eps}}\in \mathbb R^+$, when $\arg
\hat{\eps}=0$. \item Upper indices $0$ and $\infty$ will be
associated to the two parts of the modulus and to other objects. In
all cases, $0$ (resp. $\infty$) will
 be associated with $-\sqrt{\hat{\eps}}$ (resp. $\sqrt{\hat{\eps}}$).
\end{itemize}

\subsection{The Ecalle-Voronin modulus of a diffeomorphism and its unfolding}

We briefly summarize some results of \cite{MRR} on the unfoldings of
the Ecalle-Voronin invariants of a generic parabolic point of a
diffeomorphism
\begin{equation} f(z)= z + z^2 +o(z^2).\label{mod.1} \end{equation}  Since the paper \cite{MRR} only deals with
$1$-parameter families, we start by proving a \lq\lq preparation
theorem" for generic unfoldings with several parameters. The
preparation makes clear the role of the  \lq\lq canonical parameter".

The perspective of \cite{MRR} is to compare  a generic $1$-parameter
family $f_\eps$ with a \lq\lq model" family, namely the time-one map
for the family of vector fields
\begin{equation}v_\eps(z)=\frac{z^2-\eps}{1+a(\eps)z}\frac{\partial}{\partial
z}.\label{mod.3}\end{equation}  If $\mu_{\hat{\eps}}^0$ and
$\mu_{\hat{\eps}}^\infty$ are the eigenvalues at the singular points
$-\sqrt{\hat{\eps}}$ and $\sqrt{\hat{\eps}}$ of \eqref{mod.3}, then
we can remark that
\begin{eqnarray}
a(\eps)
=\frac1{\mu_{\hat{\eps}}^\infty}+\frac1{\mu_{\hat{\eps}}^0},\qquad\quad
\frac1{\sqrt{\hat{\eps}}}=\frac1{\mu_{\hat{\eps}}^\infty}-\frac1{\mu_{\hat{\eps}}^0},
\label{mod.4}\end{eqnarray} i.e. $\eps$ and $a(\eps)$ are analytic
invariants of the system \eqref{mod.3}. Moreover $a(\eps)$ depends
analytically on $\eps$. We wish to prepare our family of
diffeomorphisms so that the multipliers at the fixed points,
$\lambda^{0,\infty}_{\hat{\eps}}$, correspond to those of the
time-$1$ map of \eqref{mod.3}, and hence
$\lambda^{0,\infty}_{\hat{\eps}}=\exp(\mu^{0,\infty}_{\hat{\eps}})$.

\begin{theorem}\label{thm_prepared}
 We consider a germ of a $k$-parameter  analytic family of
diffeomorphisms $f_\eta:(\C,0)\rightarrow (\C,0)$ depending on a
multi-parameter $\eta=(\eta_1, \dots, \eta_k)$ with a double fixed
point at the origin for $\eta=0$, such that $\frac{\partial
f}{\partial \eta_j}\neq0$ for some $j\in\{1, \dots, k\}$. There
exists a germ of analytic change of coordinates and parameters
$(z,\eta)\mapsto (\ov{z}, \eps, \nu_1, \dots, \nu_{k-1})$
transforming the family to the prepared form
\begin{equation}\ov{f}_{\eps,\nu}(\ov{z}) = \ov{z}
+(\ov{z}^2-\eps)\left[1+\ov{\b}(\eps,\nu)+
\ov{A}(\eps,\nu)\ov{z}+(\ov{z}^2-\eps)\ov{Q}(\ov{z},\eps,\nu)\right],\label{prepared_k_par}\end{equation}
with the additional property that $\ov{\b}(0,0)=0$ and
$$\frac1{\sqrt{\eps}}= \frac1{\ln(\ov{f}_{\eps,\nu}'(\sqrt{\eps}))}
-\frac1{\ln(\ov{f}_{\eps,\nu}'(-\sqrt{\eps}))}.$$  The parameter
$\eps$ is unique and called the {\it canonical
parameter}.  With this choice of canonical parameter, the function
$$a(\eps)= \frac1{\ln(\ov{f}_{\eps,\nu}'(\sqrt{\eps}))}
+\frac1{\ln(\ov{f}_{\eps,\nu}'(-\sqrt{\eps}))},$$
is a formal invariant of the system which depends analytically on $\eps$.

\end{theorem} \noindent{\scshape Proof.} Since
$\frac{\partial^2 (f_\eta-id)}{\partial z^2}(0)\neq0$, the
Weierstrass preparation theorem allows to write $f_\eta(z)- z =
P_\eta(z)U_\eta(z)$ where $P_\eta(z)$ is a Weierstrass polynomial of
degree 2 and $U_\eta(z)\neq 0$ for small $(z,\eta)$. A translation
in $z$ allows to bring $P_\eta(z)$ to the form (we do not change the
name of the variable) $z^2- D(\eta)$. Moreover the genericity
implies that $\frac{\partial D(\eta)}{\partial \eta_j}\neq0$,
allowing to replace the parameter $\eta_j$ by  $\ov{\eps}=D(\eta)$.
Let $\nu= (\eta_1, \dots, \eta_{j-1},\eta_{j+1}, \dots, \eta_k)$
Using a dilatation in $z$ and $ \ov{\eps}$ (without changing their
names) we can suppose that the initial family has the two fixed
points located at $z=\pm\sqrt{\ov{\eps}}$ and that $U(0,0)=1$, i.e.
that we start with a family:
 $$f_{\eta}(z) = z + (z^2-\ov{\eps}) h(z,\ov{\eps},\nu),$$
 where
 $h(z,\ov{\eps},\nu )= 1+O(z)+O(|\ov{\eps},\nu|)$.
 By the Weierstrass division theorem we have
$$h(z,\ov{\eps},\nu )=
Q(z,\ov{\eps},\nu)(z^2-\ov{\eps})+(a_0+B(\ov{\eps},\nu))z+1+C(\ov{\eps},\nu),$$
 where $B(0,0)=C(0,0)=0$.
 The multipliers at $\pm\sqrt{\ov{\eps}}$  are given by
\begin{equation}\begin{cases}
\lam^0=f_\eta'(-\sqrt{\ov{\eps}})=
1-2\sqrt{\ov{\eps}}\left[1+C(\ov{\eps},\nu)-(a_0+B(\ov{\eps},\nu))\sqrt{\ov{\eps}}\right],\\
 \lam^\infty=f_\eta'(\sqrt{\ov{\eps}})=1+2\sqrt{\ov{\eps}}\left[1+C(\ov{\eps},\nu)+(a_0+B(\ov{\eps},\nu))
 \sqrt{\ov{\eps}}\right].
 \end{cases}\end{equation}

An additional scaling in $z$ and $\ov{\eps}$ is necessary of the
form $$(\ov{z},\eps)= (z(1+b(\ov{\eps},\nu)),
\ov{\eps}(1+b(\ov{\eps},\nu))^2),$$ with
$b(\ov{\eps},\ov{\eta})=O(|\ov{\eps},\nu|)$
 to be determined. It changes the family to the form
   $$\ov{f}_{\eps,\nu}(\ov{z}) =\ov{z} + (\ov{z}^2-\eps^2)
   \left(\frac{1+C(\eps,\nu)}{1+b(\ov{\eps},\nu)}+ \frac{a_0+B(\eps,\nu)}{(1+b(\ov{\eps},\nu))^2}\ov{z} +
   (\ov{z}^2-\eps)\ov{Q}(\ov{z},\eps,\nu)\right).
   $$
We ask that the new multipliers at $\pm\sqrt{\eps}$ satisfy
 $$ \frac1{\ln(\lam^\infty)}-\frac1{\ln(\lam^0)}=\frac1{\sqrt{\eps}}.$$ This equation is
 solvable since
 $$\left.\frac{\partial \sqrt{\eps}\left(\frac1{\ln(\lam^\infty)}-\frac1{\ln(\lam^0)}\right)}{\partial b}\right|_{\eps=0}\neq0$$
 and yields $b(\ov{\eps},\nu)=O(|\ov{\eps},\nu|)$.
The other formal invariant is given by
\begin{equation}a(\ov{\eps},\nu)=
\frac1{\ln(\lambda^0)}+\frac1{\ln(\lambda^\infty)},\label{equation_a}\end{equation}
which is clearly analytic in $\eps$ and $\nu$.  Thus, we obtain the
required form,
\begin{equation}\begin{cases}\lambda^0=\exp\left(-\frac{2\sqrt{\ov{\eps}}}{1-
\ov{a}(\ov{\eps},\nu)\sqrt{\ov{\eps}}}\right),\\
\lam^\infty=  \exp\left(\frac{2\sqrt{\ov{\eps}}}{1+
\ov{a}(\ov{\eps},\nu)\sqrt{\ov{\eps}}}\right).\end{cases}\label{cond_prepared}\end{equation}
\hfill$\Box$

\medskip

The paper \cite{MRR} describes a complete modulus of analytic
classification for one-parameter prepared families of the form
\eqref{prepared_k_par} for values of $\eps$ in a small neighborhood
of the origin. This modulus is given by an unfolding of the
Ecalle-Voronin modulus of $f_0$. Since $\eps$ is an analytic
invariant for a prepared family, it is given by a family of moduli
for each fixed value of $\eps$. However no family of moduli analytic
in $\eps$ exists in general, so the modulus must be defined in a
ramified way. Furthermore \cite{MRR} does not identify a sufficient
condition for such a family to be realizable as the modulus of an
unfolding.

\medskip

\begin{figure}
\begin{center}
\includegraphics[height=5cm]{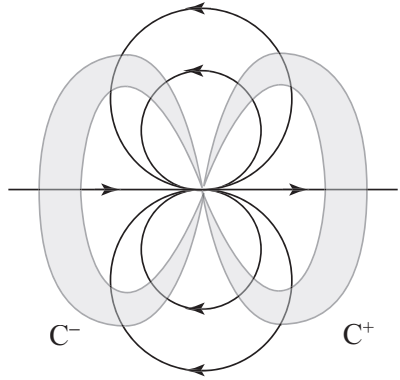}
\caption{The Ecalle-Voronin modulus } \label{modulus0}
\end{center}
\end{figure}

\noindent{\bf Description of the Ecalle-Voronin modulus ($\eps=0$).}
This modulus is effectively given by the {\it orbit space} of $f_\eps$.
We consider two fundamental domains $C^\pm$ of crescent shapes as in
Figure~\ref{modulus0}, which are given by two curves $l_{\pm}$ and
their images by $f_0$.

Each  orbit is represented by at most one point in each crescent,
but some orbits can have representatives in the two crescents. Hence
the orbit space is the union of the two crescents modulo the
identification of points of the same orbit. To give this
identification in an intrinsic way, one remarks that the two
crescents in which we identify the curves $l_\pm$ and $f(l_\pm)$
have the conformal structure of spheres $\S^\pm$, with the points 0
and $\infty$ identified. The coordinates on the spheres are unique
up to linear changes of coordinates.  Then the Ecalle-Voronin
modulus is the equivalence class of  pairs of germs
$(\psi^0,\psi^{\infty})$ of analytic diffeomorphisms, where
$\psi^0:(\S^{+},0)\rightarrow (\S^{-},0)$ and
$\psi^{\infty}:(\S^{+},\infty)\rightarrow (\S^{-},\infty)$ are
defined respectively in the neighborhoods of 0 and $\infty$, under
conjugation by linear changes of coordinates in the source and
target space. Let us define a map $f_0$ to be \emph{iterable} or
\emph{embedable} if $f_0$ is the time-one map of an analytic vector
field.  The map $f_0$ is iterable if and only if both of the germs
$\psi^0$ and $\psi^{\infty}$ are linear.

\medskip

\noindent{\bf The unfolded Ecalle-Voronin modulus.}\label{unfoldEV}
In \cite{MRR} it is proved that for any sufficiently small
neighborhood $U$ of the origin in $z$-space and for any
$\delta\in(0,\pi)$ (later we will restrict to
$\delta\in(0,\frac{\pi}2)$), there exists $\rho>0$, which is the
radius of a small sectorial neighborhood
\begin{equation}V_{\rho,\delta}=\{\hat{\eps}:|\hat{\eps}|<\rho, \arg(\hat{\eps})\in
(-\delta,2\pi+\delta)\}\cup\{\eps=0\},\label{sector_V}\end{equation}
of the origin in the universal covering of the parameter space
punctured at $0$ such that for each $\hat{\eps}\in V_{\rho, \delta}$
the orbit space is described as follows
\begin{itemize}
\item There exists two crescents $C^\pm_{\hat{\eps}}$ with endpoints at the
two singular points bounded by curves $l_{\pm,\hat{\eps}}$ and their
images $f_\eps(l_{\pm,\hat{\eps}})$ (Figure~\ref{modulus}).

\begin{figure}[!h]
\begin{center}
\includegraphics[height=10cm]{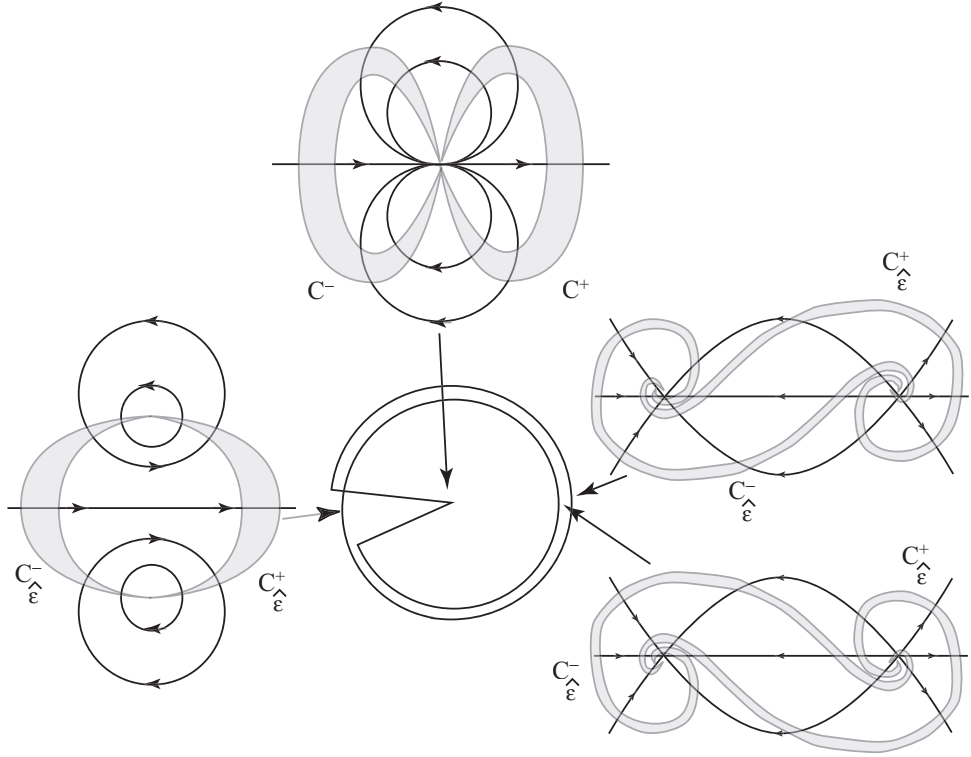}
\caption{The  modulus for the family } \label{modulus}
\end{center}
\end{figure}

\item The crescents
$C^\pm_{\hat{\eps}}$ in which we identify the curves
$l_{\pm,\hat{\eps}}$ and their images $f_\eps(l_{\pm,\hat{\eps}})$
have the conformal structure of spheres $\S_{\hat{\eps}}^{\pm}$ with
the singular point $\sqrt{\hat{\eps}}$ (resp. $-\sqrt{\hat{\eps}}$)
located at $\infty$ (resp. 0).
\item Points in the two neighborhoods of 0 and $\infty$ on
the spheres $\S_{\hat{\eps}}^{\pm}$ are identified modulo analytic
maps, $\psi_{\hat{\eps}}^0,\psi_{\hat{\eps}}^{\infty}:
\S_{\hat{\eps}}^+\rightarrow \S_{\hat{\eps}}^-$, defined in the
neighborhoods of 0 and $\infty$ respectively. These maps are
obviously uniquely defined up to the choice of coordinates on the
spheres. Hence it is natural to consider the equivalence classes of
pairs $(\psi_{\hat{\eps}}^0,\psi_{\hat{\eps}}^{\infty})$ under the
equivalence relation:
\begin{equation}(\psi_{\hat{\eps}}^0,\psi_{\hat{\eps}}^{\infty})\sim(\breve{\psi}_{\hat{\eps}}^0,
\breve{\psi}_{\hat{\eps}}^{\infty})\Longleftrightarrow \exists
c(\hat{\eps}),c'(\hat{\eps})\in\mathbb C ^*\;
\begin{cases} \breve{\psi}_{\hat{\eps}}^0(w)=c'(\hat{\eps})\psi_{\hat{\eps}}^0(c(\hat{\eps}) w) \\
 \breve{\psi}_{\hat{\eps}}^{\infty}(w)=c'(\hat{\eps})\psi_{\hat{\eps}}^{\infty}(c(\hat{\eps})
 w)\end{cases}\label{mod.6}\end{equation}
where $c(\hat{\eps}), c'(\hat{\eps})$ are analytic in $V_{\rho,\delta}\setminus\{0\}$
 with continuous non-zero limit at $0$.  Let us denote the equivalence class of the family
 $(\psi_{\hat{\eps}}^0,\psi_{\hat{\eps}}^{\infty})$ under $\sim$
 by $[\psi_{\hat{\eps}}^0,\psi_{\hat{\eps}}^{\infty}]$.
\end{itemize}

\begin{theorem} \begin{enumerate}
\item \cite{MRR} The family
$(a(\eps),[\psi_{\hat{\eps}}^0,\psi_{\hat{\eps}}^{\infty}])$ for
some choice of $V_{\rho,\delta}$ is a complete modulus of analytic
classification for the one-parameter prepared family
\eqref{prepared_k_par}, called the {\it modulus of the family
\eqref{prepared_k_par}}. \item In the case of a $k$-parameter
prepared family, the modulus
$(a(\eps,\nu),[\psi_{\hat{\eps},\nu}^0,\psi_{\hat{\eps},\nu}^{\infty}])$
has representatives which depend analytically on the additional parameters $\nu$.
\end{enumerate}
\end{theorem}

 In this
paper we will always use one degree of freedom in the equivalence
relation $\sim$ to manage that $(\psi_{\hat{\eps}}^0)'(0)=1$. To
preserve this property we will limit ourselves to $c'\equiv c$ in
\eqref{mod.6}. It follows from \cite{MRR} that we then have
$(\psi_{\hat{\eps}}^\infty)'(\infty)= \exp(4\pi^2 a(\eps))$.

In practice, we will prefer to work with other presentations of the
modulus, $(\Psi_{\hat{\eps}}^0,\Psi_{\hat{\eps}}^{\infty})$, where
$\Psi_{\hat{\eps}}^{0,\infty} = E^{-1}\circ
\psi_{\hat{\eps}}^{0,\infty} \circ E$, with $E$ defined in
\eqref{def_E}. The functions $\Psi_{\hat{\eps}}^{0,\infty}$ will
have a direct construction from the Fatou coordinates defined in
Section~\ref{section2.3} below.
\medskip

\begin{remark}\label{remark_alpha_omega}
\begin{enumerate}
\item $\delta$ is characterized by the property that for $\arg(\hat{\eps})\in (-\delta,\delta)$
(resp. $\arg(\hat{\eps})\in (2\pi-\delta, 2\pi +\delta)$) there
exists an orbit with $\alpha$-limit (resp. $\omega$-limit) in
$\sqrt{\hat{\eps}}$ (resp. $-\sqrt{\hat{\eps}}$) and $\omega$-limit
(resp. $\alpha$-limit) in $-\sqrt{\hat{\eps}}$ (resp.
$\sqrt{\hat{\eps}}$). Moreover for $\arg(\hat{\eps})\in
(-\delta,\delta)\cup (2\pi-\delta, 2\pi +\delta)$, three cases are
possible for orbits:
\begin{itemize}
\item they have $\alpha$-limit at the repellor and escape the
neighborhood;
\item they have $\omega$-limit at the attractor and the backwards orbits escape the
neighborhood;
\item they have $\alpha$-limit at the repellor and $\omega$-limit at the
attractor.\end{itemize}
\item While in \cite{MRR} it was  shown that we could take $\delta$
as close as $\pi$ as wanted provided $\rho$ be sufficiently small,
we can remark that even with very small $\delta$ we cover a whole
neighborhood of the origin in $\eps$-space. The first point of view,
namely taking $\delta$ close to $\pi$, is linked with the
$1/2$-summability properties in $\eps$ which will be shown below.
However, there will be no need to work with $\delta$ large when we
will study the compatibility condition and indeed the Figures and
estimates will be simpler if we work with $\delta \in
(0,\frac{\pi}2)$. Figure~\ref{extreme} describes the extreme
situations for $\delta$.
\item In fact, it would be natural here to re-express all our
results in terms of germs of functions with respect to the family
of sectors $V_{\rho,\delta}$.  We have not used this language here,
though it is implicit in what we do, as we wished to make clear at
each point the dependence on $\rho$ and $\delta$.  However, we will
make use of arbitrary restrictions of $\rho$ or $\delta$ in what follows
without further comment.
\end{enumerate}
\begin{figure}\begin{center}
\subfigure[$\delta$
large]{\includegraphics[width=4cm]{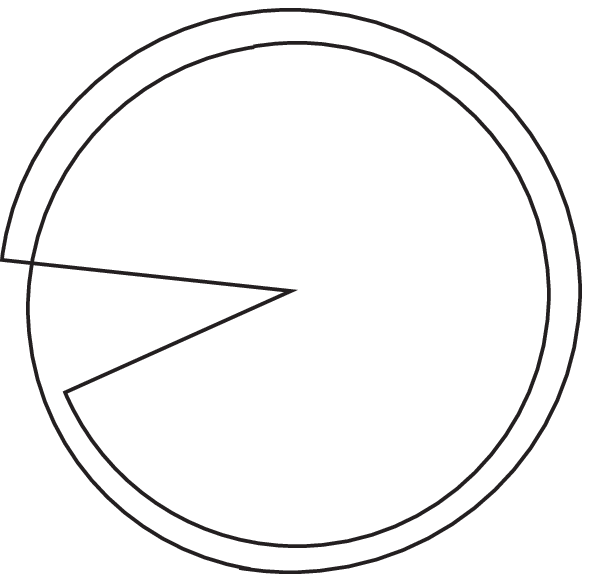}}
\qquad\qquad\subfigure[$\delta$
small]{\includegraphics[width=4cm]{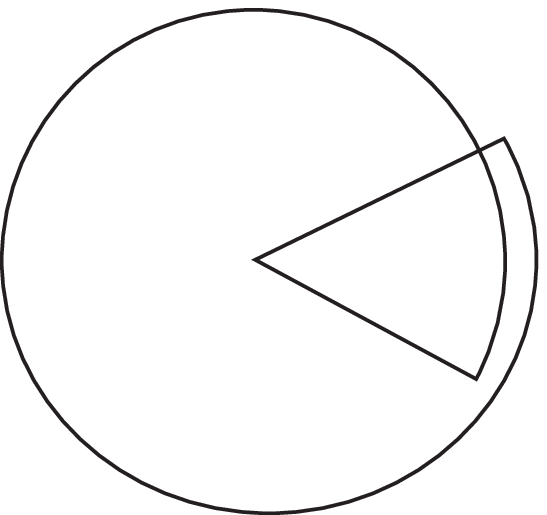}}
\caption{\label{extreme} $V_{\rho,\delta}$ for different sizes of
$\delta$}\end{center}
\end{figure}
\end{remark}

\noindent{\bf The dependence of the modulus on $\eps$.} As stated
above, it is not possible in general to define the modulus so that
its definition depends continuously on $\eps$ in a neighborhood $V$
of the origin. However, given $\del\in(0,\pi)$, we can choose $V$
sufficiently small that the sectorial neighborhood $V_{\rho,\delta}$
projects onto $V$. There exist representatives of the modulus
$\psi^{0,\infty}_{\hat{\eps}}$ which depend analytically on
$\hat{\eps}\ne0$ and continuously on $\hat{\eps}$ at $\hat{\eps}=0$.

In this way we obtain two presentations of the modulus for $\arg
\eps\in (-\del,\del)$. We compare them via the Glutsyuk modulus
defined below.

\medskip

From the unfolded modulus we can deduce the dynamics near each of
the fixed points by means of a renormalized return map when the
multiplier is on the unit circle. Otherwise the renormalized return
maps at the fixed points are linearizable.

\medskip

\noindent{\bf The renormalized return maps.} These maps are defined
on one sphere, for instance $\S^{+}_{\hat{\eps}}$. In the
neighborhood of $\sqrt{\hat{\eps}}$ (resp. $-\sqrt{\hat{\eps}}$)
which we identify to $\infty$ (resp. $0$) on $\S^{+}_{\hat{\eps}}$
we define return maps by iterating $f_\eps$ until the image is
contained in $\S^{+}_{\hat{\eps}}$: given $z\in C^{+}_{\hat{\eps}}$
in the neighborhood of $\sqrt{\hat{\eps}}$ (resp.
$-\sqrt{\hat{\eps}}$) and $w$ its coordinate on
$\S^{+}_{\hat{\eps}}$, let $n\in \mathbb N$ be minimum such that
$f_\eps^n(z)\in C^{+}_{\hat{\eps}}$ and let
$k_{\hat{\eps}}^\infty(w)$ (resp. $k_{\hat{\eps}}^0(w)$) be its
coordinate on $\S^{+}_{\hat{\eps}}$. Then $k_{\hat{\eps}}^\infty$
(resp. $k_{\hat{\eps}}^0$) is the renormalized return map in the
neighborhood of $\sqrt{\hat{\eps}}$ (resp. $-\sqrt{\hat{\eps}}$).
These return maps are given by the composition of the maps
$\psi_{\hat{\eps}}^0$ and $\psi_{\hat{\eps}}^{\infty}$ with a global
transition map $L_{\hat{\eps}}:\S^{-}_{\hat{\eps}}\rightarrow
\S^{+}_{\hat{\eps}}$, the Lavaurs map. The Lavaurs map is an
analytic map from $\CP^1$ to $\CP^1$ fixing 0 and $\infty$. Hence it
is linear, yielding that the nonlinear part of the return map comes
from the unfolding of the two components of the Ecalle-Voronin
modulus. Let us call these two return maps $k_{\hat{\eps}}^0=L_\eps
\circ \psi_{\hat{\eps}}^0$ and $k_{\hat{\eps}}^\infty=L_\eps
\circ\psi_{\hat{\eps}}^\infty$. From \cite{MRR}, they have
multipliers
$$\begin{cases}
(k_{\hat{\eps}}^0)'(0)=\exp\left(\frac{4\pi^2}{\mu_{\hat{\eps}}^0}\right),\\
(k_{\hat{\eps}}^\infty)'(\infty)=\exp\left(\frac{4\pi^2}{\mu_{\hat{\eps}}^\infty}
\right).\end{cases}$$ In the Glutsyuk domain, namely
$\arg\hat{\eps}\in(-\delta,\delta)\cup(2\pi-\delta,2\pi+\delta)$,
$(k_{\hat{\eps}}^0)'(0)$ and $(k_{\hat{\eps}}^\infty)'(\infty)$ are
exponentially small or large in $\sqrt{\hat{\eps}}$ ($\sim \exp(\pm
\frac{C}{|\sqrt{\hat{\eps}}|})$)
\medskip

 \noindent{\bf The Glutsyuk modulus.} The Glutsyuk modulus is
defined for small values of $\eps$ with $\arg{\eps}\in
(-\delta,\delta)$ and we will decide to work with $\delta\in
(0,\frac{\pi}2)$. For such $\eps$, the fixed points $\sqrt{\eps}$
and $-\sqrt{\eps}$ are respectively hyperbolic repeller and
attractor. Moreover, as stated in Remark~\ref{remark_alpha_omega},
there are orbits of $f_\eps$ in $U$ which have $\sqrt{\eps}$ (resp.
$-\sqrt{\eps}$) as $\alpha$- (resp. $\omega$-) limit set.

\begin{figure}
\begin{center}
\includegraphics[height=3.5cm]{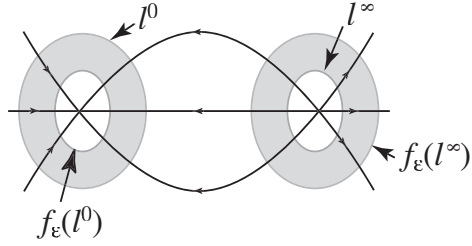}
\caption{The fundamental domains in the Glutsyuk modulus
 }
\label{fig:Glutsyuk_modulus}
\end{center}
\end{figure}

We take two closed curves $l^0$ and $l^\infty$ surrounding
$-\sqrt{\eps}$ and $\sqrt{\eps}$.  Since the fixed points are
hyperbolic, we can choose $l^{0,\infty}$ so that the region
$C_\eps^{0,\infty}$ between the curves $l^{0,\infty}$ and their
images $f_\eps(l^{0,\infty})$ are homeomorphic to annuli (see
Figure~\ref{fig:Glutsyuk_modulus}). We identify $l^{0,\infty}$ and
$f_\eps(l^{0,\infty})$ to get two tori $\mathbb T_\eps^{0,\infty}$
which represent the local orbit space of the hyperbolic fixed
points. Since $f_\eps$ has connecting orbits, we can iterate
$f_\eps$ in such a way as to identify a collar of $\mathbb
T_\eps^\infty$ with a collar in $\mathbb T_\eps^0$.  In the limit
$\eps=0$, the tori become pinched and the map between the collars
splits into two maps between the respective ends of the pinched
tori.  The moduli of the tori depend on $a(\eps)$ and $\eps$ and can
be derived directly from the multipliers of the fixed points.

\begin{figure}\begin{center}
\subfigure[$\eps=0$]{\includegraphics[width=4cm]{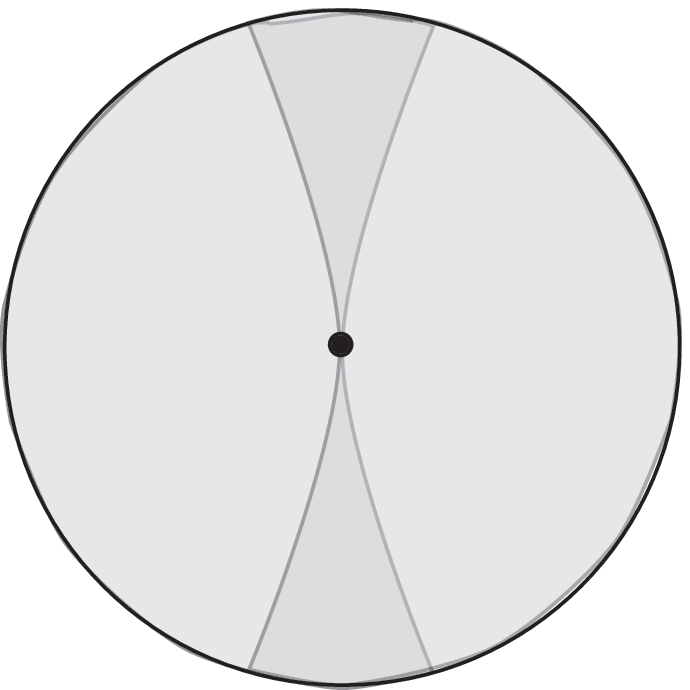}}
\qquad\qquad\subfigure[$\arg\hat{\eps}\in (
-\delta,\delta)$]{\includegraphics[width=4cm]{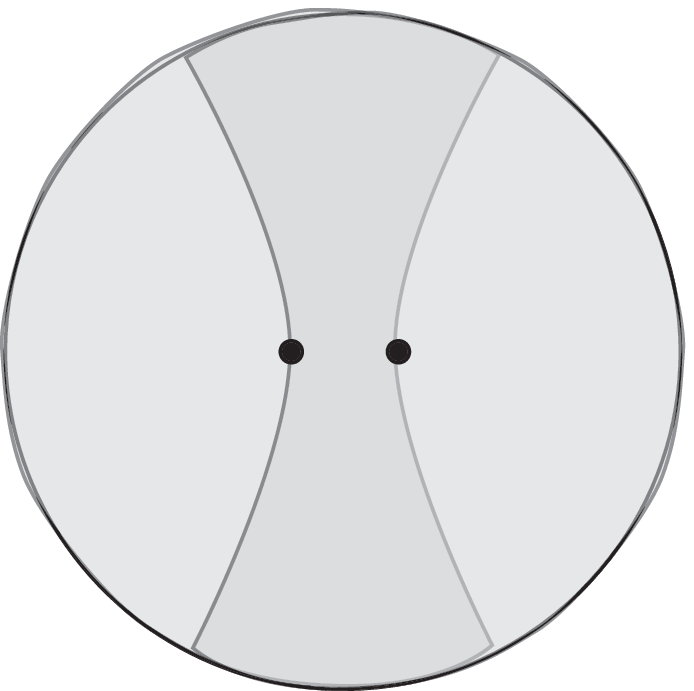}}
\caption{\label{domain_Glutsyuk} The domain of the Glutsyuk modulus
in the original coordinate $z$}\end{center}
\end{figure}

This map is one presentation of the Glutsyuk modulus. A more usual
but less geometric presentation is equivalent to the covering map of
the above construction. That is, we describe the Glutsyuk modulus in
the following way.  Since the two points are hyperbolic, there
exists in the neighborhood of each fixed point $\pm \sqrt{\eps}$ a
diffeomorphism $\varphi_{\eps}^{\pm}$ conjugating $f_\eps$ to the
model, i.e. the time one map of \eqref{mod.3}. For a sufficiently
small choice of $V_{\rho,\delta}$ the domains of
$\varphi_{\eps}^{\pm}$ overlap allowing to define the map
$$ \varphi_\eps^G=\varphi_\eps^-\circ (\varphi_\eps^+)^{-1}.$$
If we call $$V_G(\rho)=\{\eps;|\eps|<\rho, \arg\eps\in
(-\delta,\delta)\},$$ then it is easy to verify that, for
sufficiently small $ \rho$, $(\varphi_\eps^G)_{\eps\in V_G(\rho)}$
is an analytic invariant of the family $f_\eps$ under analytic
families of change of coordinates preserving the canonical
parameter. The Glutsyuk modulus is unique up to composition on the
left and on the right by time $t$ maps  $v_\eps^t$ of the vector
field \eqref{mod.3}. The family $(\varphi_\eps^G)_{\eps\in
V_\eta(\rho)}$ gives the presentation of the Glutsyuk modulus.
The domain for
$\varphi_\eps^G$ appears in Figure~\ref{domain_Glutsyuk}.
\smallskip

In practice we will also need to work with other presentations
obtained with the use of Fatou coordinates described now.

\subsection{Fatou coordinates and other presentations of the
modulus}\label{section2.3}

On $U$ we make the change of coordinate $Z=p_\eps^{-1}(z)$ defined
by
\begin{equation}
Z=p_\eps^{-1}(z)=\begin{cases}\frac1{2\sqrt{\eps}}\ln\frac{z-\sqrt{\eps}}{z+\sqrt{\eps}},&\eps\neq0,\\
-\frac1{z},&\eps=0.\end{cases}\label{eq_p}\end{equation}

In the $Z$-coordinate, the map $f_\eps$ is transformed to $F_\eps$
which is very close to the translation $T_1$. {\it Fatou
coordinates} are changes of coordinates $Z\mapsto W$ defined on
simply connected domains in $Z$-space called {\it translation
domains} and conjugating $F_\eps$ to $T_1$.

A translation domain is constructed by choosing an admissible line
$\ell$ in the image of $p_\eps^{-1}(U)$ in $Z$-space, i.e. a line
such that $\ell$ and $F_\eps(\ell)$ are disjoint and bound a strip
in $p_\eps^{-1}(U)$, and by saturating this strip under the action
of $F_\eps$,

Given an admissible line $\ell$ in $Z$-space, the associated Fatou
coordinate is uniquely defined up to left  composition with a
translation.

The corresponding presentation of the modulus is a comparison of two
Fatou coordinates.

\begin{figure}[!h]
\begin{center}
\includegraphics[width=9cm]{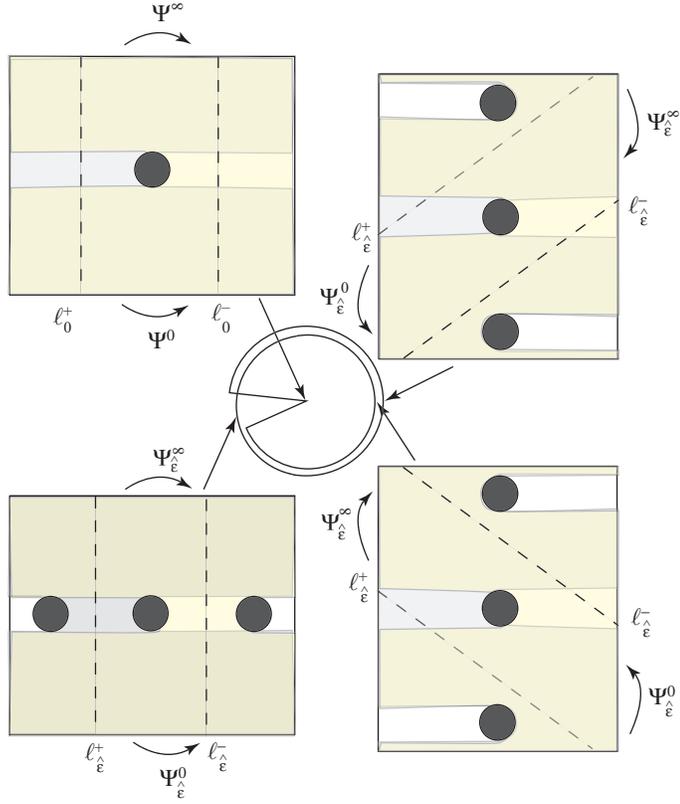}
\caption{Fatou coordinates in Lavaurs point of view }
\label{strips_Lavaurs}
\end{center}
\end{figure}

\begin{figure}[!h]
\begin{center}
\includegraphics[width=9cm]{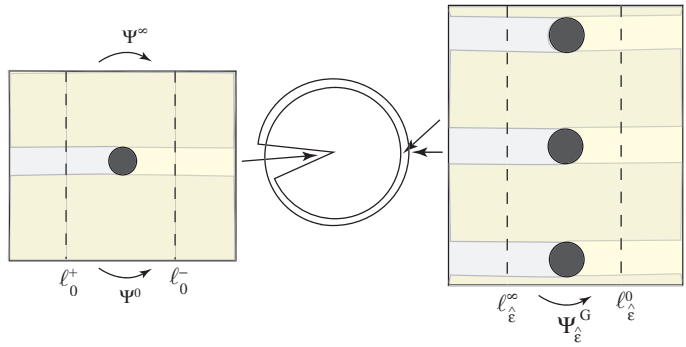}
\caption{Fatou coordinates in Glutsyuk point of view }
\label{strips_Glutsyuk}
\end{center}
\end{figure}

In the Lavaurs point of view, we compare two Fatou coordinates
$\Phi_{\hat{\eps}}^{\pm}$ defined on translation domains constructed
with slanted lines $\ell_{\hat{\eps}}^{\pm}$ passing between two
holes as in Figure~\ref{strips_Lavaurs}, while in the Glutsyuk point
of view we compare two Fatou coordinates
$\Phi_{\hat{\eps}}^{0,\infty}$ defined on translation domains
constructed with lines $\ell_{\hat{\eps}}^{0, \infty}$ parallel to
the line of holes as in Figure~\ref{strips_Glutsyuk}.

\begin{definition}
\begin{description}
\item{(1)} The modulus in the Lavaurs point of view is given by
\begin{equation}\Psi_{\hat{\eps}}=\Phi_{\hat{\eps}}^-\circ
(\Phi_{\hat{\eps}}^+)^{-1},\label{modulus_Lavaurs}\end{equation} up
to composition with a translation on the left and a translation on
the right. Since the domain is disconnected, this map is indeed
described by the two maps $\Psi_{\hat{\eps}}^0$ (resp.
$\Psi_{\hat{\eps}}^\infty$) defined for $Im(W)< -Y_0$ (resp.
$Im(W)>Y_0$). We also use the alternative presentation
\begin{equation}
\psi_{\hat{\eps}}=E\circ \Psi_{\hat{\eps}}\circ
E^{-1}.\label{eq:small_psi}\end{equation} Here the domain of
$\psi_{\hat{\eps}}$ is the union of a neighborhood of $0$ and a
neighborhood of $\infty$ on $\mathbb C\mathbb P^1$. The respective
restrictions of $\psi_{\hat{\eps}}$ to these neighborhoods are noted
$\psi_{\hat{\eps}}^0$ and $\psi_{\hat{\eps}}^\infty$. These clearly
coincide with the definitions  $\psi_{\hat{\eps}}^{0,\infty}$ given
previously when considering  the spheres $\S^{\pm}$.
\item{(2)} When $\arg \hat{\eps}\in
(-\eta,\eta)$, there exist Fatou coordinates
$\Phi_{\hat{\eps}}^{0,\infty}$ associated to translation domains
defined with lines parallel to the holes as in
Figure~\ref{strips_Glutsyuk}. (We call these Fatou coordinates the
Fatou Glutsyuk coordinates.) The modulus in the Glutsyuk point of
view is then given by
\begin{equation}\Psi_{\hat{\eps}}^G=\Phi_{\hat{\eps}}^0\circ
(\Phi_{\hat{\eps}}^\infty)^{-1},\label{modulus_Glutsyuk}\end{equation}
up to composition with a translation on the left and a translation
on the right.  The maps $\varphi^G$ (resp. $\varphi^\pm$) mentioned
previously are just the push forward of $\Psi^G$ (resp.
$\Phi^{0,\infty}$) via $p_\eps$.
 \end{description}\end{definition}

\begin{remark}
From the uniqueness of the Fatou Glutsyuk coordinates, when $\arg
\hat{\eps}\in (2\pi-\eta,2\pi+\eta)$  the Glutsyuk modulus is
defined by $\Psi_{\hat{\eps}}^{G}=\Phi_{\hat{\eps}}^\infty\circ
(\Phi_{\hat{\eps}}^0)^{-1}$. \end{remark}

\section{The local realization}\label{sect:local}

We will work with parameter values $\hat{\eps}$ in some
$V_{\rho,\delta}$, as in \eqref{sector_V}. Unless specified, we will
always suppose that the sectors $V_{\rho,\delta}$ contain $\eps=0$.
It is clear that we can extend our definition of the modulus of a
family $f_\eps$, to cover the case of a ramified prepared family
$f_{\hat{\eps}}$ defined for $\hat{\eps}\in V_{\rho,\delta}$, where
$f_{\hat{\eps}}$ is analytic in $V_{\rho,\delta}$ and locally of the
form $$f_{\hat{\eps}}(z)= z+ (z^2-\eps)(1+h_{\hat{\eps}}(z)),$$ with
$h_{\hat{\eps}}(z)=O(|\hat{\eps},z|)$.

We denote
$$\begin{cases}\mu^0(\hat{\eps})=-\frac{2\sqrt{\hat{\eps}}}{1-a(\eps)\sqrt{\hat{\eps}}}\\
\mu^\infty(\hat{\eps})=\frac{2\sqrt{\hat{\eps}}}{1+a(\eps)\sqrt{\hat{\eps}}},\end{cases}$$
and hence $\mu^{0,\infty}(e^{2\pi i }\hat{\eps})=
\mu^{\infty,0}(\hat{\eps})$ and
$$a(\eps)=\frac1{\mu^0(\hat{\eps})}+\frac1{\mu^\infty(\hat{\eps})}$$
(which is not ramified in $\eps$!).

For such ramified families we prove the following theorem:

\begin{theorem}\label{thm:local}
Let $\delta\in (0,\frac{\pi}2)$, and consider a germ of analytic
function $a(\eps)$ at the origin. Let  $V_{\rho,\delta}$ be a
sectorial neighborhood of the origin in the universal covering of
$\eps$-space punctured at the origin of the form \eqref{sector_V},
such that $a(\eps)$ has a representative on $V_{\rho,\delta}$.

Let $\Psi_{\hat{\eps}}^0(W)$ (resp. $\Psi_{\hat{\eps}}^\infty(W)$)
be families of germs of analytic diffeomorphisms at $Im(W)=-\infty$
(resp. $Im(W)=+\infty$) having representatives $\Psi_{\hat{\eps}}^0:
R^0 \rightarrow \mathbb C$ (resp. $\Psi_{\hat{\eps}}^\infty:
R^\infty \rightarrow \mathbb C$) defined for $\hat{\eps}\in
V_{\rho,\delta}$ in domains $R^0=\{Im(W) <-Y_0\} $ (resp.
$R^\infty=\{Im(W)
>Y_0\}$) for some $Y_0>0$ and such that
\begin{description}
\item{(i)} $\Psi_{\hat{\eps}}^{0,\infty}$ depend analytically on $\hat{\eps}\in V_{\rho,\delta}\setminus\{0\}$ and
have continuous limits when $\hat{\eps}\to 0$.
\item{(ii)} $\Psi_{\hat{\eps}}^{0,\infty}$ commute with $T_1$.
\item{(iii)} We have \begin{equation}\begin{cases} \Psi_{\hat{\eps}}^0(W)=W + O(exp(2\pi i W)), \ \ Im(W) << 0,\\
   \Psi_{\hat{\eps}}^\infty(W)=W- 2\pi i a(\eps)+O(exp(2\pi i W)), \ \ Im(W) >> 0.
  \end{cases}\label{shift_a}\end{equation}
\end{description}
 Then for any $\delta'\in (0,\delta)$
there exists $\rho'\in(0,\rho]$,  a neighborhood $U$ of the origin
in $\mathbb C$ containing the two points $\pm \sqrt{\hat{\eps}}$ and
a family of analytic diffeomorphisms $f_{\hat{\eps}}(z):
U\rightarrow \mathbb C$ depending on $\hat{\eps}\in
V_{\rho',\delta'}$, such that:
\begin{itemize} \item For all $\hat{\eps}\in V_{\rho',\delta'}$, $f_{\hat{\eps}}(z)$ has exactly two fixed
points located at $\pm\sqrt{\hat{\eps}}$ and  is of the form
$$f_{\hat{\eps}}(z)= z+ (z^2-\eps)(1+h_{\hat{\eps}}(z)),$$
with $h_{\hat{\eps}}(z)=O(|\hat{\eps},z|)$.
\item
$f_{\hat{\eps}}'(\sqrt{\hat{\eps}})=\exp(\mu^\infty)$ and
$f_{\hat{\eps}}'(-\sqrt{\hat{\eps}})=\exp(\mu^0)$. (So
$f_{\hat{\eps}}$ is prepared.)
\item $f_{\hat{\eps}}(z)$ depends analytically of $\hat{\eps}\in
V_{\rho',\delta'}\setminus\{0\}$ and has a continuous limit when
$\hat{\eps}\to 0$.
\item The modulus of $f_{\hat{\eps}}$ is given by
$[\Psi_{\hat{\eps}}^0,\Psi_{\hat{\eps}}^\infty]$.
\end{itemize}
If the functions $a(\eps,\nu)$ and
$\Psi_{\hat{\eps},\nu}^{0,\infty}$ depend analytically on a
multi-parameter $\nu$, then the function $f_{\hat{\eps},\nu}$
depends analytically on $\nu$.
\end{theorem}

For the proof of the theorem we will concentrate on the
one-parameter case. It will be obvious that all steps will be
analytic in extra parameters.

 The following lemma will be used in the proof and
elsewhere in the paper.

\begin{lemma}\label{lemma_estimates}\begin{description}
\item{(i)} We consider families of germs of analytic diffeomorphisms
$\Psi_{\hat{\eps}}^0(W)$ (resp. $\Psi_{\hat{\eps}}^\infty(W)$) at
$Im(W)=-\infty$ (resp. $Im(W)=+\infty$) commuting with $T_1$, having
representatives $\Psi_{\hat{\eps}}^0: R^0 \rightarrow \mathbb C$
(resp. $\Psi_{\hat{\eps}}^\infty: R^\infty \rightarrow \mathbb C$)
defined for $\hat{\eps}\in V_{\rho,\delta}$ in domains $R^0=\{Im(W)
<-Y_0\} $ (resp. $R^\infty=\{Im(W)
>Y_0\}$) for some $Y_0>0$ and such that $\Psi_{\hat{\eps}}^{0,\infty}$
depend analytically on $\hat{\eps}\in V_{\rho,\delta}$ and have
continuous limits when $\hat{\eps}\to 0$. Let \begin{equation}
\begin{cases}
 \Psi_{\hat{\eps}}^0(W)=W+ \sum_{n\leq-1}b_n(\hat{\eps})\exp(2\pi i
 n W),\\
\Psi_{\hat{\eps}}^\infty(W)=W-  2\pi i
a(\eps)+\sum_{n\geq1}c_n(\hat{\eps})\exp(2\pi i
 n W),\end{cases}\label{series_psi}\end{equation}
let $\beta>0$ be small  and let $$
 \begin{cases}M^0= \max_{Im(W)\leq-Y_0-\beta}\left|\Psi_{\hat{\eps}}^0(W)-W\right|,\\
 M^\infty=\max_{Im(W)\geq Y_0+\beta}\left|\Psi_{\hat{\eps}}^\infty(W)-W+2\pi
 i a(\eps)\right|.\end{cases}$$
 Then
 $$\begin{cases} \left|b_n(\hat{\eps})\right|<M^0\exp(-2\pi n(Y_0+\beta)), &n \leq
 -1,\\
 \left|c_n(\hat{\eps})\right|<M^\infty\exp(2\pi
 n(Y_0+\beta)),&n\geq1.\end{cases}$$
 The series  $\Psi_{\hat{\eps}}^0$ (resp. $\Psi_{\hat{\eps}}^\infty$) in \eqref{series_psi} is absolutely
 convergent for $Im(W)\leq-Y_0-\beta$ (resp. $Im(W)\geq Y_0+\beta$).
 Moreover there exists a constant $N=N(\beta)$ depending only on $\beta$ such
 that
 \begin{equation}
 \begin{cases}\left|\Psi_{\hat{\eps}}^0(W)-W\right|<M^0N(\beta)\exp(2\pi (Y_0+\beta +
 Im(W)), &Im(W)<-Y_0-2\beta,\\
\left|\Psi_{\hat{\eps}}^\infty(W)-W+2\pi i a(\eps)\right|<M^\infty
N(\beta)\exp(2\pi (Y_0+\beta-
 Im(W)),
 &Im(W)>Y_0+2\beta.\end{cases}\label{estime_Psi}\end{equation}\item{(ii)}
For any $\beta>0$, the maps $\Psi^0$ and
$\Psi^\infty$ are uniformly continuous in the region $\{|Im W|>Y_0+\beta\}\times V_{\delta,\rho}$.
\item{(iii)} The image of $\{Im W< -Y_0\}$ (resp. $\{Im W>Y_0\}$) under $\Psi_{\hat{\eps}}^0$ (resp. $\Psi_{\hat{\eps}}^
\infty$) contains some half-plane of the form
 $\{Im W< -Y_1\}$ (resp. $\{Im W> Y_1\}$).
 \end{description}
 \end{lemma}
 \noindent{\scshape Proof.}\begin{description}
 \item{(i)}
 This follows from the fact that
 $$b_n=\int_{X_0-i(Y_0+\beta)}^{X_0+1-i(Y_0+\beta)}(\Psi_{\hat{\eps}}^0-id)(X-i(Y_0+\beta))\exp(-2\pi
 i n(X-i(Y_0+\beta)))dX,$$
 and similarly for $c_n$.
 \item{(ii)}
 This follows from the fact that the maps commute with $T_1$
and have a definite limit as $|Im W|\to \infty$ or $\hat{\eps}\to 0$.
\item{(iii)} This follows from the fact that $\Psi_{\hat{\eps}}^{0,\infty}$ commute with $T_1$.\hfill $\Box$\end{description}

\medskip

In the rest of the paper we will choose our different sectors in $z$-space
 (corresponding to strips in $W$-space), so that any region where we
 need to consider $\Psi_{\hat{\eps}}^{0}$ (resp. $\Psi_{\hat{\eps}}^{\infty}$)
 is located inside $Im W <-Y_0-2\beta$ (resp. $Im W> Y_0+2\beta$) for some
 suitable $\beta$, so that the estimates of Lemma~\ref{lemma_estimates} will
 always be valid.

 \medskip

 \noindent{\scshape Proof of Theorem~\ref{thm:local}.} We choose any $\delta'\in (0,\delta)$. Working
with $\delta'$ instead of $\delta$ allows to consider
$\arg(\hat{\eps})$ to vary inside a compact set and hence to yield
uniform estimates in $\arg(\hat{\eps})$. We look for a neighborhood
$U=B(0,r)$ of the origin in $z$-space. The final choice of $r$ and
$\rho'$ considered before will be done in several steps throughout
the proof. We consider the regions $R^0$ and $R^\infty$ in $W$-space
and the multivalued mapping:

\begin{equation}W=q_{\hat{\eps}}^{-1}(z)=\begin{cases}\frac1{2\sqrt{\hat{\eps}}}
\ln\frac{z-\sqrt{\hat{\eps}}}{z+\sqrt{\hat{\eps}}} +
\frac{a(\eps)}2\ln(z^2-\eps),&\hat{\eps}\neq0,\\-\frac1{z}+a(0)\ln(z),
&\hat{\eps}=0.\end{cases}\label{function_q}\end{equation}

While the
inverse $q_{\hat{\eps}}$ exists, it cannot be described by a simple
formula.

Note that the function $q_{\hat{\eps}}^{-1}(z)$ is simply the time
of the vector field \eqref{mod.3}. The map $q_{\hat{\eps}}^{-1}$ has
the property that the restriction  of $q_{\hat{\eps}}^{-1}\circ
p_{\hat{\eps}}$ to a translation domain is a Fatou coordinate of the
model family, namely a conjugacy of $p_{\hat{\eps}}^{-1}$ applied to
the model with the translation by 1. (Recall that the model is the
time one map of the vector field \eqref{mod.3}).

The function $q_{\hat{\eps}}^{-1}(z)$ is a multi-valued analytic
function of two variables outside the set $\{(z,\eps)|z^2-\eps=0\}$. For
$\eps=0$, the function $q_{\hat{\eps}}^{-1}$ is not a global
diffeomorphism if $a(0)\neq0$. So we should not consider it over the
whole complex plane and it is better to limit ourselves to sectors
in a small neighborhood $U=B(0,r)$ of the origin in $z$-space. The
function $q_{\hat{\eps}}^{-1}$ is ramified both at
$\pm\sqrt{\hat{\eps}}$. Moreover when $a\neq 0$ a cut cannot simply
be taken between $-\sqrt{\hat{\eps}}$ and $\sqrt{\hat{\eps}}$ since
there is a global ramification when one makes a turn on $C(0,r)$.

Although it is difficult to visualize the map $q_{\hat{\eps}}^{-1}$ directly,
it can be pictured more easily when lifted to the $Z$-plane
via $p_{\hat{\eps}}$.  Here it will be a multi-valued function, whose
difference in value when continued around any of the holes in the $Z$-plane
is just $2\pi i a(\eps)$.  The absolute difference between $W$ and $Z$ in a
simply connected region is bounded by $2|a(\eps)|\ln(r)$.  Thus, if we
restrict our attention to a simply connected region, $W$-space can be thought
of as a small distortion of $Z$-space.

The distance vector between the centers of two holes is of the order
\begin{equation}\alpha=\alpha_{\hat{\eps}}=\frac{\pi i}{\sqrt{\hat{\eps}}}.\label{eq_alpha}\end{equation}
Hence, the distance between two consecutive holes is of the order of $|\alpha|$  and the
radius of holes is of the order of $\frac1{r}$ for small $\eps$.

\begin{figure}[!h]
\begin{center}
\includegraphics[width=11cm]{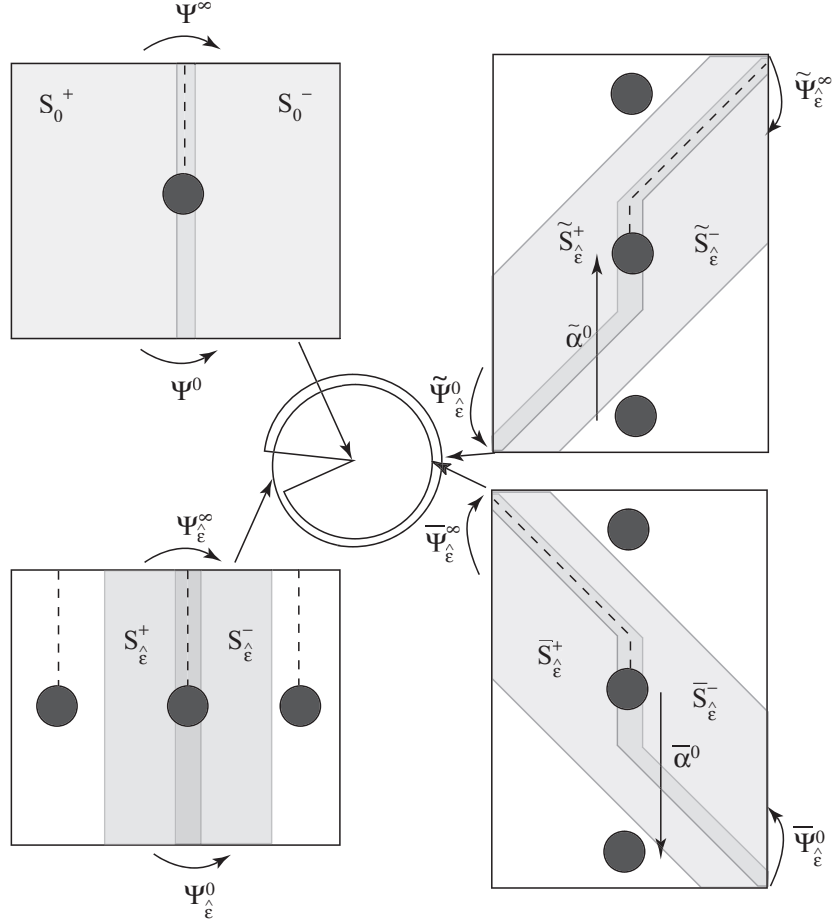}
\caption{The choice of strips.  The dotted lines represent the cuts.
} \label{strips}
\end{center}
\end{figure}

As suggested above, we will limit ourselves to simply connected
regions on which $q_{\hat{\eps}}^{-1}$ and its inverse
$q_{\hat{\eps}}$ are well defined.  We choose two strips
$S_{\hat{\eps}}^\pm$ located on each side of the principal hole as
in Figure~\ref{strips}. The choice of the strips and of $r$ and
$\rho'$ is given in the following Lemma.

\begin{lemma}\label{lemma_strips_adjusted}
For $\delta\in (0,\frac{\pi}2)$ there exists $\rho'>0$ sufficiently
small such that for $|\hat{\eps}|<\rho'$ and $\arg({\hat{\eps}})\in
(-\delta, 2\pi+\delta)$ there exist adjusted strips constructed as
follows.
\begin{itemize}
\item The total width of the union of the two strips in the direction of
the line of holes is $\frac{3|\alpha|}2$. \item The horizontal width
of the intersection is fixed and equal to $2h$ for some positive
constant $h<\frac1{2r}$ (recall that the radius of the holes is
approximately $\frac1{r}$). \item Let
\begin{equation}\theta=\frac12\left(\frac{\pi}{2}+\arg(\sqrt{\hat{\eps}})\right).\label{equation_theta}
\end{equation} The
strips are bounded on one side by a slanted line of slope
\begin{equation}
t=-\tan\theta.\label{equation_tau}\end{equation} On the other side
they are bounded by a vertical segment $Re W=\pm h$ of total height
$|\alpha|/4$.
From the two
endpoints of the segment we continue with two half lines with slope
$-\tan\theta$ as drawn in Figure~\ref{strips}.
\item The
radius $r$ is chosen sufficiently small so that the intersection
part of the strips outside the fundamental holes is located in the
region $|Im W|>Y_0+2\beta$ where we can apply the estimates of
Lemma~\ref{lemma_estimates}.\end{itemize}
\end{lemma}

\noindent{\scshape Proof.} We only discuss the
range $\theta\in(\frac{\pi}{8}, \frac{\pi}4)$ where both the
strip and the line of holes have negative slope.  The case
$\theta\in(\frac{3\pi}{8}, \frac{7\pi}8)$ is similar.

If the holes are of negligible width, then it is a simple matter of
geometry that the construction above is valid if $\tan(\theta)>1/3$,
and is therefore satisfied in our range.  By choosing $\rho'$
sufficiently small, we can make the effective size of the holes
arbitrarily small and hence the result follows.
 \hfill $\Box$
\medskip

\begin{figure}[!h]
\begin{center}
\includegraphics[width=10cm]{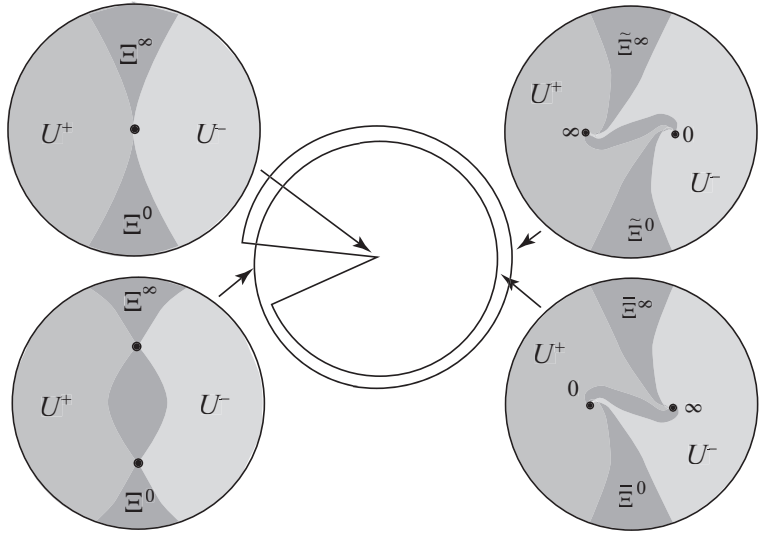}
\caption{The sectors $U_{\hat{\eps}}^\pm$ and their intersection }
\label{sectors}
\end{center}
\end{figure}

We now consider the images of the two strips, $S_{\hat{\eps}}^\pm$, under
$q_{\hat{\eps}}$.  These yield two sectors $U_{\hat{\eps}}^\pm$ whose union
is $U\setminus\{\pm\sqrt{\hat{\eps}}\}$. For $\eps=0$ the
intersection $U_0^+\cap U_0^-$ is formed of two narrow sectors
$U_0^0$ and $U_0^\infty$ with vertex at $0$ and ending on the
boundary of $U$, while for $\hat{\eps}\neq0$ the intersection is
formed of three parts: two sectors $U_{\hat{\eps}}^0$ (resp.
$U_{\hat{\eps}}^\infty$) with vertex at $-\sqrt{\hat{\eps}}$ (resp.
$\sqrt{\hat{\eps}}$) and ending on the boundary of $U$ and one
crescent $U_{\hat{\eps}}^C$ with its two endpoints at
$\pm\sqrt{\hat{\eps}}$ (Figure~\ref{sectors}).
The crescent $U_{\hat{\eps}}^C$ comes from the fact that $q_{\hat{\eps}}^{-1}$ is
multivalued and approximately periodic with a period of the order of
$\alpha_{\hat{\eps}}=\frac{\pi i}{\sqrt{\hat{\eps}}}$ and the width
of the union of the two strips in the direction of
$\alpha_{\hat{\eps}}$ is  $\frac32\alpha_{\hat{\eps}}$.

On $U_{\hat{\eps}}^{\pm}$ we can define $q_{\hat{\eps}}^{-1}$ in a
uniform way, which we call $q_{\hat{\eps},\pm}^{-1}$. The
determinations are chosen so that $q_{\hat{\eps},\pm}^{-1}$ agree on
$U_{\hat{\eps}}^0$. If we take the analytic extension of
$q_{\hat{\eps},-}^{-1}$  after making one turn in the positive
direction around $-\sqrt{\hat{\eps}}$, then the extension has the
form $T_{\frac{2\pi i }{\mu^0}}\circ q_{\hat{\eps},-}^{-1}$.

Let \begin{equation}\begin{cases}\Xi_{\hat{\eps}}^0= id
+\xi_{\hat{\eps}}^0=q_{\hat{\eps},+}\circ \Psi_{\hat{\eps}}^0\circ
q_{\hat{\eps}, +}^{-1},\\ \Xi_{\hat{\eps}}^\infty= id
+\xi_{\hat{\eps}}^\infty=q_{\hat{\eps},+}\circ
(\Psi_{\hat{\eps}}^\infty+2\pi i a(\eps))\circ
q_{\hat{\eps},+}^{-1},\end{cases}\label{def_chi}\end{equation} which
are defined respectively in regions containing
$U_{\hat{\eps}}^{0,\infty}$. For future reference, we also take
$\Xi_{\hat{\eps}}^C = id$.

We construct an abstract complex manifold $M_{\hat{\eps}}$ by gluing
$U_{\hat{\eps}}^\pm$ along their intersection.  More precisely, let
$z^\pm$ be the coordinates on $U_{\hat{\eps}}^\pm$. Then we identify
\begin{equation}
z^-=\begin{cases} z^++ \xi_{\hat{\eps}}^0(z^+)= \Xi_{\hat{\eps}}^0(z^+), & z^+\in U_{\hat{\eps}}^0,\\
z^+ + \xi_{\hat{\eps}}^\infty(z^+)=\Xi_{\hat{\eps}}^\infty(z^+), &
z^+\in U_{\hat{\eps}}^\infty,
\\z^+ =\Xi_{\hat{\eps}}^C(z^+), & z^+ \in U_{\hat{\eps}}^C,\end{cases}\label{def_Xi}\end{equation}
deleting those points in $U_{\hat{\eps}}^-$ which are in $U_{\hat{\eps}}^{0,\infty}$ but
are not in the image of $\Xi_{\hat{\eps}}^{0,\infty}$ to ensure that the space we get
is Hausdorff.

This gluing is well-defined, since near $Im(W)=\pm\infty$,
$\Psi_{\hat{\eps}}^{0,\infty}$ is close to a translation. It is easy
to take $r$ and $|\hat{\eps}|$ sufficiently small so that this
translation is very small compared to the width of the strips: the
first condition ($r$ small) ensures that the balls of
Figure~\ref{strips} are  sufficiently large, while the second
($|\hat{\eps}|$ small) guarantees that the strips and their
intersection can be chosen wide.

The map $T_1$ on the strips lifts to a well-defined holomorphic map
$F_{\hat{\eps}}$ on $M_{\hat{\eps}}$, due to the fact that
$\Psi_{\hat{\eps}}^{0,\infty}$ commute with $T_1$.  We want to show
that $M_{\hat{\eps}}$ is conformally equivalent to a disk in
$\mathbb C$, $D_{\hat{\eps}}$,  punctured at $\pm\sqrt{\hat{\eps}}$.
For this we first find a smooth map from $M_{\hat{\eps}}$ to $\C$,
and then use the Ahlfors-Bers theorem to correct this to a
holomorphic map.

Having done this, the image of the map $F_{\hat{\eps}}$ is just the diffeomorphism $f_{\hat{\eps}}$
we are seeking.  Indeed, the $W$ coordinate considered as a multi-valued function in the
$Z$-plane gives Fatou coordinates for $f_{\hat{\eps}}$, and our gluings $\Xi_{\hat{\eps}}$
can be written as
\begin{equation}\label{qm_qp}
\begin{cases}\Xi_{\hat{\eps}}^0 = id+\xi_{\hat{\eps}}^0=q_{\hat{\eps},+}\circ \Psi_{\hat{\eps}}^0\circ
q_{\hat{\eps}, +}^{-1} = q_{\hat{\eps},-}\circ \Psi_{\hat{\eps}}^0\circ
q_{\hat{\eps}, +}^{-1},\\
\Xi_{\hat{\eps}}^\infty= id+\xi_{\hat{\eps}}^\infty=q_{\hat{\eps},+}\circ
(\Psi_{\hat{\eps}}^\infty+2\pi i a(\eps))\circ
q_{\hat{\eps},+}^{-1}=q_{\hat{\eps},-}\circ\Psi_{\hat{\eps}}^\infty \circ
q_{\hat{\eps},+}^{-1}, \\
\Xi_{\hat{\eps}}^C = id = q_{\hat{\eps},+}\circ q_{\hat{\eps},+}^{-1}=
q_{\hat{\eps},-}\circ T_{\frac{2\pi i}{\mu_{\hat{\eps}}^0}}\circ
q_{\hat{\eps},+}^{-1}.
\end{cases}
\end{equation}
That is, the gluings correspond exactly to the fact that the modulus
of $f_{\hat{\eps}}$ is
$(a(\eps),[\Psi_{\hat{\eps}}^0,\Psi_{\hat{\eps}}^\infty])$.  The
punctures in the disc $D_{\hat{\eps}}$, correspond to the critical
points of the map $f_{\hat{\eps}}$, and their multipliers and thence
$a(\eps)$ can be similarly derived from
$\Xi_{\hat{\eps}}^{0,\infty,C}$. The rest of the statements of
Theorem~\ref{thm:local} follow.

\smallskip
We therefore wish to map $M_{\hat{\eps}}$ to $\C$ in a smooth way.
We express this map via the coordinate patches of $M_{\hat{\eps}}$
on the $W$-plane.  We work first with a fixed $\hat{\eps}$.

Let $\varphi: \mathbb R\rightarrow \mathbb [0,1]$ be a
$C^\infty$ monotonic increasing function such that
$$\varphi\equiv\begin{cases} 0, &x\leq 0,\\
 1,&x\geq 1.\end{cases}$$
 Hence for each $n$ there exists a constant
$C_n$ such that
\begin{equation}\label{varphibound}
  |\varphi^{(n)}|\leq C_n.
\end{equation}

Writing $W=X+iY$, we take two $C^\infty$ curves $X=\ell_i(Y)$ with
$\ell_2(Y)=\ell_1(Y)+h$ which lie within the   intersection of the
two strips outside the holes, and take
\begin{equation}N_{\hat{\eps}}(X+iY)=\varphi\left(\frac{X-\ell_1(Y)}{h}\right)\label{eq_N}
\end{equation}
and
$$\begin{cases}
\Theta_{\hat{\eps}}^{-}(x,y)= N_{\hat{\eps}}\circ q_{\hat{\eps}}^{-1},\\
\Theta_{\hat{\eps}}^{+}(x,y)=
1-\Theta_{\hat{\eps}}^{-}(x,y),\end{cases}$$ on
$q_{\hat{\eps}}^{-1}(U^+)\cup q_{\hat{\eps}}^{-1}(U^-)$.

For $m\in M_{\hat{\eps}}$, we define
\begin{equation}\chi_{\hat{\eps}}(m)=z^+\Theta_{\hat{\eps}}^++z^-\Theta_{\hat{\eps}}^-,\label{eq_def_chi}\end{equation}
where $m$ has coordinates $z^+\in U_{\hat{\eps}}^+$ and/or $z^-\in
U_{\hat{\eps}}^-$.

In this way we realize (via $\chi_{\hat{\eps}}$) $M_{\hat{\eps}}$
as a neighborhood of the origin, punctured at $\pm\sqrt{\hat{\eps}}$.
However, the conformal structure of $M_{\hat{\eps}}$ is not preserved,
but is rather expressed by the Beltrami differential
$\mu_{\hat{\eps}}=\frac{\partial \chi_{\hat{\eps}}/\partial
\ov{z}^+}{\partial \chi_{\hat{\eps}}/\partial z^+}$.
We want to show that there exists $K\in(0,1)$ such that $|\mu_{\hat{\eps}}|<K$.
We can then correct the map $\chi_{\hat{\eps}}$ to a conformal map via the
Ahlfors-Bers theorem.

We shall only study what happens on $U_{\hat{\eps}}^{0,\infty}$ as
$\mu\equiv0$ outside these sectors.  We rewrite:
$$\begin{array}{lll} \chi_{\hat{\eps}}(z_+)&=&
z^+(\Theta_{\hat{\eps}}^++\Theta_{\hat{\eps}}^-)+(z^--z^+)\Theta_{\hat{\eps}}^-\\
&=&z^++(z^--z^+)\Theta_{\hat{\eps}}^-\\
&=&z^++\xi_{\hat{\eps}}^{0,\infty}(z^+)\Theta_{\hat{\eps}}^-.\end{array}$$
Then
$$\begin{array}{lll}
\frac{\partial \chi_{\hat{\eps}}}{\partial
z^+}&=&1+\left[\Theta_{\hat{\eps}}^-\frac{\partial
\xi_{\hat{\eps}}^{0,\infty}}{\partial z^+}+
\xi_{\hat{\eps}}^{0,\infty}\frac{\partial
\Theta_{\hat{\eps}}^-}{\partial z^+}\right],\\
\frac{\partial \chi_{\hat{\eps}}}{\partial
\ov{z}^+}&=&\xi_{\hat{\eps}}^{0,\infty}\frac{\partial
\Theta_{\hat{\eps}}^-}{\partial \ov{z}^+}.\end{array}$$

The derivatives of $\Theta_\eps^\pm$ satisfy (for
$z^\pm=x+iy$) near $\pm\sqrt{\hat{\eps}}$:
\begin{equation}\left|\frac{\partial^{n_1+n_2} \Theta_\eps^{\pm}}{\partial x^{n_1}\partial
y^{n_2}}\right|\leq
K_n'\left|z^\pm\pm\sqrt{\hat{\eps}}\right|^{-\gamma
(n_1+n_2)}\label{bound_derivatives}\end{equation}  for some positive
constant $\gamma>0$. Indeed, the estimate \eqref{bound_derivatives}
comes directly from the fact that the derivatives of $\varphi$ are
uniformly bounded by \eqref{varphibound} and that
$(q_{\hat{\eps}}^{-1})'(z)=\frac{1+a(\eps)z}{z^2-\eps}$.

We start by considering $\eps=0$. It is known that
$\xi_0^{0,\infty}$ is exponentially flat in $z^+$ (see for instance
\cite{I}, but the argument is similar to the argument below for the
case $\hat{\eps}\neq0$).

We choose $r>1$ sufficiently small so that we have for $|z^+|<r$
$$\begin{cases}
\frac{\partial \chi_0}{\partial \ov{z}^+}<\frac18,\\
\frac{\partial \chi_0}{\partial z^+}>\frac78.
\end{cases}$$
Using the continuity in $\hat{\eps}$ and estimates on
$\xi_{\hat{\eps}}^{0,\infty}$, to be proved below, we will choose
$\rho'>0$ sufficiently small so that for $|\hat{\eps}|<\rho'$ we
have
$$\begin{cases}
\frac{\partial \chi_{\hat{\eps}}}{\partial \ov{z}^+}<\frac14,\\
\frac{\partial \chi_{\hat{\eps}}}{\partial z^+}>\frac34.
\end{cases}$$

For that we need to bound the functions
$\xi_{\hat{\eps}}^{0,\infty}= \Xi_{\hat{\eps}}^{0,\infty}-id$ and
their derivatives. We use the fact that the
$\Xi_{\hat{\eps}}^{0,\infty}$ are conjugate to
$\Psi_{\hat{\eps}}^{0,\infty}$ through $q_{\hat{\eps}}$. Of course
$Y_0$ can be chosen so that
$\left|\Psi_{\hat{\eps}}^{0,\infty}-id\right|$ is uniformly bounded.
Moreover we have that
$$\xi_{\hat{\eps}}^{0,\infty}(z)=
v_{\eps}^{(\Psi_{\hat{\eps}}^{0,\infty}-id) (q^{-1}(z))}(z)$$ where
$v_{\eps}^{(\Psi_{\hat{\eps}}^{0,\infty}-id) (q^{-1}(z))}$ is the
flow of $v_\eps$ (see \eqref{mod.3}) for the time
$(\Psi_{\hat{\eps}}^{0,\infty}-id) (q^{-1}(z))$, which is uniformly
bounded. It follows from the theorems on the flow and its dependence
on parameters that $\xi_{\hat{\eps}}^{0,\infty}$ and its derivative
with respect to $z$ are uniformly bounded for $|\hat{\eps}|$
sufficiently small. To show that
$|\mu_{\hat{\eps}}|=\left|\frac{\partial \chi_{\hat{\eps}}/\partial
\ov{z}^+}{\partial \chi_{\hat{\eps}}/\partial z^+}\right|<K<1$ for
$|z|<r$ and $|\hat{\eps}|<\rho$ we need to ensure that the
derivatives of $\xi_{\hat{\eps}}^{0,\infty}$ are sufficiently flat
at $\pm \sqrt{\eps}$. So we will  show that
\begin{equation}\label{estimate99}
\left|\xi_{\hat{\eps}}^{0,\infty}(z)\right|< C(\hat{\eps})
\left|z\mp\sqrt{\hat{\eps}} \right|^{\frac{A}{|\sqrt{\hat{\eps}}|}},
\end{equation}
holds for the values $z\in
U^{0,\infty}$ which correspond to values $W=q_{\hat{\eps}}^{-1}(z)$ in the
slanted part of the intersection of the strips.  Here $A$ is a positive
constant which is independent of $\hat{\eps}$.

We will prove \eqref{estimate99} for $\xi_{\hat{\eps}}^0$, the case
$\xi_{\hat{\eps}}^\infty$ being similar and only sketched.

In the slanted part of the intersection of
the strips we have $Im W< -\frac{B}{|2\sqrt{\hat{\eps}}|}$ for some
$B>0$. This yields for the corresponding part of $U^0$
\begin{equation}\label{norm_q}
\left|\exp(-2\pi i q^{-1}(z))\right|<e^{\frac{-2\pi
B}{2|\sqrt{\hat{\eps}}|}}
\end{equation}
when $Im(q^{-1}(z))< -\frac{B}{|2\sqrt{\hat{\eps}}|}$. Let
\begin{equation} g(z)=(z-\sqrt{\hat{\eps}})^{\frac{-2\pi i
(1+a\sqrt{\hat{\eps}})}{2\sqrt{\hat{\eps}}}}
(z+\sqrt{\hat{\eps}})^{\frac{2\pi i
(1-a\sqrt{\hat{\eps}})}{2\sqrt{\hat{\eps}}}}.
\label{eq_g}\end{equation} Then $g(z)=\exp(-2\pi i q^{-1}(z))$ and
$\left|g(z)\right|= \exp(2\pi Im(q^{-1}(z)))<e^{\frac{-2\pi
B}{2|\sqrt{\hat{\eps}}|}}$.

We have $$\Psi_{\hat{\eps}}^0\circ
 q_{\hat{\eps}}^{-1}(z) = q_{\hat{\eps}}^{-1}(z)+\sum_{n\leq -1}b_n g(z)^{-n},$$ yielding that
$$\left|\Psi_{\hat{\eps}}^0\circ
 q_{\hat{\eps}}^{-1}(z)-
 q_{\hat{\eps}}^{-1}(z)\right|=O\left(\left|g(z)\right|\right)= O\left(\left|\exp(-2\pi iq_{\hat{\eps}}^{-1}
 (z))\right|\right) . $$
 Since
 \begin{equation}\frac{dq}{dz}=\frac{1+az}{z^2-\eps},\label{der:q}\end{equation}
if we join two points $z_1$ and $z_2$ in
 the neighborhood of $-\sqrt{\hat{\eps}}$ by a path $\gamma(t)$, $t\in[0,1]$, of length bounded by $c|z_1-z_2|$
 for some $c>0$,
  so that $|\gamma(t)|\geq\min(|z_1+\sqrt{\hat{\eps}}|,|z_2+\sqrt{\hat{\eps}}|)$ for all $t\in [0,1]$, then
  $\left|q(z_1)-q(z_2)\right|\leq c|z_1-z_2|\,\max_{t\in [0,1]}\left|\frac{dq}{dz}(\gamma(t))\right|$.  It follows that
 $$\left|q_{\hat{\eps}}\circ\Psi_{\hat{\eps}}^0\circ
 q_{\hat{\eps}}^{-1}(z)-z\right|=O(|g(z)|). $$
 This holds uniformly in all the region because of
 Lemma~\ref{lemma_estimates} and the constant $C$ is an upper bound for $\left|(z-\sqrt{\hat{\eps}})^{\frac{-2\pi i
(1+a\sqrt{\hat{\eps}})}{2\sqrt{\hat{\eps}}}}\right|$ in the region
corresponding to the slanted part of the strip near
$-\sqrt{\hat{\eps}}$.

In the same way it is possible, using the chain rule,  to show
that the derivatives of $\xi_{\hat{\eps}}^0$ at $-\sqrt{\hat{\eps}}$
remain bounded when $\hat{\eps}\to 0$. Indeed for the derivatives of
$q$ or $q^{-1}$ we use \eqref{der:q}, while for the derivatives of
$\Psi_{\hat{\eps}}^0$ we use \eqref{norm_q}.

 For the case of $\xi_{\hat{\eps}}^\infty$ defined in
\eqref{def_chi}, the only difference with the previous one is the
presence of the translation term in $\Psi_{\hat{\eps}}^\infty$,
which comes from the comparison between the two maps,
$q_{\hat{\eps},\pm}^{-1}$, on $U_{\hat{\eps}}^\infty$. Indeed the map
$q_{\hat{\eps}}^{-1}$ corresponds to  the time for the vector field
\eqref{mod.3}. We have two times $q_{\hat{\eps},\pm}^{-1}$ defined
respectively over $U_{\hat{\eps}}^\pm$. While they can be chosen to
coincide on $U_{\hat{\eps}}^0$  we have that
$q_{\hat{\eps},-}^{-1}=q_{\hat{\eps},+}^{-1} -2\pi i a$ over
$U_{\hat{\eps}}^\infty$. Then the gluing corresponds to
$$q_{\hat{\eps},-}^{-1}=\Psi_{\hat{\eps}}^\infty (q_{\hat{\eps},+}^{-1})+2\pi i a =q_{\hat{\eps},+}^{-1}+\sum_{n\geq
1}c_n\exp(2\pi i nq_{\hat{\eps},+}^{-1}) $$ (see also
\eqref{qm_qp}). The rest of the argument is as in the case of
$\xi_{\hat{\eps}}^0$.

\medskip

Hence $\mu_{\hat{\eps}}$ is a Beltrami field which we extend by
$\mu_{\hat{\eps}}(\pm\sqrt{\hat{\eps}})=0$ in a $C^1$ way. By the
Ahlfors-Bers theorem there exists a 1-1 map $\sigma_{\hat{\eps}}:
\chi_{\hat{\eps}}(M_{\hat{\eps}})\rightarrow \mathbb C$ which is
holomorphic in the sense of this structure and whose image is the
disk $r\mathbb D$. Since this construction is continuous in
$\hat{\eps}$ up to the limit $\eps=0$, we can always suppose that
the boundary point $r$ of $M_{\hat{\eps}}$ is sent to the boundary
point $r$ of $r\mathbb D$ by the composition
$\sigma_{\hat{\eps}}\circ \chi_{\hat{\eps}}$. Then
\begin{equation}\zeta_{\hat{\eps}}=\sigma_{\hat{\eps}}\circ
\chi_{\hat{\eps}}\label{function_zeta}\end{equation} is holomorphic,
yielding that the manifold $M_{\hat{\eps}}$ is conformally
equivalent to the disk $r\mathbb D$ punctured in two points:
$\mathbb D\setminus\{x_1,x_2\}$. We conjugate with the unique
M\"obius transformation $\tau_{\hat{\eps}}$ sending $x_1$, $x_2$ and
$r$ respectively on $-\sqrt{\hat{\eps}}$, $\sqrt{\hat{\eps}}$ and
$r$. The image of $r\mathbb D$ is a disk $D_{\hat{\eps}}$ not
necessarily centered at the origin and whose boundary contains
$\{r\}$. Let us now consider the case $\eps=0$: there exists a
one-parameter family of M\"obius transformations $\tau$ sending the
double point $x_1=x_2$ and $r$ to $0$ and $r$ respectively. Each one
is uniquely determined by the derivative at $x_1$. We choose the one
such that $\zeta_0'(0)\tau'(0)=1$. Indeed we have
$$\lim_{\hat{\eps}\to 0}\frac{(\tau_{\hat{\eps}}\circ\zeta_{\hat{\eps}})
(\sqrt{\hat{\eps}})-(\tau_{\hat{\eps}}\circ\zeta_{\hat{\eps}})
(-\sqrt{\hat{\eps}})}{2\sqrt{\hat{\eps}}}\equiv1.$$

The construction of $\mu_{\hat{\eps}}$ is continuous in $\hat{\eps}$
and has a limit when $\hat{\eps}\to 0$ on radial rays, yielding the
same property for the construction above. We will show below how to
modify it slightly so as to ensure that it is also holomorphic in
$\hat{\eps}\neq0$ and with a uniform limit on all rays.

Let us start by looking at the different limits we get for $\eps=0$
along the different rays $\arg(\hat{\eps})=Const$. When constructing
an abstract manifold by charts and transition maps between charts,
the size of the charts is not intrinsic and it is possible to modify
them as long as the new transition maps are analytic extensions of
the previous ones. So we get different presentations of a unique
manifold as long as the total underlying set is the same. We must be
careful at the boundary. Indeed the outer boundary of
$U_{\hat{\eps}}^+$ is not in general sent into the outer boundary of
$U_{\hat{\eps}}^-$ under the gluing map. This is why we have taken
so much care so that the intersection of the strips be constant near
the boundary of the hole in $W$-space (see Figure~\ref{strips}).
With this property the limit is independent of $\arg(\hat{\eps})$
since the different $\xi_0^{0,\infty}$ obtained with different
slopes are all analytic extensions one of the other.

Let us now show that the map $f_{\hat{\eps}}$ depends analytically
on $\hat{\eps}$. We start by considering a small sector
$\arg{\hat{\eps}}\in (\theta_0-\eta,\theta_0+\eta)$ for some fixed
$\theta_0$ and some small $\eta$. It is possible over such a sector
to reproduce the same construction as above, but with strips having
a fixed slope (for instance that chosen for
$\arg{\hat{\eps}}=\theta_0$) and a fixed intersection domain. Since
the intersection is fixed, it is possible to choose a fixed
$N_{\hat{\eps}}$ in \eqref{eq_N}, hence depending analytically on
$\hat{\eps}$. In this way we locally get  maps $\sigma_{\hat{\eps}}$
and $\zeta_{\hat{\eps}}$ which are analytic in $\hat{\eps}$. But
these maps have just been instrumental in constructing a unique disk
$D_{\hat{\eps}}$ endowed with a unique map $f_{\hat{\eps}}$. It
follows that $f_{\hat{\eps}}$ depends analytically on $\hat{\eps}$.
The analytic dependence on the auxiliary multi-parameter $\nu$, is
an immediate application of the analytic dependence on parameters in
the Ahlfors-Bers theorem. \hfill $\Box$

\section{The compatibility condition}\label{sect:compatibility}

In Section~\ref{sect:local} we have realized the modulus
$(a(\eps),[\psi_{\hat{\eps}}^0,\psi_{\hat{\eps}}^{\infty}])_{\hat{\eps}\in
V_{\rho,\delta}}$ in a family $f_{\hat{\eps}}$ which is ramified in
$\hat{\eps}$ over some sectorial neighborhood $V_{\rho,\delta}$. We
are now interested in the condition that the family
$(a(\eps),[\psi_{\hat{\eps}}^0,\psi_{\hat{\eps}}^{\infty}])_{\hat{\eps}\in
V_{\rho,\delta}}$ must satisfy in order that there exists a
realization in a uniform family $f_\eps$ defined for $\eps\in
B(0,\rho)$.

We limit our discussion to the sector \begin{equation}V_G=
V_G(\rho)=\{\eps;0<|\eps|<\rho, \arg \eps\in
(-\delta,\delta)\},\label{V_eta}\end{equation} which is covered in
$V_{\rho,\delta}$ by two small sectors
\begin{equation}\begin{cases} \widetilde{V}=\{\hat{\eps};0<|\hat{\eps}|<\rho,
\arg\hat{\eps}\in(2\pi -
\delta, 2\pi +\delta)\}\\
\ov{V}=\{\hat{\eps};0<|\hat{\eps}|<\rho, \arg\hat{\eps}\in( -
\delta, +\delta)\}.\end{cases}\label{secteur_V_eta}\end{equation} We
remark that the Glutsyuk modulus exists for $\eps\in V_G$, and $\rho$ sufficiently
small. Depending
on the context and whether we want to concentrate on $\rho$ or not
we will use either the notation $V_G$ or $V_G(\rho)$.

A necessary condition for the existence of a uniform realization is
that the functions $f_{\hat{\eps}}$ and $f_{\hat{\eps}e^{2\pi i}}$
be conjugate. In order to simplify the notation we will write
\begin{equation}\begin{cases}
\ov{\eps}= \hat{\eps}, &\qquad\hat{\eps}\in \ov{V},\\
\tilde{\eps}=\hat{\eps} e^{2\pi i}, &\qquad\hat{\eps} e^{2\pi i}\in
\widetilde{V}.\end{cases}\end{equation} Hence $\ov{\eps}$ and
$\tilde{\eps}$ project on the same $\eps\in V_G$.  These functions
have their moduli presented in different ways. We need to find a
compatibility condition (in terms of the modulus) which expresses
the fact that the two presentations encode the same dynamics up to
conjugacy.

In order to investigate this further we use the notation
$\ov{\Psi}^{0,\infty}$ and
$\ov{\Xi}^{0,\infty}=id+\ov{\xi}^{0,\infty}$ when $\hat{\eps}\in
\ov{V}$ and $\widetilde{\Psi}^{0,\infty}$ and
$\widetilde{\Xi}^{0,\infty}= id +\tilde{\xi}^{0,\infty}$ when
$\hat{\eps}\in \widetilde{V}$. We work for a fixed value
$\ov{\eps}=\hat{\eps}\in \ov{V}$ and the corresponding
$\tilde{\eps}=\hat{\eps}e^{2\pi i}\in \widetilde{V}$. Because we
work with two fixed values of $\hat{\eps}$ we will omit mentioning
these values in the indices.
\begin{figure}[!h]
\begin{center}
\subfigure[$\hat{\eps}\in\widetilde{V}$
]{\includegraphics[width=4.5cm]{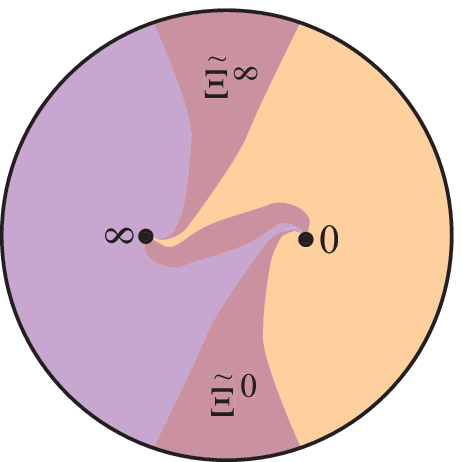}}\qquad
\subfigure[The Glutsyuk
sectors]{\includegraphics[width=4.5cm]{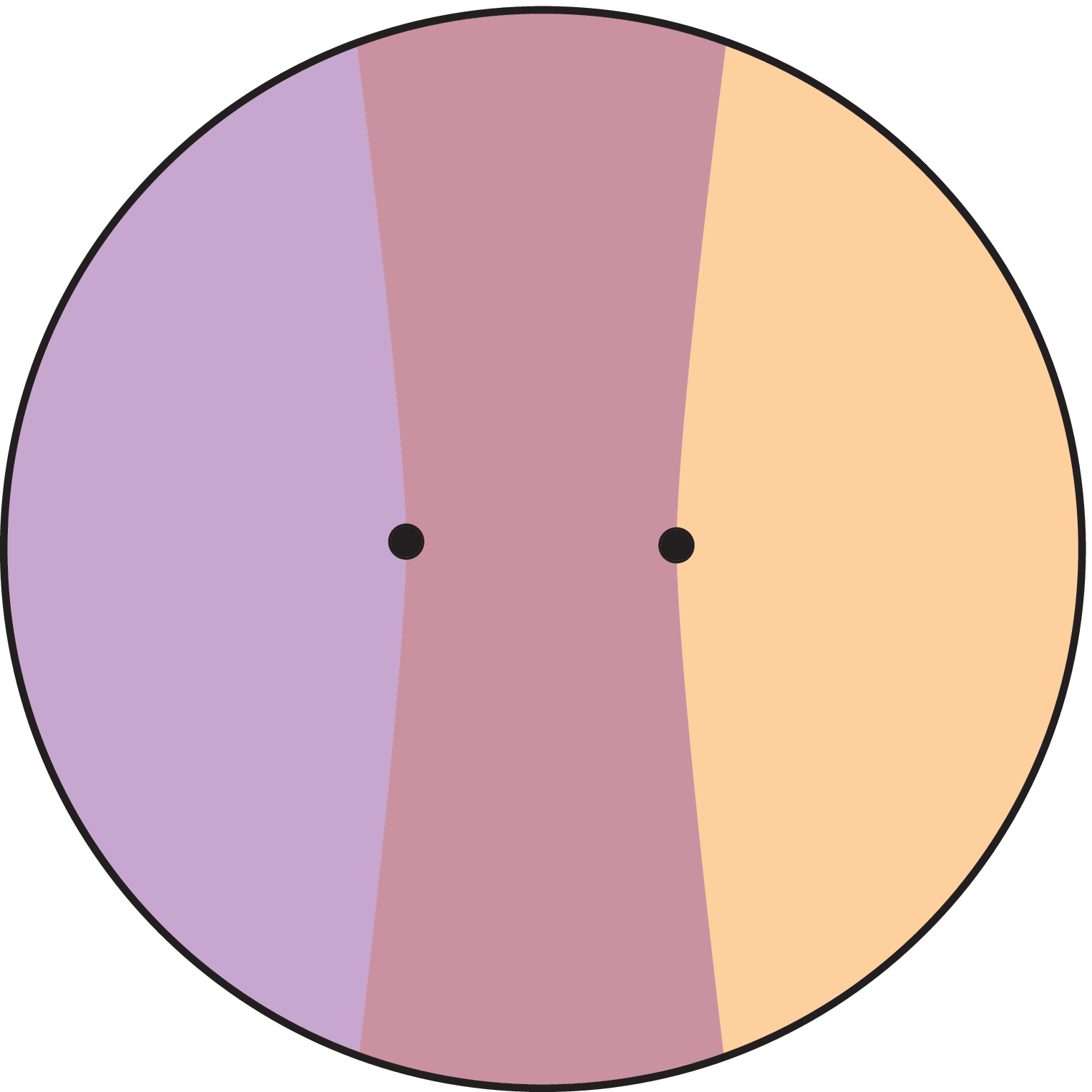}}\qquad
\subfigure[$\hat{\eps}\in\ov{V}$]{\includegraphics[width=4.5cm]{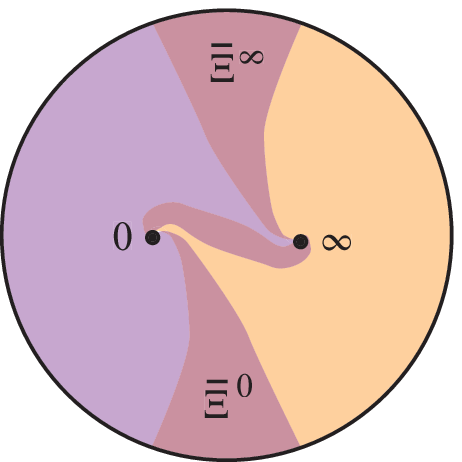}}
\caption{The different sectors} \label{compatibility}
\end{center}
\end{figure}
In the point of view corresponding to $\ov{V}$, the left (resp.
right) singular point is $-\sqrt{\hat{\eps}}$ (resp.
$\sqrt{\hat{\eps}}$) and $\ov{\Xi}^0$ (resp. $\ov{\Xi}^\infty$)
describes the gluing when turning around it. In the point of view
corresponding to $\widetilde{V}$, the left (resp. right) singular
point is $\sqrt{\hat{\eps}}$ (resp. $-\sqrt{\hat{\eps}}$) and
$\widetilde{\Xi}^\infty$ (resp. $\widetilde{\Xi}^0$) describes the
gluing when turning around it. Remark that in all cases $\infty$
(resp. $0$) will represent $\sqrt{\hat{\eps}}$ (resp.
$-\sqrt{\hat{\eps}}$).

\begin{figure}[!h]
\begin{center}
\includegraphics[width=11cm]{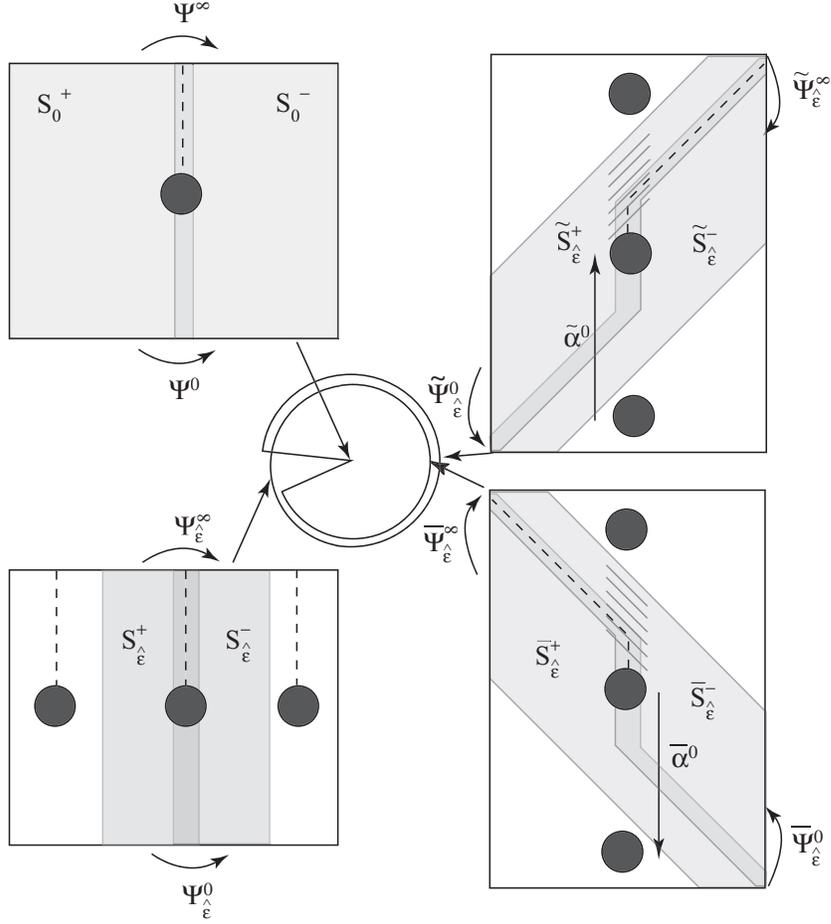}
\caption{The darkened (striped) region above the principal hole
where we compare the two points of view.}
\label{strips_adjusted_compare}
\end{center}
\end{figure}

The idea is to derive the Glutsyuk modulus from these two Lavaurs
moduli and to equate them. This is done in considering the darkened
(striped) regions of the two pictures on the right in
Figure~\ref{strips_adjusted_compare}.

We define the following quantities related to the periods of
$q_{\hat{\eps}}^\pm$ near the inverse images of $\pm
\sqrt{\hat{\eps}}$.

\begin{equation}\begin{cases}
\alpha^\infty=-\frac{2\pi i
(1+a(\eps)\sqrt{\hat{\eps}})}{2\sqrt{\hat{\eps}}}=-\frac{2\pi
i}{\mu^\infty},\\
\alpha^0=-\frac{2\pi i
(1-a(\eps)\sqrt{\hat{\eps}})}{2\sqrt{\hat{\eps}}}=\frac{2\pi
i}{\mu^0}.\end{cases}\label{eq_alpha1}
\end{equation} Hence
\begin{equation}\alpha^\infty=\alpha^0-2\pi i a. \label{eq_alpha2}
\end{equation} We will have
$\tilde{\alpha}^0$ and $\tilde{\alpha}^\infty$ over $\widetilde{V}$
and $\ov{\alpha}^0$ and $\ov{\alpha}^\infty$ over $\ov{V}$. Moreover
\begin{equation}\begin{cases}
\tilde{\alpha}^0=-\ov{\alpha}^\infty,\\
\tilde{\alpha}^\infty=-\ov{\alpha}^0.\end{cases}\label{eq_alpha3}
\end{equation} We define
\begin{equation}\begin{cases}
\widetilde{C}^{0,\infty}=\exp(-2\pi i\tilde{\alpha}^{0,\infty}),\\
\ov{C}^{0,\infty}=\exp(-2\pi
i\ov{\alpha}^{0,\infty}).\end{cases}\label{eq_C}
\end{equation}
In particular $\widetilde{C}^{0,\infty}=\exp(-2\pi
i\tilde{\alpha}^{0,\infty})$  are exponentially large in
$\sqrt{\hat{\eps}}$ while $\ov{C}^{0,\infty}=\exp(-2\pi i
\ov{\alpha}^{0,\infty})$ are exponentially small in
$\sqrt{\hat{\eps}}$.

\begin{theorem}\label{thm:compatibility}
\begin{description}
\item{(i)} There exists $Y_1>0$ such that for all $\hat{\eps}\in \widetilde{V}_\eta$ there
exists a map $\widetilde{H}^0$ defined in a region $Im(W)<-Y_1$,
commuting with $T_1$, and such that
\begin{equation}\widetilde{H}^0 \circ T_{\tilde{\alpha}^0}\circ
\widetilde{\Psi}^0 =
T_{\tilde{\alpha}^0}\circ\widetilde{H}^0.\label{equation_H_0}\end{equation}
  In the new coordinate $\widetilde{W}^0=\widetilde{H}^0(W)$ the
renormalized return map $T_{\tilde{\alpha}^0}\circ
\widetilde{\Psi}^0$ is a translation. Then $\widetilde{W}^0$ is one
Fatou Glutsyuk coordinate. Similarly there exists a map
$\widetilde{H}^\infty$ defined in the region $Im(W)>Y_1$, commuting
with $T_1$, and such that
\begin{equation}\widetilde{H}^\infty\circ T_{\tilde{\alpha}^0}\circ
\widetilde{\Psi}^\infty=
T_{\tilde{\alpha}^\infty}\circ\widetilde{H}^\infty.\label{equation_H_infty}\end{equation}
In the new coordinate
$\widetilde{W}^\infty=\widetilde{H}^\infty(W)$ the renormalized
return map is a translation and $\widetilde{W}^\infty$ is the second
Fatou Glutsyuk coordinate. The Glutsyuk modulus is then given by
$\widetilde{H}^\infty \circ (\widetilde{H}^0)^{-1}$. \item{(ii)}
Similarly there exists $Y_2>0$ such that for all
$\ov{\eps}\in\ov{V}_\eta$ there exists $\ov{H}^{0,\infty}$ commuting
with $T_1$ and such that
\begin{equation} \ov{H}^0 \circ \ov{\Psi}^0\circ T_{\ov{\alpha}^0}=
T_{\ov{\alpha}^0}\circ\ov{H}^0\label{equation_H0_ov}\end{equation}
on $Im(W)<-Y_2$ and \begin{equation}\ov{H}^\infty \circ
\ov{\Psi}^\infty\circ T_{\ov{\alpha}^0}=
T_{\ov{\alpha}^\infty}\circ\ov{H}^\infty\label{equation_Hinfty_ov}\end{equation}
on $Im(W)>Y_2$. The Glusyuk modulus is then given in this context by
$\ov{H}^0 \circ (\ov{H}^\infty)^{-1}$. Considering (i) and (ii)
together we can of course suppose that $Y_1=Y_2$.
\item{(iii)} The maps $\ov{H}^{0,\infty}$ and $\widetilde{H}^{0,\infty}$ are  unique up to left
composition with a translation. In particular they are unique if we
ask that their limits for $Im W\to \pm \infty$ be the identity.
\item{(iv)} The functions $\widetilde{H}^0$ and $\widetilde{H}^\infty$ (resp. $\ov{H}^0$ and $\ov{H}^\infty$) have
analytic extensions defined on domains which intersect.
\item{(v)} A necessary condition for the family
$(\Psi_{\hat{\eps}}^0,\Psi_{\hat{\eps}}^\infty)$ to be the modulus
of an analytic family $f_\eps$ of diffeomorphisms is that for
corresponding values of $\ov{\eps}\in \ov{V}_\eta$ and
$\tilde{\eps}\in \widetilde{V}_\eta$
 there exist constants $D_{\eps}$ and $D'_{\eps}$ (depending on $\eps$, not on $\hat{\eps}$!) such that
\begin{equation}\widetilde{H}^\infty \circ (\widetilde{H}^0)^{-1}=T_{D_{\eps}}\circ\ov{H}^0
\circ (\ov{H}^\infty)^{-1}\circ
T_{D'_{\eps}}.\label{necessary_condition}\end{equation} This
condition is called the {\em compatibility condition}.
\item{(vi)} {}The functions $\widetilde{H}^{0,\infty}$ and $\ov{H}^{0,\infty}$ can be chosen to depend
analytically on the auxiliary multi-parameter $\nu$, as can the
constants $D$ and $D'$ in \eqref{necessary_condition}.
\end{description}\end{theorem} \noindent {\scshape Proof.}
\begin{description} \item{(i) and (ii)} Conjugating
$T_{\tilde{\alpha}^0}\circ \widetilde{\Psi}^{0,\infty}$ under
$E(W)=\exp(-2\pi i W)$ yields maps $\tilde{\kappa}^{0,\infty}=E\circ
T_{\tilde{\alpha}^0}\circ \widetilde{\Psi}^{0,\infty}\circ E^{-1}$
with multiplier of modulus different from one.  Hence
$\tilde{\kappa}^0$ (resp. $\tilde{\kappa}^\infty$) is linearizable
in the neighborhood of $0$ (resp. $\infty$): there exists
$\tilde{h}^{0,\infty}$ such that
$$\tilde{h}^{0,\infty}\circ \tilde{\kappa}^{0,\infty}= L_{\exp(-2\pi i
\tilde{\alpha}^{0,\infty})}\circ \tilde{h}^{0,\infty}.$$ The maps
$\widetilde{H}^{0,\infty}$ are simply $E^{-1}\circ
\tilde{h}^{0,\infty}\circ E$. It then follows that they commute with
$T_1$.

The existence of $\ov{H}^{0,\infty}$ is similar.

\item{(iii)}  This is obvious since $\ov{h}^{0,\infty}$ and $\tilde{h}^{0,\infty}$ are
unique up to left composition with linear maps.
\item{(iv)} The relation \eqref{equation_H_0} allows to
extend $\widetilde{H}^0$ by means of $\widetilde{H}^0\circ
T_{\tilde{\alpha}^0}=T_{\tilde{\alpha}^0}\circ \widetilde{H}^0\circ
(\widetilde{\Psi}^0)^{-1}\circ T_{-\tilde{\alpha}^0}$, so its domain
becomes the image of $\widetilde{\Psi}^0$ augmented of a strip of
width $\tilde{\alpha}^0$. Similarly for $\widetilde{H}^\infty$,
$\ov{H}^0$ and $\ov{H}^\infty$.

We claim the existence of a uniform domain. The intuitive idea is
that there are no recurrent points for $f_{\hat{\eps}}$ for these
values of $\hat{\eps}$. In practice, the relations
\eqref{equation_H_0}, \eqref{equation_H_infty},
\eqref{equation_H0_ov} and \eqref{equation_Hinfty_ov} allow to
extend the maps in the direction of $\alpha$. The fact that the maps
commute with $T_1$ allows to extend them until the holes. Hence the
claim.

\item{(v)} The compatibility condition comes from the fact that each
Fatou Glutsyuk coordinate is uniquely determined up to a
translation.

\item{(vi)} This is clear from the nature of the proofs above.
\hfill $\Box$\end{description}

The compatibility condition was found independently by Reinhard
Sch\"afke \cite{Sc} in the case $\Psi_{\hat{\eps}}^\infty\equiv id$.

\begin{remark} For each translation $T_A$ there exists a unique $\widetilde{H}^0$,
such that $\lim_{Im(W)\to -\infty} \widetilde{H}^0 = W+A$. Similar
statements are valid for $\widetilde{H}^\infty$, $\ov{H}^0$,
$\ov{H}^\infty$. \end{remark}

\begin{proposition}\label{prop_constants}
We consider the modulus
$\left.\left(a(\eps),[\Psi_{\hat{\eps}}^0,\Psi_{\hat{\eps}}^\infty]\right)\right|_{\hat{\eps}\in
V_{\rho,\delta}}$ attached to a germ of one-para\-me\-ter prepared
analytic family of diffeomorphisms of the form
\eqref{prepared_k_par} and hence satisfying the compatibility
condition \eqref{necessary_condition}. Then there exists an analytic
function defined by $\hat{\eps}\mapsto \gamma_{\hat{\eps}}$ on
$V_{\rho,\delta}$, such that on $V_\eta$ we have
$$\gamma_{\tilde{\eps}}-\gamma_{\ov{\eps}}= D_{\eps}- 2 \pi i a$$
and $$\lim_{\hat{\eps}\to 0}\gamma_{\hat{\eps}}= \gamma_0,$$
for some constant $\gamma_0$.
\end{proposition}

\begin{corollary}\label{remark_D} Given a modulus $\left.\left(a(\eps),[\Psi_{\hat{\eps}}^0,\Psi_{\hat{\eps}}^\infty]\right)\right|_{\hat{\eps}\in
V_{\rho,\delta}}$, it is possible to choose a representative
$\left.(a(\eps),
\Psi_{\hat{\eps}}^0,\Psi_{\hat{\eps}}^\infty)\right|_{\hat{\eps}\in
V_{\rho,\delta}}$
  so that that $D_\eps\equiv 2\pi i a$ in
\eqref{necessary_condition} and
\begin{equation} D'_\eps=-2\pi i a +O(\exp(-2\pi i
\ov{\alpha}^0)).\label{eq_G} \end{equation}
\end{corollary}
\noindent{\scshape Proof.} It is possible to represent the modulus
by the family
\begin{equation}(\Upsilon_{\hat{\eps}}^0,\Upsilon_{\hat{\eps}}^\infty)=(T_{-\gamma_{\hat{\eps}}}\circ
\Psi_{\hat{\eps}}^0\circ
T_{\gamma_{\hat{\eps}}},T_{-\gamma_{\hat{\eps}}}\circ
\Psi_{\hat{\eps}}^\infty\circ
T_{\gamma_{\hat{\eps}}}).\label{new_modulus}\end{equation} In the
equation \eqref{necessary_condition} the maps $\widetilde{H}^0$,
$\widetilde{H}^\infty$, $\ov{H}^0$, $\ov{H}^\infty$ are then
replaced by \begin{equation}\begin{cases}
\widetilde{H}_1^{0,\infty}=T_{-\gamma_{\tilde{\eps}}}\circ
\widetilde{H}^{0,\infty}\circ T_{\gamma_{\tilde{\eps}}},\\
\ov{H}_1^{0,\infty}=T_{-\gamma_{\ov{\eps}}}\circ
\ov{H}^{0,\infty}\circ T_{\gamma_{\ov{\eps}}}.
\end{cases}\label{new_H}\end{equation} They satisfy the
compatibility condition
\begin{equation}\widetilde{H}_1^\infty \circ (\widetilde{H}_1^0)^{-1}=T_{2\pi i a}\circ\ov{H}_1^0
\circ (\ov{H}_1^\infty)^{-1}\circ
T_{D''}.\label{eq_comp2}\end{equation} We postpone the proof that
\begin{equation} D''=-2\pi i a +O(\exp(-2\pi i
\ov{\alpha}^0))\label{eq_G1} \end{equation}  after the proof of
Lemma~\ref{lemma_flatness}. \hfill $\Box$

\bigskip
 \noindent{\scshape
Proof of Proposition~\ref{prop_constants}.} When we define Fatou
coordinates we have one degree of freedom per Fatou coordinate. One
degree of freedom has been used when we asked that $\lim_{Im(W)\to
-\infty}\Psi_{\hat{\eps}}^0= id$, the other degree of freedom can
be used to fix a base point for the Fatou coordinate $\Phi_{\hat{\eps}}^-$.
Consider Figure~\ref{strips}: we can choose a base point $Z_0$ located on
the right of the principal hole and we can choose the Fatou
coordinate $\Phi_{\hat{\eps}}^-$ such that
$\Phi_{\hat{\eps}}^-(Z_0)=Z_0$. This is done via the composition
$T_{-\gamma_{\hat{\eps}}}\circ \Phi_{\hat{\eps}}^-$. Then $\Phi^+$
is completely determined by $\lim_{Im(W)\to
-\infty}\Psi_{\hat{\eps}}^0= id$. This yields the new
representative of the modulus in \eqref{new_modulus}.

Once the Lavaurs Fatou coordinates are chosen, the Fatou Glutsyuk
coordinates are completely determined by the limit conditions on the
functions $\ov{H}^{0,\infty}$ and $\widetilde{H}^{0,\infty}$. So for
the new Fatou coordinate and representative \eqref{new_modulus}, the
new  Fatou Glutsyuk coordinates are simply given in \eqref{new_H}
(i.e. by $\widetilde{H}_1^{0,\infty}(W)$ and
$\widetilde{H}_1^{0,\infty}(W)$). At the limit when $\eps=0$, the
Fatou Lavaurs and Fatou Glutsyuk coordinates coincide.

 The only thing we need to take care of is
that the darkened regions of Figure~\ref{strips} lie in different
sheets due to the sweep of the cut as $\eps$ made a full turn.
Indeed when we
adjust the constant $D$ we compare the domains of $\ov{H}^0$ and
$\widetilde{H}^\infty$. $\ov{H}^0$ conjugates $\ov{\Psi}^0\circ
T_{\ov{\alpha}^0}$ to a translation. We have $T_{\ov{\alpha}^0}:
\ov{S}^-\rightarrow \ov{S}^+$, while $\ov{\Psi}^0:
\ov{S}^+\rightarrow \ov{S}^-$. Hence $\ov{\Psi}^0\circ
T_{\ov{\alpha}^0}: \ov{S}^-\rightarrow \ov{S}^-$ and $\ov{H}^0$ is
defined on $\ov{S}^-$. On the other hand $\widetilde{H}^\infty$
conjugates $T_{\tilde{\alpha}^0}\circ \widetilde{\Psi}^\infty:
\widetilde{S}^+\rightarrow \widetilde{S}^+$ to a translation. Hence
$\widetilde{H}^\infty$ is defined on $\widetilde{S}^+$. Because of
the definition of $\ov{S}^\pm$ and $\widetilde{S}^\pm$ the passage
map $\ov{S}^-\rightarrow \widetilde{S}^+$ is $T_{2\pi i a}$.
$\Box$

\begin{lemma}\label{lemma_flatness} We consider the maps
$\widetilde{H}^0$, $\widetilde{H}^\infty$, $\ov{H}^0$,
$\ov{H}^\infty$ of Theorem~\ref{thm:compatibility}. We let
$$\begin{cases}\widetilde{\Psi}^0 = id+\widetilde{\Lambda}^0,&
\widetilde{\Psi}^\infty=T_{-2\pi i a}+\widetilde{\Lambda}^\infty,\\
\ov{\Psi}^0=id+\ov{\Lambda}^0,& \ov{\Psi}^\infty= T_{-2\pi i
a}+\ov{\Lambda}^\infty,\\
\widetilde{H}^{0,\infty}= id +\widetilde{G}^{0,\infty},&
\ov{H}^{0,\infty}= id+\ov{G}^{0,\infty}.\end{cases}$$
\begin{description}
\item{(i)} The functions $\widetilde{G}^{0,\infty}$ are given by the
following series which are absolutely convergent for $|Im
W|>Y_0+2\beta$ (see Lemma~\ref{lemma_estimates}) and $\hat{\eps}\in
\widetilde{V}$
\begin{equation}
\widetilde{G}^0=-\sum_{n=1}^\infty
\widetilde{\Lambda}^0\circ(T_{\tilde{\alpha}^0}\circ
\widetilde{\Psi}^0)^{-n},\label{sol_tilde_G0}\end{equation}
\begin{equation}
\widetilde{G}^\infty=\sum_{n=0}^\infty
\widetilde{\Lambda}^\infty\circ(T_{\tilde{\alpha}^0}\circ
\widetilde{\Psi}^\infty)^n.\label{sol_tilde_Ginfty}\end{equation}
Similarly the functions $\ov{G}^{0,\infty}$ are given by the
following series which are absolutely convergent for $|Im
W|>Y_0+2\beta$ and $\hat{\eps}\in \bar{V}$
\begin{equation}
\ov{G}^0=\sum_{n=0}^\infty\ov{\Lambda}^0\circ
T_{\ov{\alpha}^0}\circ(\ov{\Psi}^0\circ T_{\ov{\alpha}^0}
)^n,\label{sol_bar_G0}\end{equation}
\begin{equation}
\ov{G}^\infty=- \sum_{n=1}^\infty \ov{\Lambda}^\infty\circ
T_{\ov{\alpha}^0}\circ(\ov{\Psi}^\infty\circ T_{\ov{\alpha}^0}
)^{-n}.\label{sol_bar_Ginfty}\end{equation} For $\hat{\eps}\to 0$ we
have the following limits $$\begin{cases} \lim_{\hat{\eps}\to 0}
\widetilde{H}_{\hat{\eps}}^0= \lim_{\hat{\eps}\to 0} \ov{H}_{\hat{\eps}}^0=id,\\
\lim_{\hat{\eps}\to 0} \widetilde{H}_{\hat{\eps}}^\infty=T_{2\pi i a }\circ \Psi_0^\infty,\\
\lim_{\hat{\eps}\to 0} (\ov{H}_{\hat{\eps}}^\infty)^{-1}=
\Psi_0^\infty \circ T_{2\pi i a }.\end{cases}$$
\item{(ii)} For
$\hat{\eps}\in \widetilde{V}$  we have
$$\begin{cases}\widetilde{H}^0=id+O(\ov{C}^0),\\
\widetilde{H}^\infty= \widetilde{\Psi}^\infty +2 \pi i a +
O(\ov{C}^0),
\end{cases}$$ while for $\hat{\eps}\in \bar{V}$ we have
$$\begin{cases}\ov{H}^0=id+O(\ov{C}^0),\\
(\ov{H}^\infty)^{-1}= \ov{\Psi}^\infty \circ T_{2\pi i a} +
O(\ov{C}^0),
\end{cases} $$ where $$\ov{C}^0<\exp\left(-\frac{2\pi(2\pi -
\gamma^*)}{\sqrt{\hat{\eps}}}\right),$$
for some $\gamma^*\in (0, \frac12)$.
\end{description}
\end{lemma}
\noindent{\scshape Proof.} (i) Let us derive \eqref{sol_tilde_G0}.
The function $\widetilde{G}^0$ satisfies $\widetilde{G}^0\circ
T_{\tilde{\alpha}^0}\circ \widetilde{\Psi}^0 =
\widetilde{G}^0-\widetilde{\Lambda}^0$, which we rewrite
\begin{equation}\widetilde{G}^0=\widetilde{G}^0\circ(T_{\tilde{\alpha}^0}\circ \widetilde{\Psi}^0)^{-1}
-\widetilde{\Lambda}^0\circ(T_{\tilde{\alpha}^0}\circ
\widetilde{\Psi}^0)^{-1}.\label{telescopic} \end{equation} We obtain
an infinite set of equations by composing \eqref{telescopic} on the
right with $(T_{\tilde{\alpha}^0}\circ \widetilde{\Psi}^0)^{-n}$.
Adding these equations yields a telescopic sum.  The formula
\eqref{sol_bar_G0} is checked in the same manner. For the formulas
\eqref{sol_tilde_Ginfty}, and \eqref{sol_bar_Ginfty} we also use
\eqref{eq_alpha2}. To prove the convergence we use
Lemma~\ref{lemma_estimates}. Indeed let $\ov{\Psi}^0= id
+\ov{\Lambda}$. Let us look at \eqref{sol_bar_G0}. If $W=X+iY$ and
$Y< -Y_0-2\beta$,  then
$$\left|\ov{\Lambda}^0(W)\right|\leq \ov{M}^0N(\beta)\exp(2\pi(Y_0+\beta+Y))=\ov{N}^0\exp(2\pi Y),$$
where $\ov{N}_0=\ov{M}^0N(\beta)$, and $N(\beta)$ is a positive function as in Lemma~\ref{lemma_estimates}.
For $\arg(\hat{\eps})\in(-\delta,\delta)$, then
$Im(\ov{\alpha}^0)=-\frac{2\pi-
\gamma(\hat{\eps})}{|\sqrt{\hat{\eps}}|}$ for some
$\gamma(\hat{\eps})\in (0,\frac12)$. We can show by induction that
\begin{equation}\left|Im (\ov{\Psi}^0\circ T_{\ov{\alpha}^0}
)^n(W)-Im(T_{\ov{\alpha}^0})^n( W)\right|<n\ov{N}^0\exp\left(2\pi
\left(Y-\frac{n(2\pi-\gamma(\hat{\eps}))}{|\sqrt{\hat{\eps}}|}\right)\right).\label{imaginary}\end{equation}
Hence
$$Im (\ov{\Psi}^0\circ T_{\ov{\alpha}^0}
)^n(W)< Im W - \frac{n(2\pi-\gamma)}{|\sqrt{\hat{\eps}}|}+n\ov{B}^0
$$
for some positive constant $\ov{B}^0$. The convergence of $\ov{G}^0$
follows.

(ii) The fact that $\widetilde{H}=id+O(\ov{C}^0)$ comes from
\eqref{imaginary}.

To derive that $(\ov{H}^\infty)^{-1}= \ov{\Psi}^\infty \circ
T_{2\pi i a}+O(\ov{C}^0)$ we calculate $(\ov{H}^\infty)^{-1}$
directly from
\begin{equation} \ov{\Psi}^\infty\circ T_{\ov{\alpha}^0}\circ
(\ov{H}^\infty)^{-1}= (\ov{H}^\infty)^{-1}\circ
T_{\ov{\alpha}^\infty}.\label{eq_H_infty_inverse}\end{equation}\hfill $\Box$

 \medskip\noindent{\scshape End of Proof of Corollary~\ref{remark_D}.}
 We now need to prove \eqref{eq_G} (i.e. \eqref{eq_G1}).
 This follows from calculation of the constant terms on both sides
 of
 \eqref{eq_comp2}, using the fact that $\ov{H}^0$ and
 $\widetilde{H}^{0}$ are almost the identity. \hfill $\Box$

\begin{remark}It is remarkable that, although the
functions $\widetilde{H}_{\hat{\eps}}^0$,
$\widetilde{H}_{\hat{\eps}}^\infty$, $\ov{H}_{\hat{\eps}}^0$,
$\ov{H}_{\hat{\eps}}^\infty$ have no geometric meaning for
$\hat{\eps}=0$, the limits however exist.
\end{remark}

\begin{theorem}\label{thm_flatness}
We consider a family
$(\Psi_{\hat{\eps}}^0,\Psi_{\hat{\eps}}^\infty)$ for which the
compatibility condition
\begin{equation}\widetilde{H}^\infty \circ (\widetilde{H}^0)^{-1}=T_{2\pi i a} \circ \ov{H}^0
\circ (\ov{H}^\infty)^{-1}\circ
T_{D'}.\label{necessary_condition2}\end{equation} is met for
$\hat{\eps}\in \ov{V}$ and the corresponding $\hat{\eps}e^{2\pi
i}\in \widetilde{V}$ and such that
\begin{equation} D'= -2\pi i a +O(\exp(-2\pi i \ov{\alpha}^0)). \label{eq_D}\end{equation} Then if we
use the notation
$$\begin{cases}
\ov{\Psi}^{0,\infty}=\Psi_{\hat{\eps}}^{0,\infty},\\
\widetilde{\Psi}^{0,\infty}=\Psi_{\hat{\eps}e^{2\pi i}}^{0,\infty},
\end{cases}$$
we have
\begin{equation}
\ov{\Psi}^0-\widetilde{\Psi}^0=O(\exp(-2\pi i
\ov{\alpha}^0))\label{psi0_summable}\end{equation} and
\begin{equation}
\ov{\Psi}^\infty-\widetilde{\Psi}^\infty=O(\exp(-2\pi i
\ov{\alpha}^0)). \label{psi_infty_summable}\end{equation}
\end{theorem} \noindent {\scshape Proof. } We have seen
in the proof of Lemma~\ref{lemma_estimates} that
$$\widetilde{H}^\infty=\widetilde{\Psi}^\infty+2\pi i a +O(\exp(2\pi i
\ov{\alpha}^\infty))$$ and \begin{equation}\ov{H}^0=id+O(\exp(-2\pi
i \ov{\alpha}^0)).\label{estime_ovH_0}\end{equation}
 $(\ov{H}^\infty)^{-1}$ has been
calculated in Lemma~\ref{lemma_flatness}. Since
$\widetilde{H}^0=id+O(\exp(2\pi i \tilde{\alpha}^0))$, we also have
\begin{equation}(\widetilde{H}^0)^{-1} =id+O(\exp(2\pi i
\tilde{\alpha}^0)).\label{estime_h_0_tilde}\end{equation} Replacing
in \eqref{necessary_condition} we show that we get
\eqref{psi_infty_summable} and \eqref{eq_D}.

From the expression of $\widetilde{\Psi}^\infty$ (resp.
$\ov{\Psi}^\infty$) in term of $\widetilde{H}^\infty$ (resp.
$\ov{H}^\infty$) it suffices to show that
$\left|\widetilde{H}^\infty(W)-T_{2\pi i
a}\circ(\ov{H}^\infty)^{-1}\circ T_{-2\pi i
a}(W)\right|=O(\exp(-2\pi i \ov{\alpha}^0))$ follows from
\eqref{necessary_condition2}. Indeed let $\widetilde{W}=
(\widetilde{H}^0)^{-1}(W)$, then
$$\begin{array}{l}
\left|\widetilde{H}^\infty(W)-T_{2\pi i
a}\circ(\ov{H}^\infty)^{-1}\circ T_{-2\pi i a}(W)\right|\\
\qquad\qquad\qquad\leq \left|\widetilde{H}^\infty(W)-
\widetilde{H}^\infty(\widetilde{W})\right|+
\left|\widetilde{H}^\infty(\widetilde{W})- T_{2\pi ia}\circ
\ov{H}^0\circ(\ov{H}^\infty)^{-1}\circ T_{D'}(W)\right|\\
\qquad\qquad\qquad\quad+ \left|T_{2\pi ia}\circ
\ov{H}^0\circ(\ov{H}^\infty)^{-1}\circ T_{D'}(W) -T_{2\pi i
a}\circ(\ov{H}^\infty)^{-1}\circ T_{-2\pi i
a}(W)\right|.\end{array}$$ The second term vanishes from
\eqref{necessary_condition2}, and the first and third terms are
small from Lemma~\ref{lemma_estimates}(ii).

We have obtained \eqref{psi_infty_summable} by studying the
equations \eqref{equation_H_0}, \eqref{equation_H_infty},
\eqref{equation_H0_ov} and \eqref{equation_Hinfty_ov}, which come
from comparing the two presentations on the right of
Figure~\ref{strips} on a region located on top of the fundamental
hole. To obtain \eqref{psi0_summable} we instead compare on a region
located at the bottom and we replace the four equations
\eqref{equation_H_0}, \eqref{equation_H_infty},
\eqref{equation_H0_ov} and \eqref{equation_Hinfty_ov}  by the four
equations
$$
\begin{cases}
\widetilde{K}^0\circ \widetilde{\Psi}^0\circ T_{\tilde{\alpha}^0} =
T_{\tilde{\alpha}^0}\circ \widetilde{K}^0,\\
\widetilde{K}^\infty\circ \widetilde{\Psi}^\infty\circ
T_{\tilde{\alpha}^0} = T_{\tilde{\alpha}^\infty}\circ
\widetilde{K}^\infty,
\end{cases}\qquad\begin{cases}
\ov{K}^0\circ T_{\ov{\alpha}^0}\circ \ov{\Psi}^0 =
T_{\ov{\alpha}^0}\circ \ov{K}^0,\\
\ov{K}^\infty\circ T_{\ov{\alpha}^0}\circ \ov{\Psi}^\infty =
T_{\ov{\alpha}^\infty}\circ \ov{K}^\infty,
\end{cases}
$$
which have the solutions \begin{equation}\begin{cases}
\widetilde{K}^0=T_{-\tilde{\alpha}^0}\circ \widetilde{H}^0\circ T_{\tilde{\alpha}^0},\\
\widetilde{K}^\infty=T_{-\tilde{\alpha}^0}\circ\widetilde{H}^\infty\circ
T_{\tilde{\alpha}^0},
\end{cases}\qquad\begin{cases}
\ov{K}^0=T_{\ov{\alpha}^0}\circ\ov{H}^0\circ T_{-\ov{\alpha}^0},\\
\ov{K}^\infty=T_{\ov{\alpha}^0}\circ\ov{H}^\infty\circ
T_{-\ov{\alpha}^0}.\end{cases}\label{sol_K}\end{equation} We verify
that:
\begin{equation}\begin{cases}
(\widetilde{K}^0)^{-1}=\widetilde{\Psi}^0+O(\exp(2\pi i \tilde{\alpha}^0)),\\
\widetilde{K}^\infty=id+O(\exp(2\pi i \tilde{\alpha}^0)),
\end{cases}\qquad\begin{cases}
\ov{K}^0=\ov{\Psi}^0+O(\exp(-2\pi i \ov{\alpha}^0)),\\
\ov{K}^\infty=id+O(\exp(-2\pi i
\ov{\alpha}^0)).\end{cases}\label{formula_K}\end{equation}

Replacing \eqref{sol_K} in the compatibility condition
\eqref{necessary_condition2} yields
$$\widetilde{K}^\infty\circ (\widetilde{K}^0)^{-1}=
T_{-\tilde{\alpha}^0-\ov{\alpha}^0+2\pi i a}\circ
\ov{K}^0\circ(\ov{K}^\infty)^{-1}\circ
T_{\widetilde{\alpha}^0+\ov{\alpha}^0+D_\eps'}.$$ Finally using that
$\tilde{\alpha}^0+\ov{\alpha}^0=2\pi i a$ we have
$$\widetilde{K}^\infty\circ (\widetilde{K}^0)^{-1}=
 \ov{K}^0\circ(\ov{K}^\infty)^{-1}\circ T_{2\pi
 ia+D_\eps'}.$$
 Since $D_\eps'+2\pi i a=O(\exp(-2\pi i \ov{\alpha}^0))$, we get \eqref{psi0_summable}.\hfill $\Box$
\smallskip

Let us recall the following theorem which is a well-known
generalization of a corollary of the Ramis-Sibuya Theorem \cite{RS}.
This theorem will be used to show the $1/2$-summability of
$\Psi_{\hat{\eps}}^0$ and $\Psi_{\hat{\eps}}^\infty$ in $\hat{\eps}$.

\begin{theorem}\label{thm:RS}  Let $\{S_1,\ldots,S_k\}$ be a covering of a
punctured disk $\mathbb D_{\eps}=\{\eps;0<|\eps|<r\}$ by $k$ sectors
arranged so that only consecutive sectors overlap (taking
$S_{k+1}=S_1$). Let $\Psi_i(\eps,\nu)$ be holomorphic and bounded
functions defined on $S_i\times U$, where $U$ is a neighborhood of
the origin in $\nu$-space and $\nu$ is a multi-parameter. Moreover
let the functions $\Psi_j$ satisfy
$$|\Psi_i(\eps,\nu)-\Psi_{i+1}(\eps,\nu)|\leq
a\exp\left(-\frac{b}{|\eps|^s}\right)$$ on $(S_i\cap S_{i+1})\times U$, with
$a$ and $b$ positive numbers. Then there
exists a power series $$\hat{\Psi}(\eps,\nu)=
\sum_{n=0}^\infty \beta_n(\nu)\eps^n,$$
where the $\beta_n(\nu)$ are
analytic on $U$, and positive numbers $A$ and $C$ such that
\begin{enumerate}
\item for all $n\geq0$
$$|\beta_n(\nu)|\leq C A^n\; (n!)^{1/s};$$
\item for each subsector $S$ of $S_j$, $j=1,\ldots,k$, there exist
constants $A_S,C_S>0$ such that for all $\nu \in S$
$$\left|\Psi_j(\eps,\nu)-\sum_{m=0}^{N-1}\beta_n(\nu)\eps^n\right|<C_SA_S^N|\eps|^N\;(N!)^{1/s}.$$
\end{enumerate}
Moreover, if one of the $\Psi_i(\eps,\nu)$ can be extended to a
sector $S$ of opening greater than $\pi/s$, then $\hat{\Psi}$ is
$s$-summable in $\eps$ in the sector $S$.
\end{theorem}

\begin{corollary} The components $\Psi_{\hat{\eps}}^0$ and $\Psi_{\hat{\eps}}^\infty$ of the modulus
of a germ of family of diffeomorphisms normalized so that the
compatibility condition is satisfied in the form
\eqref{necessary_condition2} are $1/2$-summable in $\eps$. The
direction of non-summability is the Glutsyuk direction $\mathbb
R^+$.
\end{corollary}
\noindent{\scshape Proof.} This follows directly from Theorem~\ref{thm:RS} above
using the estimates \eqref{psi0_summable} and \eqref{psi_infty_summable}
of Theorem~\ref{thm_flatness}.\hfill $\Box$

\medskip

We can now refine Theorem~\ref{thm:local}.

\begin{theorem}\label{thm:local_revisited}
Under the hypotheses of Theorem~\ref{thm:local}, we also have
the further conclusion
\begin{itemize}\item For $\eps\in V_G$ \begin{equation}
|\ov{f}_{\ov{\eps}}(z)-\tilde{f}_{\tilde{\eps}}(z)|
<B\exp\left(-\frac{A}{|\sqrt{\hat{\eps}}|}\right).\label{flatness_f}
\end{equation}
The estimate is uniform in the $\nu_i$.
Thus $f$ is $\frac12$-summable in $\hat{\eps}$.
\end{itemize}

\end{theorem}
\noindent{\scshape Proof.} The only thing to prove is the estimate
\eqref{flatness_f}. We use the shape of the strips as in
Figure~\ref{strips_adjusted_compare} so that the functions
$\xi^{0,\infty}$ defined in \eqref{def_chi} satisfy
\begin{equation}\begin{cases} \left|\ov{\xi}^0-\tilde{\xi}^0\right|<B_1\exp\left(-\frac{A_1}{|\sqrt{\hat{\eps}}|}\right),\\
\left|\ov{\xi}^\infty-\tilde{\xi}^\infty\right|<B_2\exp\left(-\frac{A_2}{|\sqrt{\hat{\eps}}|}\right).\end{cases}
\label{estime_chi}\end{equation} The functions $\ov{\xi}^{0,\infty}$
and $\tilde{\xi}^{0,\infty}$ come from conjugating
$\ov{\Psi}^{0,\infty}$ and $\widetilde{\Psi}^{0,\infty}$. The
vertical part of the strips  are common for $\hat{\eps}\in\ov{V}$
and $\hat{\eps}\in \widetilde{V}$.
 Then it is clear that
\eqref{estime_chi} follows from the analyticity of $q_{\hat{\eps}}$
in the region corresponding to the vertical parts of the strips, so
we can use \eqref{psi0_summable} and \eqref{psi_infty_summable}. The
other parts are included in regions corresponding to $|Im
W|>\frac{Y_2}{|\sqrt{\hat{\eps}}|}$ for some $Y_2>0$ independent of
$\eps$ by Lemma~\ref{lemma_strips_adjusted} where
Lemma~\ref{lemma_estimates} allows to conclude that
$\left|\ov{\Psi}^{0,\infty}\right|,\left|\tilde{\Psi}^{0,\infty}\right|<B_2\exp\left(-\frac{A_2}{|\sqrt{\hat{\eps}}|}\right)$
from which
$\left|\ov{\xi}^{0,\infty}\right|,\left|\tilde{\xi}^{0,\infty}\right|<B_3\exp\left(-\frac{A_3}{|\sqrt{\hat{\eps}}|}\right)$
follows for some positive constants $A_j,B_j$. Indeed, we proved
before that the solution of the Beltrami equation depends
analytically on $\hat{\eps}$. Moroever the solutions of two Beltrami
equations where the Beltrami fields satisfy
$\left|\ov{\mu}-\tilde{\mu}\right|<
B_4\exp\left(-\frac{A_4}{|\sqrt{\hat{\eps}}|}\right)$ and same
values at 3 chosen points also satisfy such type of estimate, from
which the result follows. \hfill $\Box$

\begin{lemma}\label{lemma_J} Under the hypotheses of Theorem~\ref{thm:local_revisited} and
the compatibility condition \eqref{necessary_condition}, there
exists a  neighborhood  $U'$ of the origin such that for each
$\hat{\eps}\in\ov{V}$ there exists a conjugacy $J_{\hat{\eps}}$
between $\ov{f}=f_{\hat{\eps}}$ and $\tilde{f}=f_{\hat{\eps}e^{2\pi
i}}$ over $U'$. The conjugacy depends analytically on $\hat{\eps}$
and tends to the identity as $\hat{\eps}\to 0$. Moreover there
exists constants $A',B'>0$ such that $J_{\hat{\eps}}$ satisfies
$$\left|J_{\hat{\eps}}-id\right|<B'\exp\left
(-\frac{A'}{|\sqrt{\hat{\eps}}|}\right).\label{estime_J} $$
\end{lemma}
\noindent{\scshape Proof.} Let us recall that for $\arg
\hat{\eps}\in (-\delta,\delta)$ we consider $\ov{f}_{\ov{\eps}}$. We
compare with the point of view for $\arg \hat{\eps}\in (2\pi
-\delta,2\pi+\delta)$ in which we consider
$\tilde{f}_{\tilde{\eps}}$. In both cases the singular point
$-\sqrt{\hat{\eps}}$ (resp. $\sqrt{\hat{\eps}}$) is attached to the
upper index $0$ (resp. $\infty$).

\begin{figure}\begin{center}
\subfigure[Domain of $\ov{\gamma}_{\ov{\eps}}^0$ and $
\tilde{\gamma}_{\tilde{\eps}}^{\infty}$]{\includegraphics[width=4.5cm]{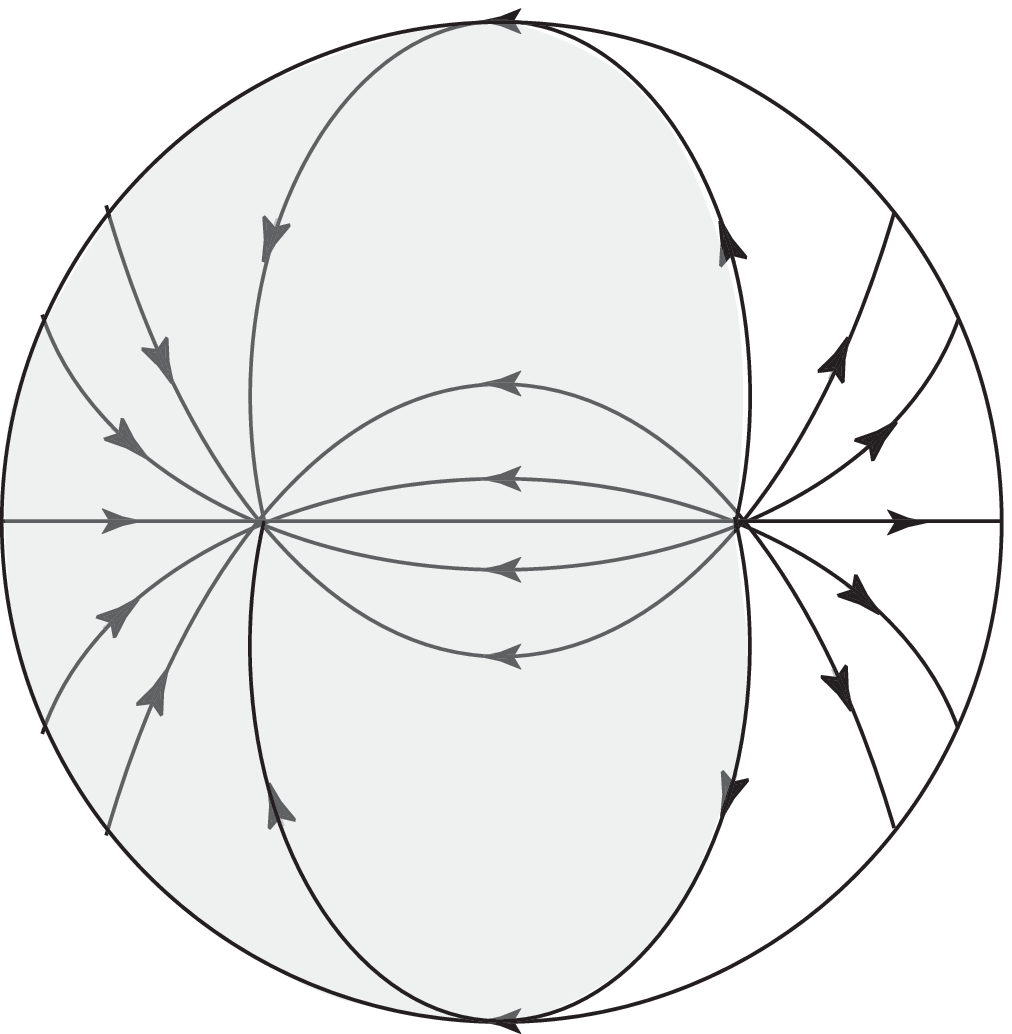}}
\qquad\qquad\qquad\subfigure[Domain of
$\ov{\gamma}_{\ov{\eps}}^\infty$ and $
\tilde{\gamma}_{\tilde{\eps}}^0$]{\includegraphics[width=4.5cm]{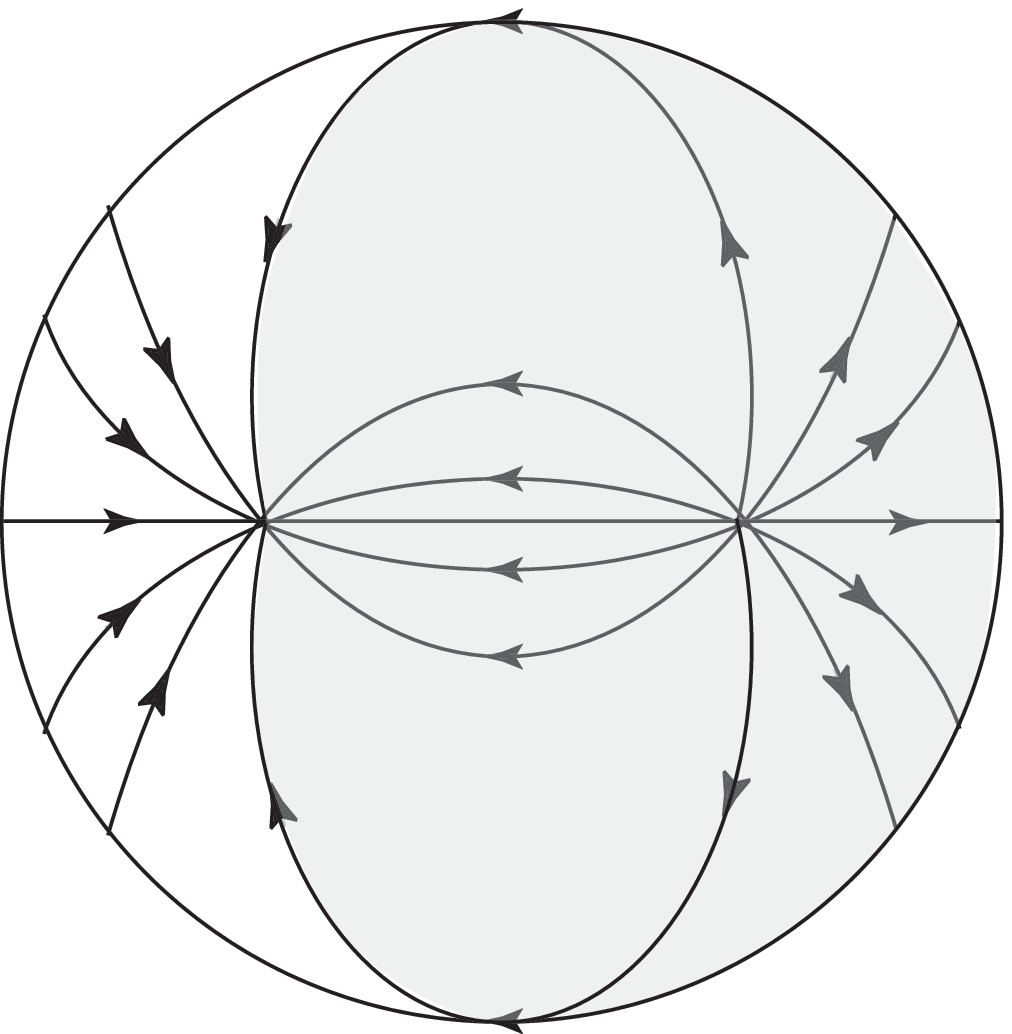}}
\caption{\label{domains} The domains of definition of the
normalizing maps}\end{center}
\end{figure}

We consider the \lq\lq normalizing  maps" in the neighborhoods of
the two singular points given by $\ov{\gamma}_{\ov{\eps}}^0$, $
\ov{\gamma}_{\ov{\eps}}^{\infty}$,  (resp.
$\tilde{\gamma}_{\tilde{\eps}}^{0}$, $
\tilde{\gamma}_{\tilde{\eps}}^{\infty}$), which are tangent to the
identity. These are the maps which transform $\ov{f}$ (resp.
$\tilde{f}$) to the model, i.e. the time one map of the vector field
\eqref{mod.3}. The advantage of these maps over the linearizing maps
is that their limits exist when $\hat{\eps}\to 0$ and that they do
not explode at the other singular point. It is known \cite{G} that
the union of the domains of $\ov{\gamma}_{\ov{\eps}}^0$ and
$\ov{\gamma}_{\ov{\eps}}^\infty$ (resp.
$\tilde{\gamma}_{\tilde{\eps}}^0 $ and
$\tilde{\gamma}_{\tilde{\eps}}^\infty $) is a whole covering of
$U_{\ov{\eps}}$ (resp. $U_{\tilde{\eps}}$) and that they overlap:
indeed the domains have a form as in Figure~\ref{domains}.

\begin{figure}\begin{center}
\subfigure[Subdomain of $\ov{\gamma}_{\ov{\eps}}^0$ and $
\tilde{\gamma}_{\tilde{\eps}}^{\infty}$]{\includegraphics[width=4.5cm]{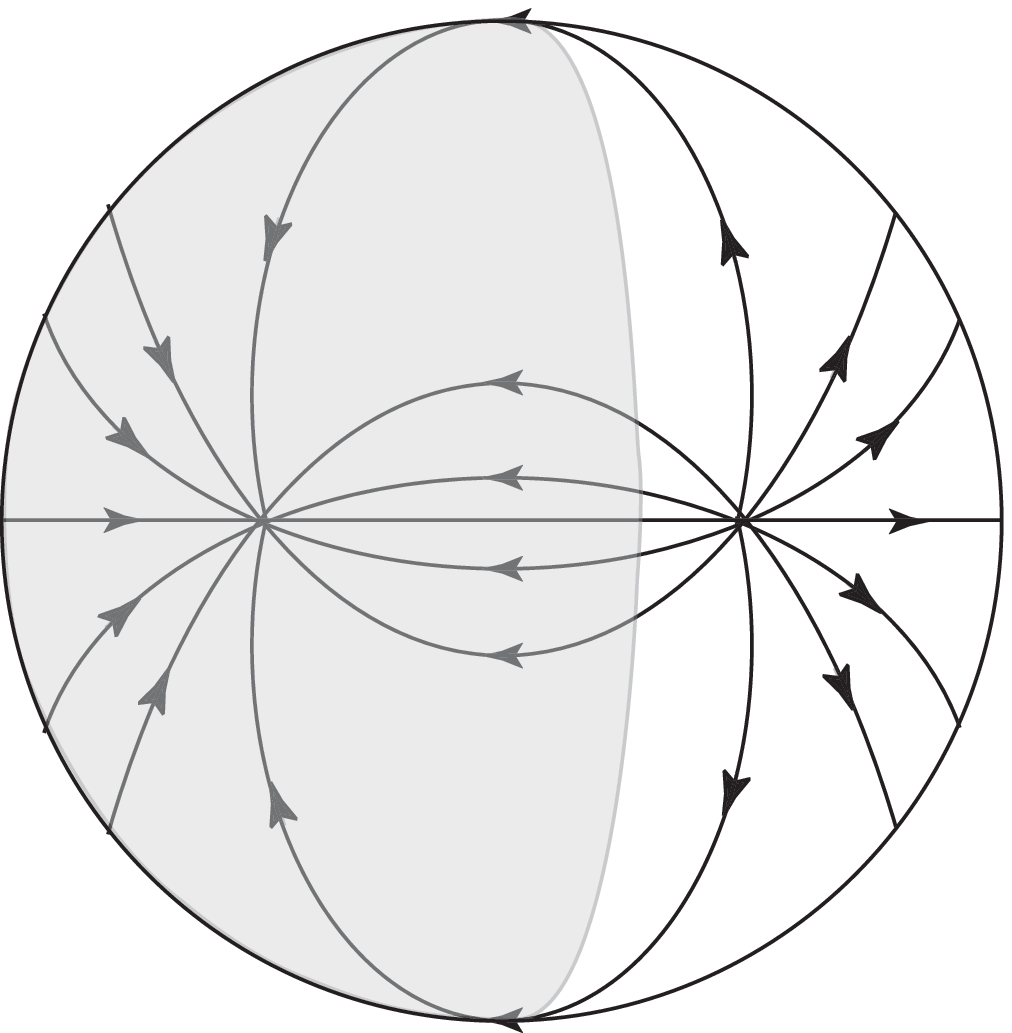}}
\qquad\qquad\qquad\subfigure[Subdomain of
$\ov{\gamma}_{\ov{\eps}}^\infty$ and $
\tilde{\gamma}_{\tilde{\eps}}^0$]{\includegraphics[width=4.5cm]{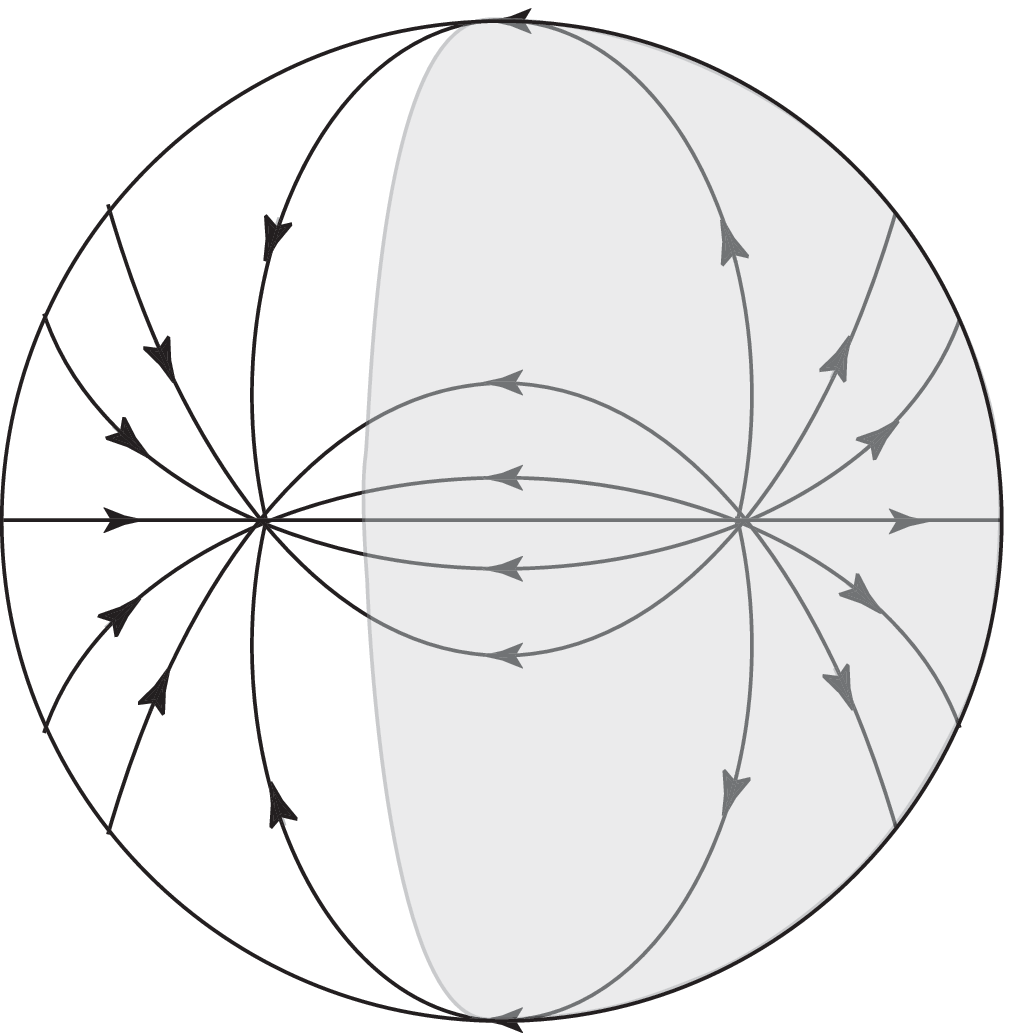}}
\caption{\label{smaller_domains} Smaller domains of definition of
the normalizing maps whose union covers $U$}\end{center}
\end{figure}

\begin{figure}\begin{center}
\subfigure[Domain of Fatou Glustyuk coordinate on the side of the
attracting
point]{\includegraphics[width=4cm]{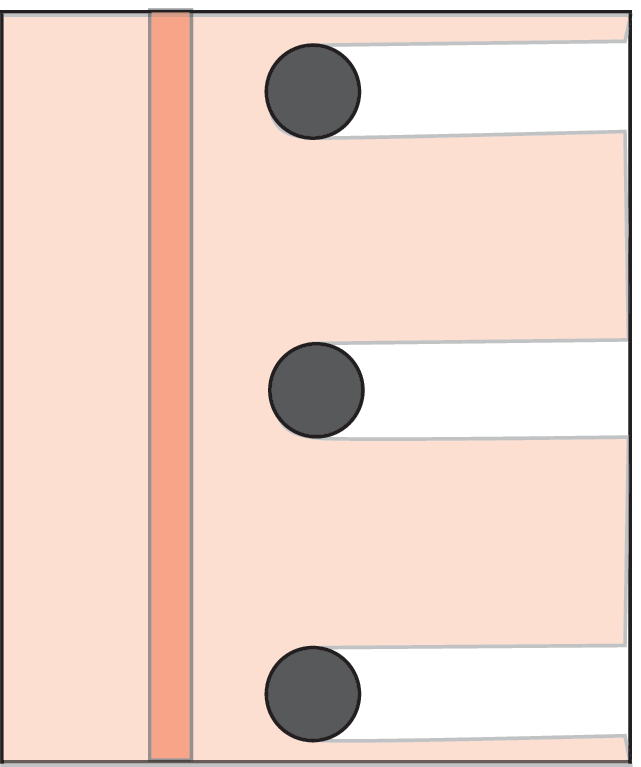}}
\qquad\qquad\qquad\subfigure[Domain of Fatou Glustyuk coordinate on
the side of the repelling
point]{\includegraphics[width=4cm]{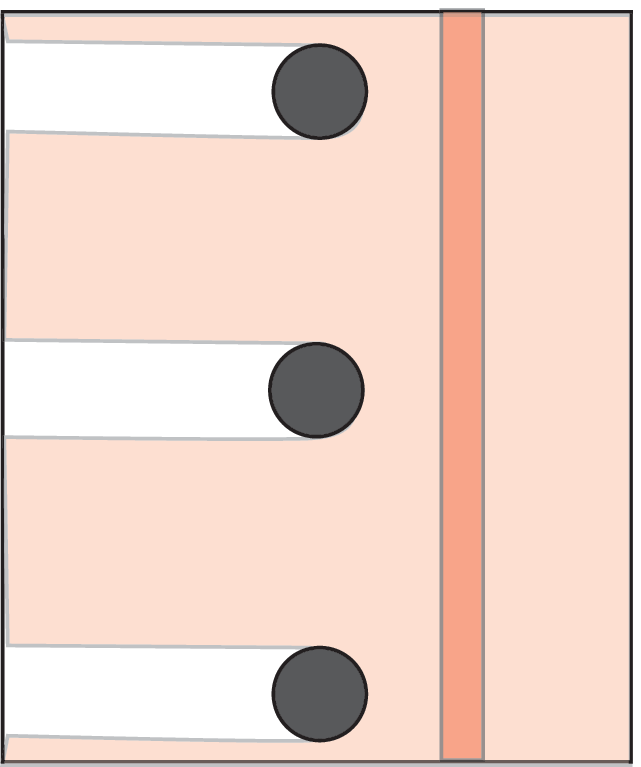}}
\caption{\label{domain_Fatou_Glutsyuk} The domains of definition
(translation domains) of the Fatou Glutsyuk coordinates. The
darkened strip is where the construction is first performed.
}\end{center}
\end{figure}

We will restrict to smaller domains as in
Figure~\ref{smaller_domains} whose union covers $U$.
On these smaller
domains we will show that there exist positive constants $A_0,B_0$
such that
\begin{equation}
\begin{cases} \left|\ov{\gamma}_{\ov{\eps}}^0(z)-
\tilde{\gamma}_{\tilde{\eps}}^\infty(z)\right|<B_0\exp\left
(-\frac{A_0}{|\sqrt{\hat{\eps}}|}\right),\\
\left|\ov{\gamma}_{\ov{\eps}}^\infty(z)-\tilde{\gamma}_{\tilde{\eps}}^0\right|<
B_0\exp\left
(-\frac{A_0}{|\sqrt{\hat{\eps}}|}\right).\end{cases}\label{estime_gamma}\end{equation}
From the maximum principle it suffices to prove that these estimates
hold on an annulus extending to the boundary of these subdomains. To
get the result we need to pass to the Fatou Glutsyuk coordinates.
Indeed these normalizing maps come from conjugating the Fatou
Glustsyuk
 coordinates with the map $q_{\hat{\eps}}^{-1}$. The
Fatou Glutsyuk coordinates are constructed as follows. We lift the
map $f_{\hat{\eps}}$ to
$$F_{\hat{\eps}}=q_{\hat{\eps}}^{-1}\circ f_{\hat{\eps}}\circ
q_{\hat{\eps}}.$$ The  Fatou Glutsyuk coordinates
$\Phi_{\hat{\eps}}^{0,\infty}$ satisfy
\begin{equation}\Phi_{\hat{\eps}}^{0,\infty}\circ
F_{\hat{\eps}}=T_1\circ
\Phi_{\hat{\eps}}^{0,\infty},\label{Fatou_eq}\end{equation} i.e.
they conjugate $F_{\hat{\eps}}$ with $T_1$ which is the time-one map
of the vector field $\frac{\partial}{\partial W}$. They are first
constructed on a strip of horizontal width $N$ and parallel to the
line of holes, and then extended to the maximal domain of definition
(called translation domain in \cite{MRR}) by means of
\eqref{Fatou_eq} (see Figure~\ref{domain_Fatou_Glutsyuk}). Both
$F_{\hat{\eps}}$ and $\Phi_{\hat{\eps}}^{0}$ (resp. $F_{\hat{\eps}}$
and $\Phi_{\hat{\eps}}^{\infty}$) commute with $T_{\alpha^0}$ (resp.
$T_{\alpha^\infty}$) on the side of $-\sqrt{\hat{\eps}}$ (resp.
$\sqrt{\hat{\eps}}$).

From the relation \eqref{flatness_f} it follows that there holds a similar relation
between $\ov{F}_{\ov{\eps}}=q_{\ov{\eps}}^{-1}\circ
\ov{f}_{\ov{\eps}}\circ q_{\ov{\eps}}$ and
$\widetilde{F}_{\tilde{\eps}}=q_{\tilde{\eps}}^{-1}\circ
\tilde{f}_{\tilde{\eps}}\circ q_{\tilde{\eps}}$:
\begin{equation}|\ov{F}_{\ov{\eps}}(W)-
\widetilde{F}_{\tilde{\eps}}(W)|<B_1\exp\left(-\frac{A_1}{|\sqrt{\hat{\eps}}|}\right)\label{bound_F}\end{equation}
for some positive constants $A_1,B_1$. An easy way to check
\eqref{bound_F} is the following: the map $f_{\hat{\eps}}$ is the
sum of a $1/2$-summable series in $\hat{\eps}$ with analytic
coefficients in $z$. Hence so is the case of its composition with
analytic maps. Moreover it is shown in \cite{MRR} that
$|F_{\hat{\eps}}(W)-W-1|<\frac14$ for $r,\rho$ sufficiently small.
Hence $|\ov{F}_{\ov{\eps}}- \widetilde{F}_{\tilde{\eps}}|$ is
bounded, from which we can conclude to the existence of a bound
independent of $W$ of the special form appearing in
 \eqref{bound_F}.

\begin{figure}\begin{center}
\subfigure[Strip of Fatou Glustyuk coordinate on the side of the
attracting
point]{\includegraphics[width=4.5cm]{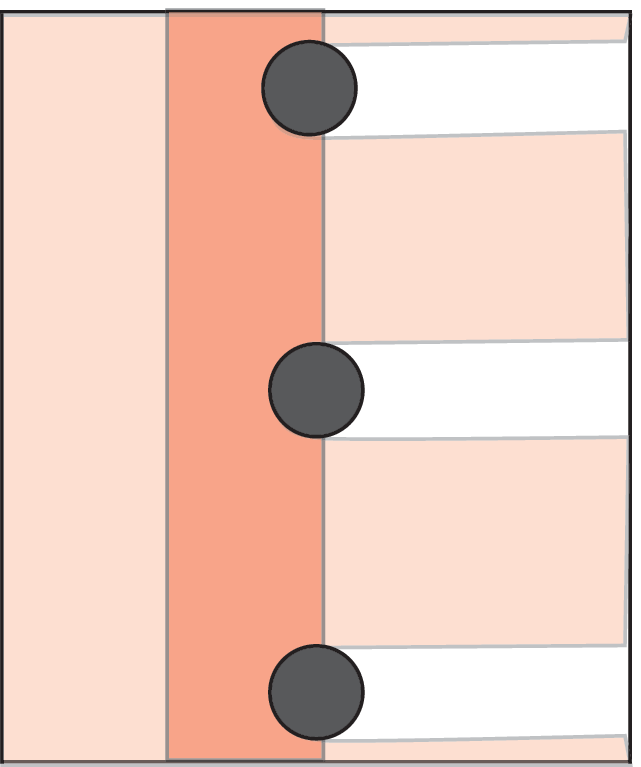}}
\qquad\qquad\subfigure[Strip of Fatou Glustyuk coordinate on the
side of the repelling
point]{\includegraphics[width=4.5cm]{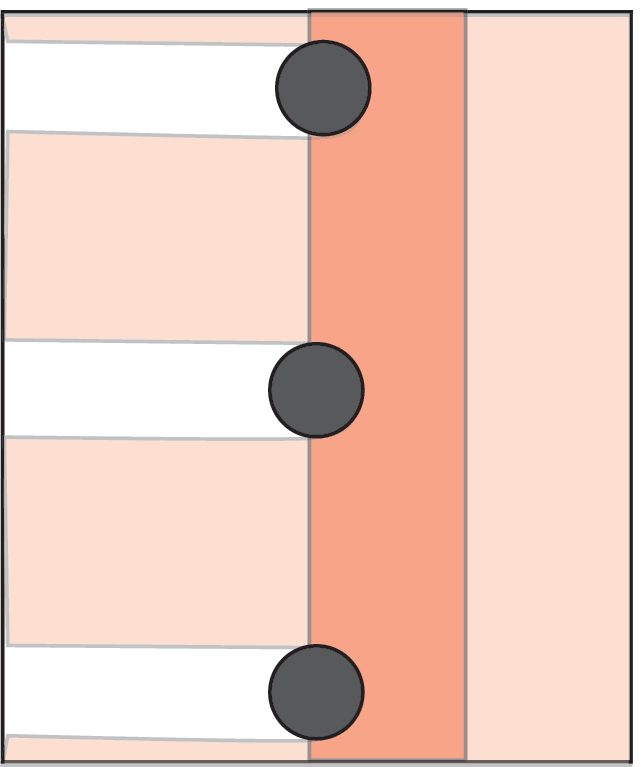}}
\caption{\label{strip_boundary} Strips whose projections by
$q_{\hat{\eps}}$ yield annular regions up to  the boundary of
subdomains as in Figure~\ref{smaller_domains} }\end{center}
\end{figure}

 The construction of
$\Phi_{\hat{\eps}}^{0,\infty}$ by means of Alhfors-Bers theorem
(\cite{MRR}) yields to estimates similar to \eqref{bound_F} for
Fatou Glutsyuk coordinates on the strips with positive constants
$A_2,B_2$:
\begin{equation}\begin{cases} |\ov{\Phi}_{\ov{\eps}}^0(W)-
\widetilde{\Phi}_{\tilde{\eps}}^\infty(W)|<B_2\exp\left(-\frac{A_2}{|\sqrt{\hat{\eps}}|}\right)
&\text{on the right strip,}\\
|\ov{\Phi}_{\ov{\eps}}^\infty(W)-
\widetilde{\Phi}_{\tilde{\eps}}^0(W)|<B_2\exp\left(-\frac{A_2}{|\sqrt{\hat{\eps}}|}\right)
&\text{on the left
strip,}\end{cases}\label{estimates_Fatou_Glutsyuk}\end{equation} as
long as we take the same normalization, for instance
$\ov{\Phi}_{\hat{\eps}}^0(Z_0)=Z_0=\widetilde{\Phi}_{\tilde{\eps}}^\infty(Z_0)$
 (resp.
$\ov{\Phi}_{\hat{\eps}}^\infty(Z_1)=Z_1=\widetilde{\Phi}_{\tilde{\eps}}^0(Z_1)$)
on the right (resp. left) strip.
The relation \eqref{Fatou_eq}
implies that for all $n\in \mathbb Z$
\begin{equation}\Phi_{\hat{\eps}}^{0,\infty}\circ
F_{\hat{\eps}}^n=T_n\circ
\Phi_{\hat{\eps}}^{0,\infty}.\label{Fatou_eq_n}\end{equation} This
in turn ensures that for any $N$ there exist constant $a_N,b_N$ such
that estimates of the form \eqref{estimates_Fatou_Glutsyuk} with
$A_2$ (resp. $B_2$) replaced by $a_N$ (resp. $b_N$) are valid in a
strip parallel to the holes of horizontal width $N$. We take $N$
sufficiently large so as to get the estimates on a strip of the form
as in Figure~\ref{strip_boundary}. Then the projection of these
strips by $q_{\hat{\eps}}$ yield annular regions up to  the boundary
of subdomains as in Figure~\ref{smaller_domains}.
Finally \eqref{estime_gamma} follows by conjugating
$\Phi_{\hat{\eps}}^{0,\infty}$ with $q_{\hat{\eps}}$.

 There exists a constant
$t(\hat{\eps})$ such
 that the map $J_{\eps}$ defined by:
\begin{equation}J_{\eps}=\begin{cases}
(\tilde{\gamma}_{\tilde{\eps}}^{\infty})^{-1}\circ v_{\hat{\eps}}^{t(\hat{\eps})}\circ \ov{\gamma}_{\ov{\eps}}^0\\
(\tilde{\gamma}_{\tilde{\eps}}^0)^{-1}\circ
\ov{\gamma}_{\ov{\eps}}^\infty\end{cases}\label{conjugacy}\end{equation}
is a conjugacy between $\ov{f}_{\ov{\eps}}$ and
$\tilde{f}_{\tilde{\eps}}$, where $v_{\hat{\eps}}^{t(\hat{\eps})}$
is the flow of the vector field \eqref{mod.3} for the time
$t(\hat{\eps})$. The compatibility condition ensures that this map
is well defined for an adequate choice of $t(\hat{\eps})$.

To determine the constant $t(\hat{\eps})$ we take a point $z_0\in
i\mathbb R^+$ on the imaginary axis close to the boundary of $U$.
Let us call
\begin{equation}\begin{cases}
J_1=(\tilde{\gamma}_{\tilde{\eps}}^{\infty})^{-1}\circ \ov{\gamma}_{\ov{\eps}}^0\\
J_2= (\tilde{\gamma}_{\tilde{\eps}}^0)^{-1}\circ
\ov{\gamma}_{\ov{\eps}}^\infty\\
J_t= (\tilde{\gamma}_{\tilde{\eps}}^{\infty})^{-1}\circ
v_{\hat{\eps}}^{t(\hat{\eps})}\circ
\ov{\gamma}_{\ov{\eps}}^0\end{cases}\label{conjugacy2}\end{equation}
The constant
 $t(\hat{\eps})$ is uniquely determined by the condition that
$J_t(z_0)=J_2(z_0)$. From their boundedness the  maps $J_1$ and
$J_2$ are uniformly continuous and equi-continuous because of the
existence of the limit when $\hat{\eps}\to 0$. Then
\eqref{estime_gamma} implies that
\begin{equation}|J_1(z)-J_2(z)|< B_5\exp\left
(-\frac{A_5}{|\sqrt{\hat{\eps}}|}\right)\label{estime_J2}\end{equation}
in the overlapping region. Moreover there exist positive constants
$A_6,B_6$ such that
$$
|t(\hat{\eps})|<B_6\exp\left
(-\frac{A_6}{|\sqrt{\hat{\eps}}|}\right).$$ The conclusion follows.
 \hfill $\Box$

\section{The global realization}\label{sec:global}

In Section~\ref{sect:local} we have shown how to realize a germ of
family
$\mathbf{\Psi}=(\Psi^0_{\hat{\eps}},\Psi^\infty_{\hat{\eps}})_{\hat{\eps}\in
V_{\rho,\delta}}$ as the modulus of a germ of family of
diffeomorphisms $f_{\hat{\eps}}$ and in Theorem~\ref{thm_flatness}
of Section~\ref{sect:compatibility} we have identified a necessary
compatibility condition so that the family $\mathbf{\Psi}$ be
realizable in a uniform family $g_{\eps}$.

We want to show that this condition is also sufficient. The idea is
the same as for the local realization: we realize the family as a
$2$-dimensional family of diffeomorphisms on an abstract
$2$-dimensional manifold and  we show that this manifold  is
holomorphically equivalent to a neighborhood of the origin minus
$\{\eps=0\}$ via the Newlander-Nirenberg theorem.

When dealing with the global realization we must work with open
sets. So we will consider open sectors in $\hat{\eps}$-space. We
consider the sector $V_{\rho',\delta'}$ constructed in the proof of
Theorem~\ref{thm:local}. Let $\delta\in (0,\delta')$ such that the
Glutsyuk modulus is defined for $\arg{\hat{\eps}}\in (-\delta,
\delta)$ and $\arg{\hat{\eps}}\in
(2\pi-\delta, 2\pi+\delta)$. We call
$$V_{\rho'}=\{\hat{\eps}\in
V_{\rho',\delta'}\setminus\{0\}|\arg{\hat{\eps}}\in(-\delta,
2\pi+\delta)\}.$$ We have the two subsectors $\ov{V}$ and
$\widetilde{V}$ defined in \eqref{secteur_V_eta}.

\begin{theorem}\label{thm_global_real} We consider a germ of function $a(\eps)$ analytic
in $\eps$ and a germ of family
$(\Psi^0_{\hat{\eps}},\Psi^\infty_{\hat{\eps}})$ for $\hat{\eps}$ in
some $V_{\rho,\delta}$ satisfying the hypotheses of
Theorem~\ref{thm:local}  and the compatibility condition
\eqref{necessary_condition2}. We suppose that $\delta$ is chosen
sufficiently small so that the conclusion of
Theorem~\ref{thm:local_revisited} holds. Then there exists a germ of
an analytic family of diffeomorphisms
\begin{equation} g_{\eps}= z+(z^2-\eps)(1+
O(\eps)+O(z))\label{germ_g_eps}\end{equation}  whose modulus is given by
$(a(\eps),[\Psi^0_{\hat{\eps}},\Psi^\infty_{\hat{\eps}}])$ in some
$V_{\rho',\delta'}$.  Moreover, if the functions $a(\eps,\nu)$ and
$\Psi_{\hat{\eps},\nu}^{0,\infty}$ depend analytically on
$(k-1)$-parameters $\nu$, then the function $g_{\eps,\nu}$ depends
analytically on $\nu$.
\end{theorem}
\noindent{\scshape Proof.} We consider the sector
$V_{\rho',\delta'}$  (with $\delta'=\delta$) constructed in the
proof of Theorem~\ref{thm:local}. We can of course suppose that
$\delta\in (0,\frac{\pi}4)$ and that $\delta$ is sufficiently small
so that the Glutsyuk modulus is defined for $\arg{\hat{\eps}}\in
(-2\delta, 2\delta)$ and
$\arg{\hat{\eps}}\in (2\pi - 2\delta, 2\pi+2\delta)$. (To realize
this requirement it suffices to take $\delta= \frac{\delta'}{2}$
where $\delta'$ is constructed in Theorem~\ref{thm:local}.)

For each $\hat{\eps}\in V_{\rho'}$ we have realized the modulus over
an open set $U_{\hat{\eps}}$ of $\mathbb C$ constructed as in the
proofs of Theorem~\ref{thm:local} and
Theorem~\ref{thm:local_revisited}. For all $\hat{\eps}\in
V_{\rho'}$, $U_{\hat{\eps}}$ contains a fixed  disk $B(0,r)$ and the
two fixed points lie inside $B(0,r)$. We can suppose $r$
sufficiently small so that $B(0,r)\subset U'$ where $U'$ is the open
neighborhood of $\pm \sqrt{\hat{\eps}}$ in Lemma~\ref{lemma_J}. So
for the rest of the proof we will suppose $U_{\hat{\eps}}=B(0,r)$.

 We consider the open
set of $\mathbb C\times \hat{\mathbb C}$ defined by
 $$\mathbb{U}=\cup_{\hat{\eps}\in
 V_{\rho'}}(U_{\hat{\eps}}\setminus\{\pm\sqrt{\hat{\eps}}\},\hat{\eps}).$$
This space is endowed with a projection $\Pi:\mathbb{U}\rightarrow
V_{\rho}$.

\begin{figure}\begin{center}
\includegraphics[width=5cm]{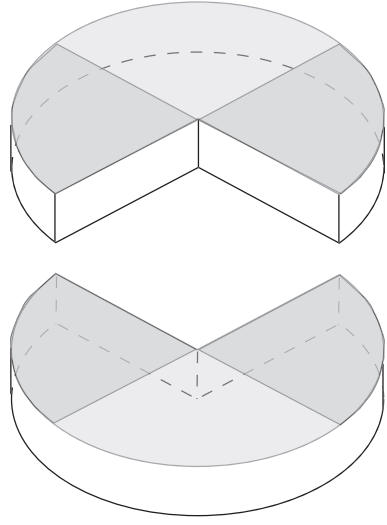}
\caption{\label{two_sectors}The two sectors $\mathbb U_1$ and
$\mathbb U_2$}
\end{center}
 \end{figure}

We cover $V_{\rho'}$ with the two sectors $V_{\rho'}^1$ and
$V_{\rho'}^2$ defined by
$$ \begin{cases}
V_{\rho'}^1=\{\hat{\eps}\in
V_{\rho'}|\arg{\hat{\eps}}\in(-\delta, \pi+\delta)\}\\
V_{\rho'}^2=\{\hat{\eps}\in V_{\rho'}|\arg{\hat{\eps}}\in(\pi
-\delta, 2\pi+\delta)\}\end{cases}$$ Their inverse images in
$\mathbb U$ are called $\mathbb U_1= \Pi^{-1}(V_{\rho'}^1)$ and
$\mathbb U_2= \Pi^{-1}(V_{\rho'}^2)$ (Figure~\ref{two_sectors}).

We construct a complex manifold $\mathcal{M}$ with atlas given by
$\{\mathbb U_1,\mathbb U_2\}$. The transition function on $\mathbb
U_1\cap \mathbb U_2$ (i.e. when $\arg \hat{\eps}\in
(\pi-\delta,\pi+\delta)$) is the identity. The other transition
function is obtained as follows: we make the gluing of
$\Pi^{-1}(\ov{V})$ with $\Pi^{-1}(\widetilde{V})$ in the following
way: we identify $(z,\ov{\eps})\in \Pi^{-1}(\ov{\eps})$ with
$(J_{\eps}(z),\tilde{\eps})\in \Pi^{-1}(\tilde{\eps})$ defined in
\eqref{conjugacy}. With this gluing $\ov{\eps}$ and $\tilde{\eps}$
simply become $\eps$. On $\mathcal{M}$ a global function $f_{\eps}$
is defined. It is given by $f_{\hat{\eps}}$ on each $\mathbb U_j$
and the definitions match because $J_{\eps}$ conjugates
$f_{\ov{\eps}}$ and $f_{\tilde{\eps}}$ .

On each of $\mathbb U_1$ and $\mathbb U_2$ we have respective
coordinates $(z_1,\eps)$ and $(z_2,\eps)$. We want to show that the
complex manifold $\mathcal{M}$ is holomorphically equivalent to a
neighborhood of the origin in $\mathbb C^2$ minus $\{\eps=0\}$.

Let $(\Theta_1,\Theta_2)$ be a partition of unity associated to the
covering $\{\mathbb U_1,\mathbb U_2\}$. As in
Theorem~\ref{thm:local}, we can suppose that the derivatives of
$\Theta_j$ grow no faster than a negative power of the variables. We
can also suppose that the $\Theta_j$ depend on $\eps$ alone. Let us
first construct a $C^\infty$-diffeomorphism
$$\Omega: \mathcal{M}\rightarrow(\mathbb
C^2,0)\setminus\{\eps=0\}$$ defined by $$\Omega=
\Theta_1\cdot(z_1,\eps)+\Theta_2\cdot(z_2,\eps)=(\Theta_1z_1+\Theta_2z_2,\eps).$$
This map is $C^\infty$. We will extend it by the identity on
$\eps=0$. To show that the extension is $C^\infty$ we use the fact that
the map $(z,\hat{\eps})\mapsto J_{\hat{\eps}}(z)$ has
$J_{\hat{\eps}}-id$ exponentially small in $\sqrt{\hat{\eps}}$
near $\hat{\eps}=0$ (see Lemma~\ref{lemma_J}). This endows
$\Omega(\mathcal{M})$ of two complex coordinates $(Z,\eps)$ where
\begin{equation}Z=\Theta_1z_1+\Theta_2z_2.\label{def_Z}\end{equation}

We now show that $\Omega$ induces an integrable almost complex
structure on $\Omega(\mathcal{M})$. Let us recall that an almost
complex structure is given by two forms  $\omega,\xi$ which are
$\mathbb C$-linear in the sense of this structure.

The almost complex structure is integrable when there exist
coordinates $(w_1,w_2)$ such that
$$\langle dw_1,dw_2\rangle_{\mathbb C} = \langle
\omega,\xi\rangle_{\mathbb C}.$$ In that case there exists a
$2\times 2$ invertible matrix $A$ whose entries are $C^\infty$
functions in $(Z,\eps)$ such that
$$\left(\begin{matrix}
\omega\\ \xi\end{matrix}\right) = A\left(\begin{matrix}dw_1\\
dw_2\end{matrix}\right)= A\, dw.$$ In particular,
$d\left(\begin{matrix} \omega\\ \xi\end{matrix}\right)= dA \wedge
dw$ contains no $(0,2)$ component. The Newlander-Nirenberg Theorem
asserts that this necessary condition is also sufficient for
integrability.

For the
second form of the complex structure we take $\xi=d\eps$.
 The other
form $\omega$ should play the role of $dZ$. It will be given by
\begin{equation}\omega=(\Omega^{-1})^*(\widetilde{\omega})\label{def_omega}\end{equation} for some form
$\widetilde{\omega}$ defined on $\mathcal{M}$. The form
$\widetilde{\omega}$ is given by $\widetilde{\omega}_j$ on the chart
$\mathbb U_j$. On $\mathbb U_2$ we take $\widetilde{\omega}_2=dz_2$.
On $\mathbb U_1 \cap \mathbb U_2$ we have $dz_1=dz_2$. So we want
$\widetilde{\omega}_1= dz_1$ on $\mathbb U_1 \cap \mathbb U_2$.
 On the region of the gluing  we have
$$\begin{array}{lll} dz_2&=&\frac{\partial J}{\partial \eps}d\eps + \frac{\partial
J}{\partial z_1}dz_1\\
&=& \tau_{\eps,1}d\eps+ (1+\tau_{\eps,2})dz_1,\end{array}$$ where
the two functions $\tau_{\eps,j}$ are exponentially flat in
$|\sqrt{\eps}|^{-1}$ near $\eps=0$. The gluing is done in the
following way: $\delta$ has been  chosen sufficiently small so that
$J_\eps$, and then $\tau_{\eps,j}$ exist for $\arg(\eps)\in
(-2\delta,2\delta)$. We take an increasing $C^\infty$ function
$\varphi: \mathbb R\rightarrow [0,1]$ such that
$$
\varphi(x)\equiv \begin{cases}0&x<-2\delta\\
1&x>-\delta.\end{cases}$$ Then
$$\widetilde{\omega}_1= dz_1+
\varphi(\arg{\eps})(\tau_{\eps,1}d\eps+ \tau_{\eps,2}dz_1).$$ From
its construction the form $\widetilde{\omega}=\widetilde{\omega}_j$
on $\mathbb U_j$ is well defined on $\mathcal{M}$, $C^\infty$ and of
type $(1,0)$.

Let us now remark that the difference $\omega-dZ$  decreases
exponentially fast as $\eps\to 0$. This comes from the fact that
$\tau_{\eps,j}$, $j=1,2$, are exponentially flat in
$|\sqrt{\eps}|^{-1}$ near $\eps=0$.

This allows to extend the almost complex structure
$\{d\eps,\omega\}$ to $\eps=0$, by taking the two forms $d\eps$ and
$dz$. The resulting almost complex structure is $C^\infty$ in a
neighborhood of the origin in $\mathbb C^2$.

To show that this complex structure satisfies the necessary
condition for integrability we need to show that
$\{d(d\eps),d\omega\}$ contains no terms of type $(0,2)$. Obviously
$d(d\eps)=0$, so we only need to study $d\omega$. From its
construction $d\widetilde{\omega}$ has no terms of type $(0,2)$. But
$\omega$ is obtained from the pull-back of $\widetilde{\omega}$.
Note that no terms containing $d\ov{\eps}$ may exist outside the
region $\arg\eps\in (-2\delta, \delta)$, since $\varphi$ is constant
there and either the $\Theta_j\equiv1$ or $z_1=z_2$. In the region
$\arg\eps\in (-2\delta, \delta)$ the maps $\tau_{\eps,j}$ are
holomorphic in $Z$ and the maps $\Theta_j$ depend on $\eps$ alone so
there is no possibility of a term in $d\ov{Z}$.

 Since the
almost complex structure satisfies the necessary condition for
integrability, we can apply the Newlander-Nirenberg Theorem
\cite{NN} to the manifold $\ov{\Omega(\mathcal{M})}$, where
$\ov{\Omega(\mathcal{M})}$ is the closure of $\Omega(\mathcal{M})$
obtained by adding $\eps=0$, $z\in U_0$. Indeed the complex
structure is integrable on $\Omega(\mathcal{M})$ and hence on
$\ov{\Omega(\mathcal{M})}$ by continuity. Then the local charts
which are holomorphic in the sense of this complex structure are
$C^\infty$. Hence there exists a diffeomorphism
$\Gamma:\ov{\Omega(\mathcal{M})} \cap \mathcal{U} \rightarrow\mathbb
C^2$, where $\mathcal{U}$ is a neighborhood of the origin in
$\mathbb C^2$,  which is holomorphic with respect to this structure
and whose image is a neighborhood of the origin in $\mathbb C^2$.
From the form of the complex structure it is clear that $\eps$ can
be taken as one of the complex coordinates. So we can suppose that
$\Gamma$ preserves $\eps$. The composition $\Gamma\circ\Omega$ is an
analytic diffeomorphism of an open set of $\mathcal{M}$ with a
neighborhood of the origin in $\mathbb C^2$. The map $\Gamma$ is not
unique. We can always choose it in such a way that it sends the
curve $z^2-\eps=0$ to the same curve.

We now conjugate the map $(f_\eps,\eps)$ with $\Gamma\circ\Omega$
yielding
$$(g_\eps,\eps)=(\Gamma\circ\Omega) \circ (f_\eps,\eps)\circ
(\Gamma\circ\Omega)^{-1}. $$ Since $g_\eps$ is bounded in the
neighborhood of $\eps=0$, it is possible to extend it to $\eps=0$ in
an analytic way. For each fixed $\eps$ the map $g_\eps$ is
conjugated to $f_\eps$ defined on the slice $\mathcal{M}_\eps$. By
continuity it is clear that $g_0$ is conjugated to $f_0=
\lim_{\hat{\eps}\to 0}f_{\hat{\eps}}$ where $f_{\hat{\eps}}$ was the
family of Theorem~\ref{thm:local}.  \hfill $\Box$

\section{Examples}\label{sect:Riccati}

In this section we consider the realization problem for a family
$(\Psi_{\hat{\eps}}^0,\Psi_{\hat{\eps}}^\infty)$ which is conjugate
under the map $w=E(W)=\exp(-2\pi i W)$ to a family of functions
\begin{equation}\begin{cases}
\psi_{\hat{\eps}}^0(w) = m_{A(\hat{\eps}), n}(w) =
\frac{w}{\left(1+A(\hat{\eps})w^n\right)^{1/n}}\\
\psi_{\hat{\eps}}^\infty(w)= L_{\exp(-4\pi^2 a(\eps))}\circ
T_{B(\hat{\eps}),n'}(w)=\exp(-4\pi^2
a(\eps))\left(w^{n'}+B(\hat{\eps})\right)^{1/n'}.
\end{cases}\end{equation}
When $n=1$, we drop the index n. For $n=n'=1$, such a modulus is
obtained for instance in the modulus of the holonomy of an unfolding
of a Riccati equation with a saddle-node (\cite{R2} or \cite{LR}),
so we will call it the \lq\lq Riccati case".

\subsection{The general case}

 Let
\begin{equation} \beta=\exp(-4\pi^2 a(\eps)),\end{equation} and
\begin{equation}\begin{cases} \ov{C}=\exp(-2\pi i \ov{\alpha}^0),\\
\widetilde{C}=\exp(-2\pi i
\tilde{\alpha}^0).\end{cases}\end{equation} Then we have
\begin{equation}
\ov{C}\beta=(\widetilde{C})^{-1}.\label{rel_C}\end{equation}

 As before, we compare the modulus
at values  $\ov{\eps}=\hat{\eps}$ and
$\tilde{\eps}=\hat{\eps}e^{2\pi i}$, which we denote by
$$\begin{cases}\ov{\psi}^0= m_{\ov{A},n},\\
\ov{\psi}^\infty=L_{\beta}\circ T_{\ov{B},n'},\end{cases}\qquad \begin{cases}\widetilde{\psi}^0= m_{\widetilde{A},n},\\
\widetilde{\psi}^\infty=L_{\beta}\circ
T_{\widetilde{B},n'}.\end{cases}$$

Let $$\begin{cases} \tilde{h}^0=E\circ \widetilde{H}^0\circ E^{-1},\\
\tilde{h}^\infty=E\circ \widetilde{H}^\infty\circ E^{-1},
\end{cases}
\qquad\begin{cases}
\ov{h}^0=E\circ \ov{H}^0\circ E^{-1},\\
\ov{h}^\infty=E\circ \ov{H}^\infty\circ E^{-1}.\end{cases}$$ They
satisfy respectively \begin{equation}\begin{cases} \tilde{h}^0\circ
L_{\widetilde{C}}\circ
\widetilde{\psi}^0=L_{\widetilde{C}}\circ\tilde{h}^0,\\
\tilde{h}^\infty\circ L_{\widetilde{C}}\circ
\widetilde{\psi}^\infty=L_{\widetilde{C}\beta}\circ\tilde{h}^\infty,\end{cases}
\qquad\begin{cases} \ov{h}^0\circ
\ov{\psi}^0\circ L_{\ov{C}}=L_{\ov{C}}\circ\ov{h}^0,\\
\ov{h}^\infty\circ \ov{\psi}^\infty\circ
L_{\ov{C}}=L_{\ov{C}\beta}\circ\ov{h}^\infty.\end{cases}\label{eq_h}\end{equation}

To calculate $\tilde{h}^0$, $\tilde{h}^\infty$, $\ov{h}^0$, and
$\ov{h}^\infty$ we use the following proposition

\begin{proposition}\label{calcul_Mobius} The  functions $m_{A,n}$ and $T_{B,n}$ satisfy:
\begin{description}
\item{(i)} $m_{A,n}\circ L_C = L_C\circ m_{AC^n,n}$;
\item{(ii)} $T_{B,n}\circ L_C= L_C\circ T_{B/C^n,n}$;
\item{(iii)} $m_{A,n}\circ m_{A',n}=m_{A+A',n}$;
\item{(iv)} $T_{B,n}\circ T_{B',n}=T_{B+B',n}$.
\end{description}\end{proposition}

\begin{theorem}
\begin{description} \item{(i)} The maps $\tilde{h}^0$, $\tilde{h}^\infty$, $\ov{h}^0$,
$\ov{h}^\infty$ are given by
$$\begin{cases}
\tilde{h}^0=m_{\tilde{d},n},&\quad\text{with}\qquad
\tilde{d}=\frac{\widetilde{A}}{1-\widetilde{C}^n},\\
\tilde{h}^\infty=T_{\tilde{e},n'},&\quad\text{with}\qquad
\tilde{e}=\frac{(\widetilde{C}\beta)^{n'}\widetilde{B}}{(\widetilde{C}\beta)^{n'}-1}=\frac{\widetilde{B}}{1-\ov{C}^{n'}},\\
\ov{h}^0=m_{\ov{d},n},&\quad\text{with}\qquad
\ov{d}=\frac{\ov{A}\,\ov{C}^n}{1-\ov{C}^n},\\
\ov{h}^\infty=T_{\ov{e},n'},&\quad\text{with}\qquad
\ov{e}=\frac{\beta^{n'}\ov{B}}{(\ov{C}\beta)^{n'}-1}=\frac{(\widetilde{C}\beta)^{n'}\ov{B}}{1-\widetilde{C}^{n'}}.\end{cases}$$
\item{(ii)} The compatibility condition can only be satisfied when either $A$ or $B$ vanish or we have $n=n'$.
The compatibility condition in the latter case is given by the
condition $\widetilde{A}\widetilde{B}=\ov{A}\,\ov{B}$, so that the
analytic invariant $AB$ depends analytically on $\eps$. The linear
changes of Glutsyuk coordinates $L_F$ and $L_G$ allowing to realize
the compatibility condition
\begin{equation}\tilde{h}^\infty\circ (\tilde{h}^0)^{-1}= L_{F}\circ
\ov{h}^0\circ (\ov{h}^\infty)^{-1}\circ
L_{G}\label{compatibility_mobius}\end{equation} are given by
$$\begin{cases}
F^n=-\frac{\widetilde{B}}{\ov{B}}\;\frac{(1-\ov{C}^n)((\beta\ov{C})^n
-1)-AB(\ov{C}\beta)^n}{\beta^n
(1-\ov{C}^n)^2},\\
G^n=-\frac{\widetilde{A}}{\ov{A}}\;\frac{\beta^{n}\left[(1-\ov{C}^n)((\beta\ov{C})^n
-1)-AB(\ov{C}\beta)^n\right]}{((\beta\ov{C})^n -1)^2}.\end{cases}
$$
Then
$$F^nG^n
= 1+2AB\widetilde{C}^{-n}(1+O(\ov{C}^n))=
1+2AB\widetilde{C}^{-n}+o(\widetilde{C}^{-n}).$$ This yields a
geometric interpretation of the analytic invariant $AB$ as a shift
between the two constants $F$ and $G$.
\item{(iii)} If $A(\hat{\eps})\equiv0$ (resp. $B(\hat{\eps}
)\equiv0$), then the
compatibility condition is given by $\ov{B}/\widetilde{B}$ (resp.
$\ov{A}/\widetilde{A}$) bounded and bounded away from 0. In
particular $\ov{B}$ and $\widetilde{B}$ (resp. $\ov{A}$ and
$\widetilde{A}$) vanish at the same values of $\eps$, with same
multiplicity.  In that case $FG=1$.
\end{description}\end{theorem}
\noindent{\scshape Proof.}
\begin{description} \item{(i)}
The result follows by applying Proposition~\ref{calcul_Mobius} in
\eqref{eq_h} and using \eqref{rel_C}.\item{(ii)} The compatibility
condition is that there exist nonzero constants $F$ and $G$ such
that $\tilde{h}^\infty\circ (\tilde{h}^0)^{-1}= L_{F}\circ
\ov{h}^0\circ (\ov{h}^\infty)^{-1}\circ L_{G}$, i.e
$$T_{\tilde{e},n'}\circ m_{-\tilde{d},n}=L_{FG}\circ m_{G^n\ov{d},n}\circ
T_{-\ov{e}/G^n,n'}.$$

Such an equation can obviously only be satisfied for $n=n'$, unless
$A\equiv0$ or $B\equiv 0$.

Let us calculate both sides when $n=n'$:
$$T_{\tilde{e},n}\circ
m_{-\tilde{d},n}(w)=\left(\frac{w^n(1-\tilde{d}\tilde{e})+\tilde{e}}{1-\tilde{d}w^n}\right)^{1/n}$$
and
$$L_{FG}\circ m_{G^n\ov{d},n}\circ
T_{-\ov{e}/G^n,n}=\left(\frac{\frac{F^nG^n}{1-\ov{d}\ov{e}}w^n-\frac{F^n\ov{e}}{1-\ov{d}
\ov{e}}}{1+\frac{G^n\ov{d}}{1-\ov{d}\ov{e}}w^n}\right)^{1/n}.$$ Then
the compatibility conditions become
\begin{equation}\begin{cases}
1-\tilde{d}\tilde{e}=\frac{F^nG^n}{1-\ov{d}\ov{e}},\\
\tilde{e}=-\frac{F^n\ov{e}}{1-\ov{d} \ov{e}},\\
\tilde{d}=-\frac{G^n\ov{d}}{1-\ov{d}\ov{e}}.\end{cases}\label{compatib_example}\end{equation}From
this we get
$$\begin{cases}
F^n=-\frac{\widetilde{B}}{\ov{B}}\;\frac{\ov{C}^n[(1-\ov{C}^n)(1-\widetilde{C}^n)-\ov{A}\,\ov{B}]}{(1-\ov{C}^n)^2},\\
G^n=-\frac{\widetilde{A}}{\ov{A}}\;\frac{(1-\ov{C}^n)(1-\widetilde{C}^n)-\ov{A}\,\ov{B}}{\ov{C}^n(1-\widetilde{C}^n)^2},\end{cases}
$$
and the compatibility condition linking $\ov{A}\,\ov{B}$ and
$\widetilde{A}\widetilde{B}$ becomes
$\tilde{d}\tilde{e}=\ov{d}\ov{e}$ which is equivalent to
$$\ov{A}\,\ov{B}=\widetilde{A}\widetilde{B}.$$ Since this product is an
invariant, we can simply note it by $AB$. Note that
$$F^nG^n=(1-\tilde{d}\tilde{e})^2=1+2AB\widetilde{C}^{-n}(1+O(\ov{C}^n))=
1+2AB\widetilde{C}^{-n}+o(\widetilde{C}^{-n}).$$

In the particular case $F=1/\beta$, i.e. the modulus family has been
normalized so as to satisfy \eqref{necessary_condition2}, then we
get that $G = \beta + O(\ov{C})$, which ensures
$\tilde{A}-\ov{A}=O(\ov{C})$ and similarly
$\tilde{B}-\ov{B}=O(\ov{C})$ as proved in
Theorem~\ref{thm_flatness}.

\item{(iii)} If $A\equiv0$, then
$\ov{d}=\tilde{d}=0$ in \eqref{compatib_example}, from which the
conclusion follows.\hfill $\Box$\end{description}

\begin{corollary}\label{ex:non_real}
No family $\left.\left(a(\eps), [m_{A(\hat{\eps}),n}, L_\beta\circ
T_{B(\hat{\eps}),n'}]\right)\right|_{\hat{\eps}\in V_{\rho,\delta}}$
is realizable as the modulus of a prepared family unfolding a
diffeomorphism with a parabolic fixed point when $n\neq n'$ and
neither $A(\eps)$ or $B(\eps)$ are identically zero.
\end{corollary}

\begin{remark} The Corollary~\ref{ex:non_real} shows the strength of the
compatibility condition. Indeed, while $(a(0), [m_{A(0),n},
L_\beta\circ T_{B(0),n'}])$ is realizable as the modulus of a single
diffeomorphism, its unfolding can never keep this simple form.
\end{remark}

\begin{theorem}\label{thm:Riccati} We consider a realizable family of triples
 $\left.\left(a(\eps), [m_{A(\hat{\eps}),n}, L_\beta\circ T_{B(\hat{\eps}),n}]\right)\right|_{\hat{\eps}\in
 V_{\rho,\delta}}$. It is possible to
choose analytic representatives of the modulus. The different
equivalence classes have a unique representative composed of a
triple of germs of analytic functions $(a(\eps), A(\eps),B(\eps))$,
with $a(\eps)$ arbitrary and $A(\eps)$, $B(\eps)$ of one of the
following type for some choice of $N_A$, $N_B \in \N =
\{0,1,\dots\}$.
\begin{description}
\item{(i)} $A(\eps)=\eps^{N_A}$, $B(\eps)=\eps^{N_B}B_1(\eps)$, with $B_1$ analytic satisfying
$B_1(0)\neq0$;
\item{(ii} $A(\eps)\equiv 0$, $B(\eps)=\eps^{N_B}$;
\item{(iii)} $A(\eps)=\eps^{N_A}$,  $B(\eps)\equiv 0$;
\item{(iv)} $A(\eps)=B(\eps)\equiv 0$.
\end{description}\end{theorem}
\noindent{\scshape Proof.} The compatibility condition shows that $AB$ is analytic in $\eps$.
Moreover we have shown
in Theorem~\ref{thm:local_revisited} that $A(\eps)$ and $B(\eps)$
can be chosen to have $1/2$-summable power series in $\eps$. These power series have
sums that are analytic in the sector $V_{\rho,\delta}$ with
continuous limit at $\eps=0$. When they are not identically zero,
they have the form $\eps^Nc(\eps)$ with $c(\eps)$ nonzero, analytic
in the sector with continuous nonzero limit at $\eps=0$. Dividing
$A$ by such a function (and multiplying $B$ by the same amount) is allowed
in the equivalence class for the modulus. Thus, in the case when $A\neq 0$
we can take a scaling so that $A\equiv \eps^{N_A}$, for some $N_A\in \N$.
This gives cases (i) or (iii).  In the case where $A\equiv 0$ we can perform
a similar division on $B$ to give (ii) or (iv).  It is clear that no more
scalings are allowed within the equivalence classes, and so the representations
(i) to (iv) are unique.\hfill $\Box$

\subsection{The Riccati case}

Here we use the following notation
$$\begin{cases}
m_A=m_{A,1},\\
T_B=T_{B,1}. \end{cases}$$

\begin{theorem}\label{thm:Riccati2} For any germs of
analytic functions $a(\eps), A(\eps),B(\eps)$, the modulus
$(a(\eps), [m_{A(\eps)}, L_\beta\circ T_{B(\eps)}])$ can be realized
as the modulus of the unfolding of the holonomy of the strong separatrix
of a Riccati equation
\begin{equation}\begin{array}{lll}
\dot{x} &=& x^2-\eps,\\
\dot{y} &=& f_{0,\eps}(x) + yf_{1,\eps}(x)
+y^2f_{2,\eps}(x),\end{array}\label{Riccati2}\end{equation} with
$f_{ j,\eps}$ a germ of analytic family of functions in $x$.
\end{theorem}
\noindent{\scshape Proof.} It is proved in \cite{R2} that the
modulus of the unfolding of the holonomy of such a family is formed
by M\"obius functions, hence by analytic functions $m_{A(\eps)}$,
$T_{B(\eps)}$ as in Theorem~\ref{thm:Riccati}. It is also shown
there that the spherical coordinates (called $w$) on the fundamental
domains of Figure~\ref{modulus} can be obtained by first integrals
of the saddle-node model family
\begin{equation}\begin{array}{lll}
\dot{x} &=& x^2-\eps\\
\dot{y} &=& y(1+ax),
\end{array}\label{sad_node_model}\end{equation}
which is the point of view in \cite{I} and \cite{RT}.  For this
reason, we will be brief with the details. We intend to treat in
full detail the general case of a saddle-node in a forthcoming
paper.

We first discuss the local realization of a family with modulus
$(a(\eps), [m_{A(\eps)},L_\beta\circ T_{B(\eps)}])$, i.e. of
 a ramified (in $\hat{\eps}$) family realizing this modulus. For the local
construction (local in $\hat{\eps}$), we consider the two same
sectors $U_{\hat{\eps}}^\pm$ of Figure~\ref{sectors} and their
intersection which is formed of the three (resp. two) sectors
$U_{\hat{\eps}}^{0, \infty,C}$ (resp. $U_{0}^{0, \infty}$) for
$\hat{\eps}\neq0$ (resp. $\eps=0$). Note that $r$ can be chosen
arbitrarily large since $\psi_{\hat{\eps}}^{0,\infty}$ are global
diffeomorphisms. Let
$$\mathcal{U}_{\hat{\eps}}^{\#} = U_{\hat{\eps}}^{\#}\times \CP^1$$
for $\#\in \{+, -, 0, \infty, C\}$. On each
$\mathcal{U}_{\hat{\eps}}^\pm$ we take the model family
\eqref{sad_node_model} in coordinates $(x,y^\pm)$. We glue together
the two models over $\mathcal{U}_{\hat{\eps}}^{\#}$, $\#\in
\{0,\infty,C\}$. Over $\mathcal{U}_{\hat{\eps}}^\pm$ we have first
integrals $H_{\hat{\eps}}^\pm(x,y^\pm)= y^\pm g_{\hat{\eps}} (x)$
with $g_{\hat{\eps}}$ given in \eqref{eq_g}. We need to write the
change of coordinates over $\mathcal{U}_{\hat{\eps}}^{0,\infty,C}$.
It comes from the change in first integral
\begin{equation}
H_{\hat{\eps}}^-=\begin{cases} \psi_{\hat{\eps}}^0
(H_{\hat{\eps}}^+)= m_{A(\eps)}(H_{\hat{\eps}}^+),
&\text{on} \quad \mathcal{U}_{\hat{\eps}}^0,\\
\psi_{\hat{\eps}}^\infty (H_{\hat{\eps}}^+)=
L_{\beta(\eps)}T_{B(\eps)}(H_{\hat{\eps}}^+), &\text{on} \quad
\mathcal{U}_{\hat{\eps}}^\infty,\\
L_{\ov{C}(\hat{\eps})} (H_{\hat{\eps}}^+), &\text{on} \quad
\mathcal{U}_{\hat{\eps}}^C,\end{cases}\end{equation} and yields
\begin{equation}
(x,y^-)=\begin{cases}
\left(x,\frac{y^+}{1+\frac{A(\eps)}{g_{\hat{\eps}} (x)}y^+}\right),
&\text{on} \quad \mathcal{U}_{\hat{\eps}}^0,\\
\left(x,y^+ + B(\eps) g_{\hat{\eps}} (x)\right), &\text{on} \quad
\mathcal{U}_{\hat{\eps}}^\infty,\\
(x,y^+), &\text{on} \quad
\mathcal{U}_{\hat{\eps}}^C.\end{cases}\label{mapeqn1}
\end{equation} Note that
$g_{\hat{\eps}}(\sqrt{\hat{\eps}})=0$ and
$1/g_{\hat{\eps}}(-\sqrt{\hat{\eps}})=0$, so we can glue in the two
lines $\{\pm\sqrt{\hat{\eps}}\}\times \CP^1$ to obtain a $C^\infty$
manifold. We show that this manifold is analytic. For this it
suffices to see that a cylindrical neighborhood of each line
$\{\pm\sqrt{\hat{\eps}}\}\times \CP^1$ minus the corresponding line
 is analytically isomorphic to the product of a pointed disk with $\CP^1$. Let us now write the details for
  a neighborhood
of the line $\{\sqrt{\hat{\eps}}\}\times \CP^1$. We consider
$\check{U}$ a small disk centered at $\sqrt{\hat{\eps}}$ that does
not contain $-\sqrt{\hat{\eps}}$ and $\check{U}^*$ the pointed disk.
We look for global coordinates $(x,Y)$ on $\check{U}^* \times
\CP^1$. For this, we look for functions $k^\pm(x)$ such that
\begin{equation} Y^\pm = y^\pm+k_{\hat{\eps}}^{\pm}
(x)\end{equation} and $Y^+\equiv Y^-$ over
$\mathcal{U}_{\hat{\eps}}^+\cap \mathcal{U}_{\hat{\eps}}^-$. Then
$k_{\hat{\eps}}^\pm$ must satisfy
\begin{equation}
k_{\hat{\eps}}^+(x)-k_{\hat{\eps}}^-(x)= \begin{cases} 0,&x\in U_{\hat{\eps}}^0\cup U_{\hat{\eps}}^C,\\
B(\eps) g_{\hat{\eps}}(x), &x\in
U_{\hat{\eps}}^\infty.\end{cases}\end{equation} There are just found
as solutions of the Cousin problem. The explicit formula for the
solution allows to show that they have a limit at
$\sqrt{\hat{\eps}}$.  Since $g_{\hat{\eps}}(\sqrt{\hat{\eps}})=0$,
they can be taken such that $k^\pm (\sqrt{\hat{\eps}})=0$. The
global coordinate we are looking for is given by  $Y= Y^\pm (x,y)$
on $\mathcal{U}_{\hat{\eps}}^{\pm}\cap (\check{U}^* \times \CP^1$
with analytic extension to $\check{U}\times \CP^1$.

A similar proof can be done in a neighborhood of the line
$\{-\sqrt{\hat{\eps}}\}\times \CP^1$. It can be reduced to the
previous proof if we use the change $Y\pm\mapsto 1/Y^\pm$.

So the manifold we have constructed is a 2-dimensional complex
analytic manifold which is fibred over a disk with a fiber given by
$\CP^1$.  Since any vector bundle over a noncompact Riemann surface
is holomorphically trivial (see for instance \cite{F}), this bundle
must also be holomorphically trivial since it is clear that it can
be constructed as the projectivization of a vector bundle, using
a suitable lift of the maps \eqref{mapeqn1}.

Of course, we would have
obtained the same result if we had used the Newlander-Nirenberg
theorem. There, we could have included $\hat{\eps}$ as a parameter
and obtained that the construction depends analytically on $\hat{\eps}$.
And it is of course possible to manage that the limit exists for
$\hat{\eps}=0$

\smallskip \noindent{\bf Correction to a uniform family.} The family we have realized is
defined over $B(0,r)\times \CP^1$ for values of $\hat{\eps}$ in a
sector of radius $\rho$ and of opening greater than $2\pi$. For this
correction, we use the Newlander-Nirenberg theorem as in
Section~\ref{sec:global}. Indeed the vector field for $\ov{\eps}\in
\ov{V}$ is conjugate to that for $\tilde{\eps}\in \widetilde{V}$.
Let $(x,\ov{Y},\ov{\eps})\mapsto
(x,\Xi(x,\ov{Y},\ov{\eps}),\tilde{\eps})$ be this conjugating map.
This map can be used to glue the family of vector fields over
$\ov{V}$ with the family of vector fields over $\widetilde{V}$. So
we realize a family of vector fields over a 3-dimensional analytic
manifold $\mathcal{M}$. We glue in $B(0,r)\times \CP^1\times
\{\eps=0\}$, thus obtaining a $C^\infty$-manifold. We must recognize
that this manifold is of the form $V\times B(0,r)\times \CP^1$. For
this we endow it of an integrable almost complex structure. Two of
the forms are given by $dx$ and $d\eps$. A form playing the role of
$dY$ is constructed as in the proof of
Theorem~\ref{thm_global_real}. The variables $x$ and $\eps$ remain
holomorphic in the new coordinates, and give a projection from the
image of the corrected manifold onto a neighborhood of
$(x,\eps)=(0,0)$. The inverse image of each point $(x,\eps)$ close
to $(0,0)$ is clearly isomorphic to the Riemann sphere.
We conclude by applying the Fisher-Grauert theorem to conclude
that the bundle has a local trivialization \cite{FG}. \hfill
$\Box$

 \bigskip

For generic $a(\eps), A(\eps),B(\eps)$, the triple $(a(\eps),
[m_{A(\eps)},L_\beta\circ T_{B(\eps)}])$ can be realized as the
modulus of the unfolding of the holonomy of the strong separatrix of
a Riccati equation given by a quadratic vector field.  Since this
proof is completely elementary, we add it for completeness.

\begin{theorem}\label{thm:Riccati3} Given germs of analytic functions $A(\eps)$ and $B(\eps)$,
then for most $a_0$ and for a corresponding germ of analytic
function $a(\eps)$ yielding a realizable family of triple $(a(\eps),
[m_{A(\eps)}, L_\beta\circ T_{B(\eps)}])$ as in Theorem~\ref{thm:Riccati}, there
exists analytic functions $c(\eps)$ and $d(\eps)$ such that the
triple can be realized as the moduli of the unfolding of the
holonomy of the strong separatrix of a Riccati equation of the form
\begin{equation}\begin{array}{lll}
\dot{x} &=& x^2-\eps\\
\dot{y} &=& c(\eps)(x^2-\eps) + y(1+a(\eps)x) + d(\eps)y^2.
\end{array}\label{Riccati}\end{equation}
There is no restriction on $a(\eps)$ when $A(0)B(0)\neq0$. Also,
when $A(0)=B(0)=0$ and $a(0)$ is not an integer, then the
triple $(a(\eps), [m_{A(\eps)}, L_\beta\circ T_{B(\eps)}])$ can be realized.
\end{theorem}

\noindent {\scshape Proof.}
For the system
\begin{equation}\begin{array}{lll}
\dot{x} &=& x^2-\eps\\
\dot{y} &=& \a(\eps)\b(\eps)(x^2-\eps) + y(1+(1-\a(\eps)-\b(\eps))x)
+ y^2,
\end{array}\label{Riccati3}
\end{equation}
it is shown in \cite{LR}, that the moduli are given (up to a scaling
of the form $(A,B)\to (Ak,B/k)$ with $k$ bounded and bounded away
from $0$) by
$$
    A(\eps) = \frac{2\pi i}{\G(1-\a)\G(1-\b)},\qquad  B(\eps) = \frac{- 2\pi i e^{\pi i (1-\a-\b)}}{\G(\a)\G(\b)}.
$$

We first take $d(\eps)=1$ and $c(\eps)=\a(\eps)\b(\eps)$ in \eqref{Riccati} to obtain \eqref{Riccati3}
where $a(\eps)=1-\a(\eps)-\b(\eps)$.
Thus,
$$
    A(\eps) = \frac{2\pi i}{\G(a+\b)\G(1-\b)},\qquad  B(\eps) = \frac{- 2\pi i e^{\pi i a}}{\G(1-a-\b)\G(\b)}.
$$
If $a(0)$ is not an integer, it is clear that we can choose $\a$ and
$\b$ to obtain any values of the parameters we wish (making sure
that we have $\b(0)\neq 0$), except for the cases where $A$ and $B$
both have a zero at $\eps = 0$. (Recall, that $A$ and $B$ are only
defined up to an inessential scaling.) If $a(0)$ is an integer, we
can only realize $A(\eps)$ and $B(\eps)$ when $A(0)B(0)\neq 0$.

\medskip
To discuss now the cases $A(0)=B(0)=0$,  we consider \eqref{Riccati}
with $d(0) = 0$ but $d(\eps)\not\equiv 0$. For $0 \neq \eps << 1$ we
can substitute $d(\eps)y\mapsto y$ to obtain
\begin{equation}\begin{array}{lll}
\dot{x} &=& x^2-\eps\\
\dot{y} &=& \c(\eps)d(\eps)(x^2-\eps) + y(1+a(\eps)x) + y^2,
\end{array}\label{Riccati4}
\end{equation}
and take $\c(\eps)=\a(\eps)\ov\b(\eps)$, and denote $\b(\eps)=\ov\b(\eps)d(\eps)$,
to obtain \eqref{Riccati3} where $a = 1-\a-\b$ as before.

However, this calculation is only for $\eps\neq0$ and we need to make sure that
the scaling factor is correct in the limit as $\eps$ tends to $0$.

The values of $A$ and $B$ in \cite{LR} are obtained from the first integral, $H$ say,
of \eqref{Riccati3} which is of the form
$$
    \kappa\, \frac{w_2 y + (x^2-\eps) w_2'}{w_3 y + (x^2-\eps)w_3'},
$$
where $\kappa = (2\sqrt{\eps})^{1-\a-\b}e^{\pi i (\frac{a+b-1}{2}+\frac{1}{2\sqrt{\eps}})}$,
and $w_2$ and $w_3$ are given by hypergeometric functions, and
in particular,
$$
  w_3 = _2F_1\left(\a,\b,\frac{1+\a+\b}{2}-\frac1{2\sqrt{\eps}},1-\frac{x}{\sqrt{\eps}}\right).
$$
In our case, we have $\b = \ov\b d$, and hence $d$ divides each term in $w_3'$.  Thus,
in original coordinates, we need to replace $H$ by
$$
  \ov{H} = H d = \kappa\, \frac{w_2 d\,Y + (x^2-\eps) w_2'}{w_3 Y + (x^2-\eps)w_3'/d},
$$
to achieve a uniform limit as $\eps$ tends to zero.
This means a scaling of $d(\eps)$ in the modulus given in \cite{LR}, which gives
$$
    A(\eps) = \frac{d}{\G(a+\ov\b d)\G(1-\ov\b d)},
    \qquad  B(\eps) = \frac{- 2\pi i e^{\pi i (1-a-b)}}{\G(1-a-\ov\b d)\G(\ov\b d)d}.
$$
We note that $(\G(\ov\b d)d)^{-1} = \ov\b + o(\ov\b,d)$, and hence,
if $a(0)$ is not an integer, we can clearly choose $\ov\b$
and $d$ to obtain any germs of functions $A$ and $B$ with
$A(0)=B(0)=0$. \hfill $\Box$

\begin{remark} The triple $(a(\eps), [m_{A(\eps)},L_\beta\circ T_{B(\eps)}])$
cannot be realized in a family of type \eqref{Riccati3}, when
$a(0)=2$, $A(0)=0$ and $B(\eps)\neq0$. \end{remark}

\subsection{The only families with continuous representative
$\psi_{\eps}^{0,\infty}$ of the modulus}

We propose the following conjecture which we prove in a special
case.

\begin{conjecture}\label{conjecture}
The only families with representative $\psi_{\eps}^{0,\infty}$ of
the modulus which are analytic in $\eps$ are the ones presented in
Theorem~\ref{Riccati}.\end{conjecture}

\begin{theorem}\label{thm:conjecture} The conjecture~\ref{conjecture} is valid in the
subcase where either $\psi_{\eps}^\infty$ (or $\psi^0$) is linear.
\end{theorem}

This has been proved in the case $a=0$ by Reinhard Sch\"afke
\cite{Sc}.\smallskip

 \noindent {\scshape Proof of Theorem~\ref{thm:conjecture}.}
 We make the proof in the case
where $\psi_{\eps}^\infty$ is linear, and thus $\ov{h}^\infty=id$ and $\tilde{h}^\infty=id$.
Using the
notation of Section~\ref{sect:Riccati}, the compatibility condition
is given by
\begin{equation}(\tilde{h}^0)^{-1}= L_{F}\circ
\ov{h}^0\circ L_{G},\label{compatibility_mobius0}\end{equation}
where
\begin{equation}\begin{cases} \tilde{h}^0\circ
L_{(\ov{C}\beta)^{-1}}\circ
\widetilde{\psi}^0=L_{(\ov{C}\beta)^{-1}}\circ\tilde{h}^0,
\\
\ov{h}^0\circ \ov{\psi}^0\circ
L_{\ov{C}}=L_{\ov{C}}\circ\ov{h}^0.\end{cases}\label{eq0_h}\end{equation}
We note that $G=F^{-1}$ because $(\ov{h}^0)'(0)=(\tilde{h}^0)'(0)=1$ in \eqref{compatibility_mobius0}.

If $\psi_{\eps}^0$ depends analytically on $\eps$, then
$\ov{\psi}^0=\widetilde{\psi}^0$.
From \eqref{eq0_h}, we have
\begin{equation}\begin{cases}
\widetilde{\psi}^0= L_{\ov{C}\beta}\circ(\tilde{h}^0)^{-1}\circ L_{(\ov{C}\beta)^{-1}}\circ\tilde{h}^0,\\
\ov{\psi}^0=(\ov{h}^0)^{-1}\circ L_{\ov{C}}\circ\ov{h}^0\circ
L_{(\ov{C})^{-1}}.\end{cases}\end{equation} Since
$\ov{\psi}^0=\widetilde{\psi}^0$, this yields, after some rearrangement,
\begin{equation}
\ov{h}^0\circ L_{(\ov{C})^{-1}}\circ (\tilde{h}^0)^{-1}\circ
L_{\ov{C}\beta} =L_{(\ov{C})^{-1}}\circ\ov{h}^0\circ
L_{\ov{C}\beta}\circ
(\tilde{h}^0)^{-1}.\label{EQUATION}\end{equation}
Substituting \eqref{compatibility_mobius0} yields
\begin{equation}
\ov{h}^0\circ L_{F(\ov{C})^{-1}}\circ \ov{h}^0\circ L_{\ov{C}\beta}
=L_{(\ov{C})^{-1}}\circ\ov{h}^0\circ L_{F\ov{C}\beta}\circ
\ov{h}^0.\label{EQUATION2}\end{equation}

We now substitute $\ov{h}^0(w)=w+\sum_{j\geq2}b_jw^j$ and equate
coefficients of $w^j$ in \eqref{EQUATION2}. Let $b_s$ be the
first nonzero coefficient. Then we need to choose
$F^{s-1}=\frac{(\ov{C}\beta)^{1-s} -1}{(\ov{C})^{1-s}-1}$. We note that
$F\beta=1$ if and only if $a \in \frac1{2\pi i}\mathbb Z$ (i.e. $\beta=1$).
For $j>s$, the coefficient of $w^j$ is a polynomial in $b_2, \dots,
b_j$, where the only monomial in $b_j$  is of the form $c_jb_j$ with
$$c_j=F\beta
\left[1-(F\beta)^{j-1} + (\ov{C}\beta)^{j-1}(F^{j-1}-1)\right]$$
which does not vanish for $|\eps|<<1$ as soon as $F\beta\neq 1$. This means that
the solution is unique. Since $m_{(1-s)b_s,s-1}$ is one solution, it is
the only one.

We now need to treat the case $a \in \frac1{2\pi i}\mathbb Z$. In
this case $\beta = 1$ and \eqref{EQUATION2} gives
\begin{equation}
\ov{h}^0\circ L_{(\ov{C})^{-1}}\circ \ov{h}^0\circ L_{\ov{C}}
-L_{(\ov{C})^{-1}}\circ\ov{h}^0\circ L_{\ov{C}}\circ
\ov{h}^0=0.\label{EQUATION3}\end{equation} If $b_2\neq 0$, the terms
of degree $2$ and $3$ yield no constraints, but the terms of
degree 4 give $b_3-b_2^2=0$, and the terms of degree $j$ give
$b_{j-1}$ uniquely in terms of $b_2,\ldots,b_{j-2}$.  Since
$\ov{h}^0=m_{-b_2,1}$ is a solution, it must be the solution.

If $b_2$ vanishes,
then we take $b_s$ to be the first non-vanishing coefficient of $\ov{h}^0-id$ as above.
Suppose $b_j \neq 0$ is the first term non-zero coefficient for which $j-1$ is not
divisible by $s-1$.  We consider the term in $w^{j+s-1}$ in \eqref{EQUATION3} to
yield a contradiction.  Thus, all terms in $\ov{h}^0$ only contains terms in
$b_{1+k(s-1)}x^{1+k(s-1)}$.  We conjugate \eqref{EQUATION3} by $z\to w^{s-1}$ to
obtain an equation of the form \eqref{EQUATION3} with $\ov{C}$ replaced by
$\ov{C}^{s-1}$ and $\ov{h}$ replaced by a power series $\ov{h}'$ in $z$ with the coefficient
of $z^2$ given by $(s-1)b_s$.  Therefore, we can proceed as above to obtain $\ov{h}' = m_{(1-s)b_s,1}(z)$,
and hence $\ov{h_0}=M_{(1-s)b_s,s-1}$.\hfill$\Box$

\section*{Acknowledgements} We are grateful to Reinhard Sch\"afke and
Lo\"ic Teyssier for helpful discussions.

\end{document}